\documentclass[11pt,a4paper,usenames,dvipsnames]{article}

\usepackage{url,amsmath,amssymb,latexsym,mathrsfs,comment,amsthm,enumerate,geometry,calc,ifthen,mathdots,textcomp,multicol}

\usepackage{cases}

\usepackage[noadjust]{cite}

\usepackage[colorlinks]{hyperref}

\usepackage[margin=10pt,font=small,labelfont=bf, labelsep=period]{caption}

\usepackage{makecell}

\usepackage[all,cmtip,2cell]{xy}

\usepackage[x11names, svgnames, rgb,table]{xcolor}
\usepackage{hhline}
\usepackage{multirow}

\usepackage{pst-node}
\usepackage{tikz-cd}

\usepackage{enumitem}

\usepackage[utf8]{inputenc}
\usepackage[T1]{fontenc}

\usepackage{tikz}
\newcommand\nc\newcommand
\nc\rnc\renewcommand

\nc\fd{\textup{fd}}
\nc\tc{\textup{tc}}
\nc\Pnfd{\P_n^{\fd}}
\nc\Pntc{\P_n^{\tc}}
\nc\PPnfd{\PP_n^{\fd}}
\nc\Tf{T^\flat}
\nc\D{\mathcal D}
\nc\xra\mr
\nc\tb{\ol t}
\nc\eb{\ol e}
\rnc\sb{\ol s}
\nc\J{\mathcal J}
\nc\F{\mathcal F}
\nc\XSn{X_{\S_n}}
\nc\XFn{X_{\F_n}}
\nc\XTn{X_{\T_n}}
\nc\XDn{X_{\D_n}}
\nc\XOn{X_{\O_n}}
\nc\RSn{R_{\S_n}}
\nc\RFn{R_{\F_n}}
\nc\RTn{R_{\T_n}}
\nc\RDn{R_{\D_n}}
\nc\ROn{R_{\O_n}}
\nc\oijn{1\leq i<j\leq n}
\nc\wh\widehat
\nc\Psim[1]{\sim_\eqref{P#1}}
\nc\Nsim[1]{\sim_\eqref{N#1}}
\nc\sP{\mathscr P}
\nc\suc{\operatorname{suc}}
\nc\sgap{\;\!}

\nc\lc{^\sharp_L}
\nc\rc{^\sharp_R}

\nc\Sgp{\operatorname{\mathsf{Sgp}}}
\nc\Mon{\operatorname{\mathsf{Mon}}}

\nc\trans[1]{\left(\begin{smallmatrix}#1\end{smallmatrix}\right)}

\nc\FIG{\operatorname{\mathsf{IG}}}
\nc\FPG{\operatorname{\mathsf{PG}}}
\nc\bbB{\mathbb B}
\nc\N{\mathbb N}
\nc\Z{\mathbb Z}
\nc\TL{\mathcal T\!\mathcal L}
\rnc\S{\mathcal S}
\nc\Eq{\mathfrak{Eq}}
\nc\PEq{\mathscr P\mathfrak{Eq}}

\newcommand{\uv}[1]{\fill (#1,2)circle(.17);}
\newcommand{\lv}[1]{\fill (#1,0)circle(.17);}

\newcommand{\uvs}[1]{{\foreach \x in {#1} { \uv{\x}}}}
\newcommand{\lvs}[1]{{\foreach \x in {#1} { \lv{\x}}}}
\newcommand{\darcx}[3]{\draw(#1,0)arc(180:90:#3) (#1+#3,#3)--(#2-#3,#3) (#2-#3,#3) arc(90:0:#3);}
\newcommand{\darc}[2]{\darcx{#1}{#2}{.4}}
\newcommand{\uarcx}[3]{\draw(#1,2)arc(180:270:#3) (#1+#3,2-#3)--(#2-#3,2-#3) (#2-#3,2-#3) arc(270:360:#3);}
\newcommand{\uarc}[2]{\uarcx{#1}{#2}{.4}}
\newcommand{\stline}[2]{\draw(#1,2)--(#2,0);}
\newcommand{\stlinex}[3]{\draw[#3](#1,2)--(#2,0);}
\nc{\uarcs}[1]{{\foreach \x/\y in {#1}{ \uarc{\x}{\y} }}}
\nc{\darcs}[1]{{\foreach \x/\y in {#1}{ \darc{\x}{\y} }}}
\newcommand{\stlines}[1]{{\foreach \x/\y in {#1} { \stline{\x}{\y} }}}
\nc\udotted[2]{\draw[dotted](#1+.5,2)--(#2-.5,2);}
\nc\ldotted[2]{\draw[dotted](#1+.5,0)--(#2-.5,0);}
\nc{\udotteds}[1]{{\foreach \x/\y in {#1}{ \udotted{\x}{\y} }}}
\nc{\ldotteds}[1]{{\foreach \x/\y in {#1}{ \ldotted{\x}{\y} }}}
\nc\vertlab[2]{\node()at(#1,2.7){\tiny$#2$};}
\nc\uuline[2]{\draw(#1,2)--(#2,2);}
\nc\ddline[2]{\draw(#1,0)--(#2,0);}
\nc\uulinex[3]{\draw[#3](#1,2)--(#2,2);}
\nc\ddlinex[3]{\draw[#3](#1,0)--(#2,0);}

\nc\hlines[2]{\stline{#1}{#1}\stline{#2}{#2}\uuline{#1}{#2}}
\nc\hlinesx[3]{\stlinex{#1}{#1}{#3}\stlinex{#2}{#2}{#3}\uulinex{#1}{#2}{#3}}
\nc\vlines[2]{\stline{#1}{#1} \stline{#1+.5}{#1+.5} \stline{#2}{#2} \stline{#2-.5}{#2-.5} \udotted{#1+.2}{#2-.2}}

\newcommand{\darcxcol}[4]{\draw[#4](#1,0)arc(180:90:#3) (#1+#3,#3)--(#2-#3,#3) (#2-#3,#3) arc(90:0:#3);}
\newcommand{\darccol}[3]{\darcxcol{#1}{#2}{.4}{#3}}
\nc{\darccols}[2]{{\foreach \x/\y in {#1}{ \darccol{\x}{\y}{#2} }}}
\newcommand{\uarcxcol}[4]{\draw[#4](#1,2)arc(180:270:#3) (#1+#3,2-#3)--(#2-#3,2-#3) (#2-#3,2-#3) arc(270:360:#3);}
\newcommand{\uarccol}[3]{\uarcxcol{#1}{#2}{.4}{#3}}
\nc{\uarccols}[2]{{\foreach \x/\y in {#1}{ \uarccol{\x}{\y}{#2} }}}

\nc{\uvert}[1]{\fill (#1,2)circle(.2);}
\rnc{\lvert}[1]{\fill (#1,0)circle(.2);}
\nc{\custpartn}[3]{{\lower1.4 ex\hbox{
\begin{tikzpicture}[scale=.3]
\foreach \x in {#1}
{ \uvert{\x}  }
\foreach \x in {#2}
{ \lvert{\x}  }
#3 \end{tikzpicture}
}}}

\newcounter{ncols}
\newcounter{incols}
\newenvironment{partn}[1]{
  \setcounter{ncols}{#1} \setcounter{incols}{\thencols - 1}\setlength{\arraycolsep}{1pt}
  \Bigl( \hspace{-1.5truemm}\scriptsize 
    \begin{array}{@{\hskip 3pt} c *{\theincols}{|c} @{\hskip 3pt}  }
}{
     \end{array}
     \normalsize \hspace{-1.5truemm}\Bigr)\setlength{\arraycolsep}{6pt}
}

\nc\la\langle
\nc\ra\rangle

\nc\ol\overline
\nc\ul\underline
\nc\da\downarrow
\nc\pr\bullet
\nc\ccirc{\circ'}
\nc\bc\odot
\nc\sm\setminus

\nc\ds\displaystyle
\nc\ts\textstyle

\nc\bn{{\bf n}}
\nc\bk{{\bf k}}

\nc\al\alpha
\nc\be\beta
\nc\ga\gamma
\nc\de\delta
\nc\ve\varepsilon
\nc\lam\lambda
\nc\om\omega
\nc\ze\zeta
\nc\si\sigma
\nc\Si\Sigma
\nc\Ga\Gamma
\nc\De\Delta
\nc\Om\Omega
\nc\io\iota

\nc\nab\nabla

\nc\vt\vartheta
\rnc\th\theta
\nc\Th\Theta

\nc\BY{\qquad\text{by}\qquad}
\nc\GIVENBY{\qquad\text{given by}\qquad}
\nc\AND{\qquad\text{and}\qquad}
\nc\ANDSIM{\qquad\text{and similarly}\qquad}
\nc\ANDSO{\qquad\text{and so}\qquad}
\nc\ANd{\quad\text{and}\quad}
\nc\COMMA{,\qquad}
\nc\COMMa{,\quad}
\nc\WHERE{\qquad\text{where}\qquad}

\rnc\iff{\ \Leftrightarrow\ }
\nc\IFf{\quad \Leftrightarrow\quad }
\nc\Iff{\ \ \Leftrightarrow\ \ }
\nc\IFF{\qquad \Leftrightarrow\qquad }
\rnc\implies{\ \Rightarrow\ }
\nc\IMPLIES{\qquad \Rightarrow\qquad }

\nc\set[2]{\{#1:#2\}}
\nc\bigset[2]{\big\{#1:#2\big\}}
\nc\pres[2]{\la#1:#2\ra}

\nc\bit{\begin{itemize}[label=\textbullet, leftmargin=5mm]}
\nc\eit{\end{itemize}}
\nc\ben{\begin{enumerate}[label=\textup{(\roman*)},leftmargin=10mm]}
\nc\bena{\begin{enumerate}[label=\textup{(\alph*)},leftmargin=10mm]}
\nc\een{\end{enumerate}}
\nc\bmc{\begin{multicols}}
\nc\emc{\end{multicols}}

\nc\pf{\begin{proof}}
\nc\epf{\end{proof}}
\nc\pfclaim{\begin{quote}\begin{proof}}
\nc\epfclaim{\end{proof}\end{quote}}
\nc\epfres{\hfill\qed}
\nc\epfreseq{\tag*{\qed}}
\let\oldproofname=\proofname
\renewcommand{\proofname}{\rm\bf{\oldproofname}}
\nc{\pfitem}[1]{\medskip \noindent #1.}
\nc{\firstpfitem}[1]{#1.}
\nc{\pfcase}[1]{\medskip\noindent {\bf Case #1.}}
\nc\aftercases{\medskip\noindent}

\renewcommand{\H}{\mathrel{\mathscr H}}
\renewcommand{\L}{\mathrel{\mathscr L}}
\newcommand{\R}{\mathrel{\mathscr R}}
\newcommand{\K}{\mathbb K}
\nc\rH{\mathrel{\H}}
\nc\rL{\mathrel{\L}}
\nc\rR{\mathrel{\R}}
\nc\rD{\mathrel{\D}}
\nc\gJ{\mathrel{\mathscr J}}
\nc\rK{\mathrel{\K}}
\nc\rsi{\mathrel{\si}}

\nc\leqL{\leq_{\L}}
\nc\leqR{\leq_{\R}}
\nc\leqJ{\leq_{\gJ}}
\nc\geqL{\geq_{\L}}
\nc\geqR{\geq_{\R}}
\nc\geqJ{\geq_{\gJ}}
\nc\leqH{\leq_{\H}}
\nc\leqF{\leq_{\F}}
\nc\geqF{\geq_{\F}}

\newcommand{\Sing}{\operatorname{Sing}}

\newcommand{\id}{\operatorname{id}}
\newcommand{\dom}{\operatorname{dom}}
\newcommand{\codom}{\operatorname{codom}}
\newcommand{\coker}{\operatorname{coker}}
\newcommand{\rank}{\operatorname{rank}}

\nc\mr\mathrel
\nc\pc[2]{(#1,#2)^\sharp}

\nc\U{\mathcal U}
\nc\V{\mathcal V}
\nc\G{\mathcal G}
\rnc\iff{\ \Leftrightarrow\ }
\rnc\implies{\ \Rightarrow\ }
\nc\Implies{\quad \Rightarrow\quad }
\nc\C{\mathscr C}
\nc\M{\mathcal M}
\nc\CC{\mathcal C}
\nc\DD{\mathcal D}
\nc\FF{\mathcal F}
\nc\I{\mathcal I}
\rnc\O{\mathcal O}
\rnc\P{\mathcal P}
\nc\PP{\mathscr P\mathcal P}
\nc\PE{\mathscr P\mathcal E}
\nc\PT{\mathscr P\mathcal T}
\nc\T{\mathcal T}
\nc\p{\mathfrak p}
\nc\q{\mathfrak q}
\rnc\r{\mathfrak r}
\nc\s{\mathfrak s}
\rnc\t{\mathfrak t}
\nc\bd{{\bf d}}
\nc\br{{\bf r}}
\nc\op\oplus
\nc\lra{\mathrel\leftrightarrow}
\nc\es\varnothing
\nc\rev{\textup{rev}}
\nc\corestt{{\upharpoonleft}}
\nc\restt{{\upharpoonright}}
\nc\corest{{\downharpoonleft}}
\nc\rest{{\upharpoonright}}
\nc\WHERe{\quad\text{where}\quad}
\nc\ot\leftarrow
\rnc\a{\mathfrak a}
\rnc\b{\mathfrak b}
\rnc\c{\mathfrak c}
\rnc\d{\mathfrak d}
\nc\sub\subseteq
\nc\mt\mapsto
\nc\im{\operatorname{im}}
\nc\B{\mathcal B}
\nc\E{\mathcal E}

\numberwithin{equation}{section}

\newtheorem{thm}[equation]{Theorem}
\newtheorem{lemma}[equation]{Lemma}
\newtheorem{cor}[equation]{Corollary}
\newtheorem{prop}[equation]{Proposition}

\theoremstyle{definition}

\newtheorem{defn}[equation]{Definition}
\newtheorem{rem}[equation]{Remark}

\newcounter{caseco}

\newcounter{subcaseco}

\newcounter{stepco}

\newcounter{stageco}

\begin{document}

\title{\vspace{-1.4cm}Presentations for semigroups of full-domain partitions\vspace{-0.2cm}}
\date{}

\author{Luka Carroll\footnote{The first author was supported by the Australian Mathematical Sciences Institute (AMSI).},\ James East, Matthias Fresacher\\[3mm]
{\it\small Centre for Research in Mathematics and Data Science,}\\
{\it\small Western Sydney University, Locked Bag 1797, Penrith NSW 2751, Australia.}\\[3mm]
{\tt\small 17453699@student.WesternSydney.edu.au}, {\tt\small J.East@WesternSydney.edu.au}, \\ {\tt\small M.Fresacher@WesternSydney.edu.au}}

\maketitle

\vspace{-1.0cm}

\begin{abstract}
\noindent
The full-domain partition monoid $\Pnfd$ has been discovered independently in two recent studies on connections between diagram monoids and category theory.  It is a right restriction Ehresmann monoid, and contains both the full transformation monoid and the join semilattice of equivalence relations.  In this paper we give presentations (by generators and relations) for $\Pnfd$, its singular ideal, and its planar submonoid.  The latter is not an Ehresmann submonoid, but it is a so-called grrac monoid in the terminology of Branco, Gomes and Gould.  In particular, its structure is determined in part by a right regular band in one-one correspondence with planar equivalences.

\medskip

\medskip


\end{abstract}

\tableofcontents

\section{Introduction}\label{sect:intro}

Diagram algebras are defined in terms of set partitions that come with a natural graphical representation and multiplication.  Originally arising in representation theory and physics, they have become important tools in many parts of mathematics and science, including topology, invariant theory, combinatorics, logic, category theory, computer science, statistical mechanics and biology; see for example the recent surveys \cite{Martin2008,Koenig2008}.  Diagram algebras are all constructed as twisted semigroup algebras of associated diagram monoids, such as partition, (partial) Brauer, Temperley--Lieb and Motzkin monoids~\cite{Martin1994,Jones1994_2, Brauer1937,MM2014, TL1971, BH2014, Wilcox2007}.  Some key applications of this fact include proofs of cellularity \cite{Wilcox2007} and presentations by generators and relations \cite{JEgrpm,JEgrpm2,JErook,East2021}.  But diagram algebras have also `given back' to semigroup theory; for some recent articles covering topics such as (pseudo)varieties, identities, Krohn--Rhodes complexity, congruence lattices, endomorphism monoids and combinatorial invariant theory; see for example \cite{EMRT2018,EG2017,EG2021,MM2007,KM2006,Maz2002,LF2006,FL2011,BDP2002, KV2019,KV2023,CHKLV2019, ACHLV2015,ADV2012_2,Auinger2014,Auinger2012,CHKLV2019}.

The current article concerns the \emph{full-domain partition monoids} $\Pnfd$.  These were recently discovered independently in studies of Ehresmann structures on partition monoids \cite{EG2021}, and categorical completions of transformation monoids \cite{Stokes2022b}.  Both papers \cite{EG2021,Stokes2022b} identified natural \emph{right restriction} structures \cite{Gould_notes} on $\Pnfd$, arising from suitable domain/range maps.  As noted in~\cite{EG2021},~$\Pnfd$ can be thought of as a categorical dual to the partial transformation monoid $\P\T_n$.  The latter is a \emph{left} restriction monoid that extends the full transformation monoid $\T_n$ by the meet semilattice of subsets of $\{1,\ldots,n\}$, whereas the right restriction monoid $\Pnfd$ extends $\T_n$ by the join semilattice of equivalence relations.  

The monoid $\Pnfd$ turns out to have a very interesting structure.  Certain aspects were investigated in \cite{EG2021}, including its regularity and ideal structure, a characterisation of Green's divisibility relations, links to Ehresmann categories (in the sense of Lawson \cite{Lawson1991}), and applications to representation theory.  Its free idempotent-generated semigroup, and an associated topological/combinatorial complex, is studied in \cite{EGM2025}.

In this paper we continue the study of $\Pnfd$, moving in the direction of presentations by generators and relations.  Presentations are key tools in many studies of (diagram) monoids and algebras, as they allow for a concise combinatorial summary of the structures, and also allow for the construction of morphisms/representations and actions from minimal data.  En route to obtaining a presentation for $\Pnfd$ itself, we first obtain one for the singular ideal ${\Sing(\Pnfd) = \Pnfd\sm\S_n}$.  Here~$\S_n$ is the symmetric group, which is the group of units of $\Pnfd$.  We also give a presentation for the submonoid $\PPnfd$ consisting of all \emph{planar} full-domain partitions, i.e.~those that can be drawn with no edge crossings.  The planar submonoid has a more complicated structure.  It is not closed under the Ehresmann operations of $\Pnfd$, and is not a right restriction monoid, but it still has a natural structure as a unary semigroup, namely a \emph{grrac monoid} in the sense of \cite{BGG2010}; rather than a semilattice, the range of the unary operation is a right regular band.  

The article is organised as follows.  After giving some necessary preliminaries in Section \ref{sect:prelim}, we state our main results in Section \ref{sect:main}.  These are Theorems \ref{thm:SingPnfd}, \ref{thm:Pnfd} and \ref{thm:PPnfd}, which respectively give presentations for $\Sing(\Pnfd)$, $\Pnfd$ and $\PPnfd$.  These theorems are proved using techniques developed in \cite{CDEGZ2023} for working with semigroup products arising from so-called \emph{action pairs}, which are reviewed in Section \ref{sect:AP}.  Next, in Section \ref{sect:Ehresmann} we cover the required definitions and results from general Ehresmann theory \cite{Lawson1991}, and the Ehresmann structure on the partition monoid \cite{EG2021}, ultimately identifying natural action pairs $(\E_n,\T_n)$ and $(\E_n,\Sing(\T_n))$ in $\Pnfd$ and $\Sing(\Pnfd)$, where~$\E_n$ is the join semilattice of equivalence relations.  These are then used to prove Theorems~\ref{thm:SingPnfd} and~\ref{thm:Pnfd} in Sections \ref{sect:SingPnfd} and \ref{sect:Pnfd}, respectively.  Finally, Section \ref{sect:PPnfd} gives the proof of Theorem~\ref{thm:PPnfd}.  This is far more intricate than the proofs of the other two theorems.  It still involves an action pair,~$(\D_n,\O_n)$, where here $\O_n$ is the monoid of \emph{order-preserving} transformations, and $\D_n$ is an apparently new right regular band in one-one correspondence with the planar equivalences.  A key step is to obtain a presentation for $\D_n$ itself, which we do in Theorem~\ref{thm:Dn}.  We conclude with Proposition~\ref{prop:grrac}, which establishes the above-mentioned grrac structure of~$\PPnfd$.

\section{Preliminaries}\label{sect:prelim}

In this section we give the preliminary definitions on semigroups (Section \ref{subsect:S}) and partition monoids (Section \ref{subsect:Pn}) that will enable us to state our main results in Section \ref{sect:main}.  For more background, see for example \cite{Howie1995,CPbook,EMRT2018}.

\subsection{Semigroups and presentations}\label{subsect:S}

A \emph{semigroup} is a set with an associative binary operation, usually denoted by juxtaposition.  A \emph{monoid} is a semigroup with a (two-sided) identity element.  For a semigroup $S$, we denote by~$S^1$ the \emph{monoid completion} of $S$; if $S$ is a monoid then $S^1=S$, or otherwise $S^1=S\cup\{1\}$, where~$1$ is a symbol not belonging to $S$, acting as an adjoined identity.  

An \emph{idempotent} of a semigroup is an element $a$ satisfying $a^2=a$.  A \emph{band} is a semigroup consisting entirely of idempotents.  A \emph{semilattice} is a commutative band.  Any semilattice is partially ordered by $a\leq b \iff a = ab(=ba)$; the product of any two elements is their greatest lower bound (meet).  We will also encounter \emph{right regular bands}, which are defined by the identities $x^2=x$ and $xyx=yx$.  Clearly semilattices are right regular bands.

\emph{Green's pre-orders} on a semigroup $S$ are defined by
\[
x\leqL y \iff x\in S^1y \COMMA x\leqR y \iff x\in yS^1 \AND x\leqJ y \iff x\in S^1yS^1.
\]
These induce \emph{Green's equivalences}
\[
{\L} = {\leqL}\cap{\geqL} \COMMA {\R} = {\leqR}\cap{\geqR} \AND {\gJ} = {\leqJ}\cap{\geqJ}.
\]

An equivalence relation $\si$ on a semigroup $S$ is a \emph{left congruence} if it is \emph{left-compatible}, in the sense that $(a,b)\in\si \implies (sa,sb)\in\si$ for all $a,b,s\in S$.  For example, $\R$ is a left congruence.  For a set of pairs $\Om\sub S\times S$, we write $\Om\lc$ for the left congruence on $S$ generated by $\Om$, i.e.~the intersection of all left congruences containing $\Om$.   In the case that $\Om$ contains a single pair $(a,b)$, we write $(a,b)\lc = \Om\lc$, and call such a left congruence \emph{principal}.  
\emph{Right congruences} are defined dually.

A \emph{(two-sided) congruence} is a left- and right-compatible equivalence $\si$.  The \emph{quotient semigroup} $S/\si$ consists of all $\si$-classes, under the induced operation.  This is well-defined because of compatibility, which is easily seen to be equivalent to having $(a,b),(s,t)\in\si \implies (as,bt)\in\si$ for all $a,b,s,t\in S$.  
For $\Om\sub S\times S$, we write $\Om^\sharp$ for the (two-sided) congruence generated by $\Om$.  

The \emph{kernel} of a semigroup morphism $\phi:S\to T$ is the congruence ${\ker(\phi) = \set{(a,b)}{a\phi=b\phi}}$ on $S$.  The Fundamental Homomorphism Theorem for semigroups states that $S/\ker(\phi)$ is isomorphic to $\im(\phi)$, the image of $\phi$.

For a set $X$, the \emph{free semigroup} $X^+$ consists of all words over $X$, under the operation of concatenation.  For $R\sub X^+\times X^+$, we say a semigroup $S$ has \emph{semigroup presentation} ${\Sgp\pres XR}$ if $S \cong X^+/R^\sharp$.  This is equivalent to there existing a surmorphism (surjective morphism)~${\phi:X^+\to S}$ with kernel $R^\sharp$, in which case we say $S$ has presentation $\Sgp\pres XR$ \emph{via~$\phi$}.  Elements of $X$ and $R$ are typically called \emph{generators} and \emph{relations}, and a relation $(u,v)$ is often displayed as an equality: $u=v$.  A \emph{normal form function} associated to a surmorphism~${\phi:X^+\to S}$ is a function $N:S\to X^+$ such that $N(s)\phi = s$ for all $s\in S$.

The \emph{free monoid} $X^*$ is defined to be $X^+\cup\{\io\}$, where $\io$ denotes the \emph{empty word}, the identity of~$X^*$.  We can also speak of monoid presentations $\Mon\pres XR$, with the appropriate adaptations.

Certain well known families of semigroups will play an important role in the paper.  For an integer $n\geq0$ we write $\bn=\{1,\ldots,n\}$, interpreting this to be empty of $n=0$.  The \emph{full transformation monoid} $\T_n$ is the set of all mappings $\bn\to\bn$, under composition.  The group of units of $\T_n$ is the \emph{symmetric group} $\S_n$, which consists of all permutations of $\bn$.  A transformation $f\in\T_n$ is \emph{order-preserving} if $x\leq y \implies xf\leq yf$ for all $x,y\in\bn$.  The set $\O_n$ of all order-preserving transformations is a submonoid of $\T_n$.

We will also encounter the join semilattice $\Eq_n$ of all equivalence relations on $\bn$.  The \emph{join} $\ve\vee\eta$ of two equivalences $\ve,\eta\in\Eq_n$ is the least equivalence containing $\ve\cup\eta$.  The identity element of $\Eq_n$ is the trivial relation $\De_{\bn} = \set{(x,x)}{x\in\bn}$, and the zero element is the universal relation~${\nab_{\bn} = \bn\times\bn}$.  For $\oijn$ we write $\eta_{ij}\in\Eq_n$ for the equivalence whose only non-trivial class is $\{i,j\}$.
If $\si$ is an equivalence relation on a set $X$, then for any $x\in X$ we write~$[x]_\si$ for the $\si$-class of $x$.

\subsection{Partition monoids}\label{subsect:Pn}

\begin{figure}[ht]
\begin{center}
\begin{tikzpicture}[scale=.5]
\begin{scope}[shift={(0,0)}]	
\uvs{1,...,6}
\lvs{1,...,6}
\uarcx14{.6}
\uarcx23{.3}
\uarcx56{.3}
\darc12
\darcx26{.6}
\darcx45{.3}
\stline34
\draw(0.6,1)node[left]{$a=$};
\draw[->](7.5,-1)--(9.5,-1);
\end{scope}
\begin{scope}[shift={(0,-4)}]	
\uvs{1,...,6}
\lvs{1,...,6}
\uarc12
\uarc34
\darc56
\darc23
\stline31
\stline55
\draw(0.6,1)node[left]{$b=$};
\end{scope}
\begin{scope}[shift={(10,-1)}]	
\uvs{1,...,6}
\lvs{1,...,6}
\uarcx14{.6}
\uarcx23{.3}
\uarcx56{.3}
\darc12
\darcx26{.6}
\darcx45{.3}
\stline34
\draw[->](7.5,0)--(9.5,0);
\end{scope}
\begin{scope}[shift={(10,-3)}]	
\uvs{1,...,6}
\lvs{1,...,6}
\uarc12
\uarc34
\darc56
\stline31
\stline55
\darc23
\end{scope}
\begin{scope}[shift={(20,-2)}]	
\uvs{1,...,6}
\lvs{1,...,6}
\uarcx14{.6}
\uarcx23{.3}
\uarcx56{.3}
\darc15
\darc56
\stline2{.95}
\darcx23{.2}
\draw(6.4,1)node[right]{$=ab$};
\end{scope}
\end{tikzpicture}
\caption{Multiplication of partitions $a,b\in\P_6$, with the product graph $\Pi(a,b)$ shown in the middle step.}
\label{fig:P6}
\end{center}
\end{figure}
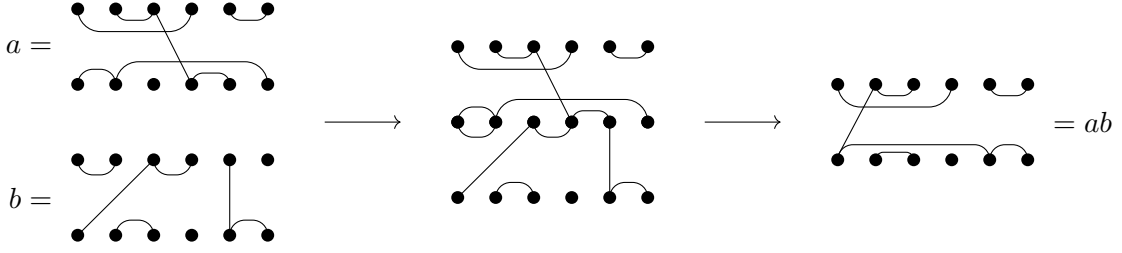

Fix an integer $n\geq0$, and again write $\bn = \{1,\ldots,n\}$.  Also fix two disjoint copies $\bn' = \{1',\ldots,n'\}$ and ${\bn'' = \{1'',\ldots,n''\}}$ of $\bn$.  The \emph{partition monoid} $\P_n$ consists of all set partitions of $\bn\cup\bn'$ under a multiplication we describe shortly.  First, we identify a partition $a\in\P_n$ with any graph on vertex set $\bn\cup\bn'$ whose connected components are the blocks of $a$.  When depicting such graphs, we typically draw vertices $\bn$ and $\bn'$ on an upper and lower row, in increasing order from left to right: $1<\cdots<n$ and $1'<\cdots<n'$.  See Figure \ref{fig:P6} for a depiction of two partitions from $\P_6$, one of them being
\[
a = \big\{ \{1,4\}, \{2,3,4',5'\}, \{5,6\}, \{1',2',6'\}, \{3'\} \big\}.
\]

To multiply partitions $a,b\in\P_n$, we first define three new graphs:
\bit
\item $a^\downarrow$, on vertex set $\bn\cup\bn''$, obtained by changing each lower vertex $x'$ of $a$ to $x''$,
\item $b^\uparrow$, on vertex set $\bn''\cup\bn'$, obtained by changing each upper vertex $x$ of $b$ to $x''$,
\item $\Pi(a,b)$, on vertex set $\bn\cup\bn''\cup\bn'$, whose edge set is the union of the edge sets of $a^\downarrow$ and $b^\uparrow$.
\eit
We call $\Pi(a,b)$ the \emph{product graph} of $a$ and $b$, and depict it with the vertices $\bn''$ on a middle row.  The product $ab\in\P_n$ is then defined to be the unique partition such that $x,y\in\bn\cup\bn'$ belong to the same block of $ab$ if and only if they are in the same connected component of $\Pi(a,b)$.  An example calculation is given in Figure \ref{fig:P6}.

The identity of $\P_n$ is the partition $1 = \bigset{\{x,x'\}}{x\in\bn} = \custpartn{1,2,4}{1,2,4}{\stline11\stline22\stline44\udotted24\ldotted24}$.  A unit of $\P_n$ has the form $\bigset{\{x,(xf)'\}}{x\in\bn}$ for some permutation $f\in\S_n$.  We identify this partition with $f$ itself, so in this way the group of units of $\P_n$ is $\S_n$.

A partition $a\in\P_n$ is \emph{planar} if it can be represented by a graph (as above), with all edges inside the rectangle spanned by the vertices, and with no edge-crossings.  For example, $b\in\P_6$ in Figure \ref{fig:P6} is planar, but $a\in\P_6$ is not.  The set
\[
\PP_n = \set{a\in\P_n}{a\text{ is planar}}
\]
is a submonoid of $\P_n$, called the \emph{planar partition monoid}.  For any subset $\Si$ of $\P_n$ we write
\[
\sP \Si = \Si\cap\PP_n = \set{a\in\Si}{a\text{ is planar}}  
\]
for the set of all planar elements of $\Si$.  When $\Si$ is a submonoid of $\P_n$, so too is $\sP\Si$.

The \emph{domain} of a partition $a\in\P_n$, denoted $\dom(a)$, is the set of all elements $x\in\bn$ for which the block of $a$ containing $x$ contains at least one element of $\bn'$.  For example, we have $\dom(a)=\{2,3\}$ and $\dom(b)=\{3,4,5\}$ for $a,b\in\P_6$ in Figure \ref{fig:P6}.  The set
\[
\Pnfd = \set{a\in\P_n}{\dom(a) = \bn}
\]
forms a submonoid of $\P_n$, called the \emph{full-domain partition monoid}, with group of units $\S_n$.

\section{Statement of the main results}\label{sect:main}

We can now state the main results of the paper, which are presentations for:
\bit
\item the full-domain partition monoid $\Pnfd$ (Theorem \ref{thm:Pnfd}),
\item the singular ideal $\Sing(\Pnfd) = \Pnfd\sm\S_n$ (Theorem \ref{thm:SingPnfd}), and 
\item the planar submonoid $\PPnfd = \Pnfd\cap\PP_n$ (Theorem \ref{thm:PPnfd}).
\eit
Note that the only unit of $\PPnfd$ is the identity element $1$, so the singular ideal is simply $\Sing(\PPnfd) = \PPnfd\sm\{1\}$, and there is no need to state a presentation for this.  To avoid trvialities, we assume henceforth that $n\geq2$.

Define an alphabet
\begin{equation}\label{eq:X}
X = \set{e_{ij} = e_{ji}}{\oijn} \cup \set{t_{ij}, t_{ji}}{\oijn},
\end{equation}
and a morphism $\phi: X^+\to\Sing(\Pnfd) : x\mt \ol x$, where 
\begin{equation}\label{eq:tijeij}
\eb_{ij} = \custpartn{1,3,4,5,7,8,9,11}{1,3,4,5,7,8,9,11}{\stline11\stline33\stline44\stline55\stline77\stline88\stline99\stline{11}{11}\uarc48\darc48\udotted13\udotted57\udotted9{11}\ldotted13\ldotted57\ldotted9{11} \vertlab11\vertlab4i\vertlab8j\vertlab{11}n}
\AND
\tb_{ij} = 
\begin{cases}
\custpartn{1,3,4,5,7,8,9,11}{1,3,4,5,7,8,9,11}{\stline11\stline33\stline44\stline55\stline77\stline84\stline99\stline{11}{11}\udotted13\udotted57\udotted9{11}\ldotted13\ldotted57\ldotted9{11} \vertlab11\vertlab4i\vertlab8j\vertlab{11}n} &\text{if $i<j$,}\\[5mm]
\custpartn{1,3,4,5,7,8,9,11}{1,3,4,5,7,8,9,11}{\stline11\stline33\stline48\stline55\stline77\stline88\stline99\stline{11}{11}\udotted13\udotted57\udotted9{11}\ldotted13\ldotted57\ldotted9{11} \vertlab11\vertlab4j\vertlab8i\vertlab{11}n} &\text{if $j<i$.}
\end{cases}
\end{equation}
Note that each $\ol x$ ($x\in X$) is an idempotent.  Also let $R$ be the following set of relations over $X$:
\begin{multicols}{2}
\noindent
\begin{align}
\label{TT1} t_{ij}^2 = t_{ij} &= t_{ji}t_{ij},  \\
\label{TT2} t_{ij}t_{kl} &= t_{kl}t_{ij},  \\
\label{TT3} t_{ik}t_{jk} &= t_{ik},  \\
\label{TT4} t_{ij}t_{ik} = t_{ik}t_{ij} &= t_{jk}t_{ij},  \\
\label{TT5} t_{ki}t_{ij}t_{jk} &= t_{ik}t_{kj}t_{ji}t_{ik},  \\
\label{TT6} t_{ki}t_{ij}t_{jk}t_{kl} &= t_{ik}t_{kl}t_{li}t_{ij}t_{jl}, 
\end{align}
\begin{align}
\label{EE1} e_{ij}^2 &= e_{ij},  \\
\label{EE2} e_{ij}e_{kl} &= e_{kl}e_{ij},  \\
\label{EE3} e_{ij}e_{jk} &= e_{jk}e_{ki},  \\
\label{ET1}	 e_{ij}t_{ij} &= t_{ij},  \\
\label{ET2}	 e_{jk}t_{ij} &= t_{ij}e_{ik},  \\
\label{ET3}	 e_{kl}t_{ij} &= t_{ij}e_{kl},  \\
\label{ET4}	 t_{ij}e_{ij} &= e_{ij},
\end{align}
\end{multicols}
\noindent where $i,j,k,l\in\bn$ are distinct in all of the above relations, except for \eqref{EE2} which only requires $i\ne j$ and $k\ne l$. Here is our first main result:

\begin{thm}\label{thm:SingPnfd}
The semigroup $\Sing(\Pnfd) = \Pnfd\sm\S_n$ has presentation $\Sgp\pres XR$ via $\phi$.  Consequently, $\Sing(\Pnfd)$ is idempotent-generated.
\end{thm}

Next, define an alphabet $Y = \{s_1,\ldots,s_{n-1},e,t\}$, and a morphism
\[
\psi:Y^*\to\Pnfd:y\mt\ul y,
\]
where
\begin{equation}\label{eq:siet}
\ul s_i = \custpartn{1,3,4,5,6,8}{1,3,4,5,6,8}{\stline11\stline33\stline45\stline54\stline66\stline88\udotted13\udotted68\ldotted13\ldotted68 \vertlab11\vertlab4i\vertlab8n}
\COMMA
\ul t = \ol t_{12} = \custpartn{1,2,3,5}{1,2,3,5}{\stline11\stline21\stline33\stline55\udotted35\ldotted35\vertlab11\vertlab5n}
\AND
\ul e = \ol e_{12} = \custpartn{1,2,3,5}{1,2,3,5}{\stline11\stline22\uuline12\ddline12\stline33\stline55\udotted35\ldotted35\vertlab11\vertlab5n}.
\end{equation}
(We use over- and under-line notation for the morphisms $\phi$ and $\psi$ in order to reduce confusion in Section \ref{sect:Pnfd}, where both morphisms will be considered at the same time.)
Also let $Q$ be the following set of relations over $Y$, in which we write $w=s_2s_3s_1s_2$:
\begin{multicols}{2}
\noindent
\begin{align}
\label{P1} s_i^2 &= \iota, \\
\label{P2} s_is_j &= s_js_i &&\hspace{-3mm}\text{ if } |i-j| > 1,  \\
\label{P3} s_is_js_i &= s_js_is_j &&\hspace{-3mm}\text{ if } |i-j| = 1,  \\
\label{P4} t^2 = t&=et=s_1t ,  \\
\label{P5} e^2=e = te  &= s_1e = es_1,  \\
\label{P6} s_it &= ts_i &&\hspace{-3mm}\text{ if } i>2,  \\
\label{P7} s_ie &= es_i &&\hspace{-3mm}\text{ if } i>2,  
\end{align}
\begin{align}
\label{P8} ts_1s_2t &= ts_1s_2s_1, \\
\label{P9} ts_2ts_2 &= s_2ts_2t ,  \\
\label{P10} es_2es_2 &= s_2es_2e,  \\
\label{P11} ts_2es_2 &= s_2es_2t ,  \\
\label{P12} tw tw &= w tw t ,  \\
\label{P13} ew ew &= w ew e ,  \\
\label{P14} tw ew &= w ew t . 
\end{align}
\end{multicols}
\noindent Here is our second main result:

\begin{thm}\label{thm:Pnfd}
The monoid $\Pnfd$ has presentation $\Mon\pres YQ$ via $\psi$.
\end{thm}

Finally, define an alphabet $Z = \set{f_i,g_i,h_i}{1\leq i<n}$, and a morphism
\[
\xi:Z^*\to\PPnfd:z\mt \ol z,
\]
where
\begin{equation}\label{eq:figihi}
\ol f_i = \custpartn{1,3,4,5,6,8}{1,3,4,5,6,8}{\stline11\stline33\stline44\stline54\stline66\stline88\udotted13\udotted68\ldotted13\ldotted68 \vertlab11\vertlab4i\vertlab8n}
\COMMA
\ol g_i = \custpartn{1,3,4,5,6,8}{1,3,4,5,6,8}{\stline11\stline33\stline45\stline55\stline66\stline88\udotted13\udotted68\ldotted13\ldotted68 \vertlab11\vertlab4i\vertlab8n}
\AND
\ol h_i = \custpartn{1,3,4,5,6,8}{1,3,4,5,6,8}{\stline11\stline33\stline44\stline55\uuline45\ddline45\stline66\stline88\udotted13\udotted68\ldotted13\ldotted68 \vertlab11\vertlab4i\vertlab8n}.
\end{equation}
These are again all idempotents.  Also let $O$ be the following set of relations over $Z$:
\bmc2
\noindent
\begin{align}
\label{O1} x_iy_i &= y_i &&\hspace{-22mm} \text{for $x,y\in\{f,g,h\}$,}\\
\label{O2} f_if_j &= f_jf_i &&\hspace{-12.5mm} \text{if $|i-j|>1$,}\\
\label{O3} g_ig_j &= g_jg_i &&\hspace{-12.5mm} \text{if $|i-j|>1$,}\\
\label{O4} h_ih_j &= h_jh_i ,\\
\label{O5} f_if_{i+1}f_i &= f_{i+1}f_if_{i+1} = f_{i+1}f_i,\\
\label{O6} g_ig_{i+1}g_i &= g_{i+1}g_ig_{i+1} = g_ig_{i+1},
\end{align}
\begin{align}
\label{O7} f_ig_j &= g_jf_i &&\text{if $j\not\in\{i,i+1\}$,}\\
\label{O8} h_if_j &= f_jh_i &&\text{if $j\not\in\{i,i-1\}$,}\\
\label{O9} h_ig_j &= g_jh_i &&\text{if $j\not\in\{i,i+1\}$,}\\
\label{O10} f_ig_{i+1} &= f_i,\\
\label{O11} g_{i+1}f_i &= g_{i+1},\\
\label{O12} h_ig_{i+1} &= h_{i+1}f_i.
\end{align}
\emc
\noindent Here is our third main result:

\begin{thm}\label{thm:PPnfd}
The monoid $\PPnfd$ has presentation $\Mon\pres ZO$ via $\xi$.  Consequently, $\PPnfd$ is idempotent-generated.
\end{thm}

We will prove Theorems \ref{thm:SingPnfd}, \ref{thm:Pnfd} and \ref{thm:PPnfd} in Sections \ref{sect:SingPnfd}, \ref{sect:Pnfd} and \ref{sect:PPnfd}, respectively.  The proofs require some background from \cite{CDEGZ2023} and \cite{EG2021}, which we outline in Sections \ref{sect:AP}--\ref{sect:Pn}.

\begin{rem}
The defining relations in the above presentations are far from being irredundant.  For example, whenever a presentation contains a pair of relations of the form $uv=v$ and $vu=u$ (for words $u$ and $v$), the relations $u^2=u$ and $v^2=v$ are easily seen to be consequences.  Thus, all of the relations in the three presentations declaring certain generators to be idempotents could be removed.
\end{rem}

\section{Strong action pairs}\label{sect:AP}

The paper \cite{CDEGZ2023} identified a natural class of semigroup products, and provided tools for constructing presentations for these semigroups under certain conditions.  The semigroups we are interested in arise in this way, so here we review the relevant results from \cite{CDEGZ2023}.

\subsection{Definitions and basic properties}

\begin{defn}[cf.~{\cite[Definition 4.5]{CDEGZ2023}}]\label{defn:AP}
A \emph{strong (right) action pair} in a monoid $M$ is a pair $(U,S)$ consisting of a submonoid $U$ of $M$ and a subsemigroup $S$ of $M$, such that
\begin{enumerate}[label=\textup{\textsf{(A\arabic*)}},leftmargin=9mm]
\item \label{A1} $Us \sub sU$ for all $s\in S$,
\item \label{A2} $su=tv \implies u=v$ for all $u,v\in U$ and $s,t\in S$.
\end{enumerate}
\end{defn}

Note that \cite[Definition 4.5]{CDEGZ2023} concerns the dual notion of \emph{left} action pairs, and allows $U$ to be a sub\emph{semigroup} (rather than submonoid); the properties \ref{A1} and \ref{A2} need to be (slightly) modified in this more general case; the simpler version above suits our purposes here.  We now gather some information from \cite{CDEGZ2023}.  

Fix a strong action pair $(U,S)$ in a monoid $M$.  The set
\[
SU = \set{su}{s\in S,\ u\in U}
\]
is a subsemigroup of $M$, and contains $S$ (as $1\in U$).  If $S$ happens to be a submonoid, then $SU$ is a submonoid of $M$, and also contains~$U$.  We also note that
\[
U\cap S = \begin{cases}
\{1\} &\text{if $S$ is a submonoid of $M$,}\\
\es &\text{otherwise.}
\end{cases}
\]
This is because from $s\in U\cap S$ and $s\cdot s = s^2\cdot 1$, it follows from \ref{A2} that $s=1$.

For $u\in U$ and $s\in S$, it follows from \ref{A1} and \ref{A2} that $us = sv$ for a unique $v\in U$, which we denote by $v = u^s$.  That is,
\begin{equation}\label{eq:act}
u\cdot s = s\cdot u^s \qquad\text{for $u\in U$ and $s\in S$.}
\end{equation}
This defines a right action of~$S$ on $U$ by monoid morphisms, meaning that
\[
u^{st} = (u^s)^t \COMMA (uv)^s = u^sv^s \AND 1^s = 1 \qquad\text{for all $u,v\in U$ and $s,t\in S$.}
\]
If $S$ happens to be a submonoid of $U$, then the action is monoidal, meaning that
\[
u^1 = u \qquad\text{for all $u\in U$.}
\]
Note that the operation in the semigroup $SU$ obeys the rule
\[
su \cdot tv = st \cdot u^tv \qquad\text{for $u,v\in U$ and $s,t\in S$.}
\]
This leads to a realisation of $SU$ as a concrete homomorphic image of a semidirect product $S\ltimes U$; see \cite[Theorem 4.53]{CDEGZ2023}.

Each element $u\in U$ determines a left congruence on $S^1$, which we denote by
\begin{equation}\label{eq:th}
\th_u = \set{(s,t)\in S^1\times S^1}{su=tu}.
\end{equation}
Generating sets for these left congruences will be important in what follows.  The next result gives a simple criterion for a left congruence $\th_u$ to be principal, i.e.~generated by a single pair.

\begin{lemma}\label{lem:1s}
If $(U,S)$ is a strong action pair, and if $u\in U$ and $s\in S^1$ are such that $us=s$ and $su=u$, then $\th_u = (1,s)\lc$.
\end{lemma}

\pf
Write $\si=(1,s)\lc$.  Since $1u=u=su$ we have $(1,s)\in\th_u$, and so $\si\sub\th_u$.  Conversely, suppose $(a,b)\in\th_u$, so that $au=bu$.  We then have $as=aus=bus=bs$, and so
\[
a = a1 \mr\si as = bs \mr\si b1 = b,
\]
showing that $\th_u\sub\si$.
\epf

Principal left congruences also feature in the following technical lemma, which will be useful in Section \ref{sect:PPnfd}.  The statement refers to the join $\bigvee_{v\in W}\th_v$ of right congruences, which coincides with their join as equivalence relations, i.e.~the least equivalence containing $\bigcup_{v\in W}\th_v$.

\begin{lemma}\label{lem:bigvee}
Suppose $(U,S)$ is a strong action pair, and let $u\in U$.  Suppose there exists an element $s\in S^1$, a subset $W\sub U$, and elements $s_v\in S^1\ (v\in W)$ such that all of the following hold:
\ben\bmc2
\item \label{bigvee1} $\th_u = (1,s)\lc$,
\item \label{bigvee2} $u\leqR v$ (in $U$) for all $v\in W$,
\item \label{bigvee3} $(1,s_v)\in\th_v$ for all $v\in W$, and 
\item \label{bigvee4} $s\in\pres{s_v}{v\in W}$.
\emc\een
Then $\th_u = \bigvee_{v\in W}\th_v$.
\end{lemma}

\pf
Write $\si = \bigvee_{v\in W}\th_v$.  To show that $\si\sub\th_u$, it is enough to show that
\[
\th_v\sub\th_u \qquad\text{for all $v\in W$.}
\]
To see this, fix some $v\in W$.  By \ref{bigvee2} we have $u=vw$ for some $w\in U$, and then
\[
(s,t)\in\th_v \Implies sv=tv \Implies su=(sv)w=(tv)w=tu \Implies (s,t)\in\th_u.
\]
To establish the inclusion $\th_u\sub\si$, it suffices by assumption \ref{bigvee1} to show that $(1,s)\in\si$.  But note that by assumption \ref{bigvee3} and $\th_v\sub\si$ we have $s_v\in[1]_\si$ for all $v\in W$.  Since $[1]_\si$ is a submonoid of $S^1$ (see \cite[Lemma 2.1]{CDEGZ2023}), it follows from assumption \ref{bigvee4} that $s\in[1]_\si$, i.e.~that $(1,s)\in\si$.
\epf

\subsection{Presentation: monoid case}\label{subsect:first}

The main results of \cite[Chapter 6]{CDEGZ2023} show how to construct presentations for a semigroup $SU$ arising from an action pair $(U,S)$, starting with presentations for $S$ and $U$.  Here we just state the relevant results for \emph{strong} action pairs, as that is sufficient for our purposes.  In our applications, certain special simplifications also arise, so we build these into the stated results from \cite{CDEGZ2023}.  

Fix a strong (right) action pair $(U,S)$ in a monoid $M$ (cf.~Definition \ref{defn:AP}), and assume here that $S$ is also a submonoid of $M$ (the case in which $S$ is not a submonoid will be treated in Section \ref{subsect:second}).
Now suppose we have presentations 
\[
\Mon\pres{X_U}{R_U} \AND \Mon\pres{X_S}{R_S} 
\]
for $U$ and $S$, respectively, via surmorphisms
\[
\phi_U:X_U^*\to U \AND \phi_S:X_S^*\to S.
\]
Without loss of generality, we assume that:
\bit
\item $X_U$ and $X_S$ are disjoint,
\item no element of $X_U$ or $X_S$ maps to $1$, and 
\item $\phi_U$ and $\phi_S$ are injective on $X_U$ and $X_S$.
\eit
We also fix normal form functions
\[
N_U:U\to X_U^* \AND N_S:S\to X_S^*,
\]
assuming that $N_U(1) = \io = N_S(1)$, and that $N_U(x\phi_U) = x$ and $N_S(y\phi_S) = y$ for all $x\in X_U$ and $y\in X_S$.

For each $x\in X_U$ and $y\in X_S$, write $x^y = N_U((x\phi_U)^{y\phi_S}) \in X_U^*$, and define
\begin{align}
\label{eq:R1} R_1 &= \set{(xy,y\cdot x^y)}{x\in X_U,\ y\in X_S}.
\intertext{For each $u\in U$, fix a generating set $\Om_u \sub S\times S$ for the right congruence $\th_u$ on $S(=S^1)$ from~\eqref{eq:th}.  For any subset $V\sub U$, define}
\label{eq:R2} R_2(V) &= \set{(N_S(s)N_U(v),N_S(t)N_U(v))}{v\in V,\ (s,t)\in\Om_v}.
\intertext{Note that when $V = X_U\phi_U$, the set $R_2(V)$ has the simpler form}
\label{eq:R2XU} R_2(X_U\phi_U) &= \set{(N_S(s)x,N_S(t)x)}{x\in X_U,\ (s,t)\in\Om_{x\phi_U}}.
\end{align}
The next result follows from (the left-right dual of) Theorem 6.5 of \cite{CDEGZ2023}, combined with Remark~6.6 (see in particular equations (6.8) and (6.9) in that remark).

\begin{thm}\label{thm:SU1}
Suppose $(U,S)$ is a strong action pair in a monoid $M$, with $S$ a submonoid of~$M$.  Suppose also that one of the following two assumptions holds:
\ben
\item \label{SU11} there exists a subset $V\sub U$ such that
\[
\th_u = \bigvee_{v\in V,\atop u\leqR \hspace{0.5mm} v} \th_v \qquad\text{for all $u\in U$, \qquad or} 
\]
\item \label{SU12} $U$ is commutative, and $\th_{uv} = \th_u\vee\th_v$ for all $u,v\in U$.  In this case, let $V = X_U\phi_U$.
\een
Then the monoid $SU$ has presentation
\[
\Mon\pres{X_U\cup X_S}{R_U\cup R_S\cup R_1\cup R_2(V)}
\]
via
\[
\phi:(X_U\cup X_S)^*\to SU:x\mt \begin{cases}
x\phi_U &\text{if $x\in X_U$,}\\
x\phi_S &\text{if $x\in X_S$.}
\end{cases}
\]
\end{thm}

Note that the set $V$ in assumption \ref{SU11} of Theorem \ref{thm:SU1} need not be a generating set for $U$ in general.  In particular, this will be the case when we apply Theorem \ref{thm:SU1} to the monoid $\PPnfd$ in Section \ref{sect:PPnfd}.

\subsection{Presentation: semigroup case}\label{subsect:second}

The case in which $S$ is not a submonoid of the over-monoid requires a little more care in the set-up, but with a similar end result.  Because there are a few extra assumptions, it will be convenient to state these up front.  Specifically, we fix an action pair $(U,S)$ in a monoid $M$, and assume that:
\bit
\item $S$ is \emph{not} a submonoid of $M$,
\item $U\sm\{1\}$ is a subsemigroup of $U$, and $U\sm\{1\} \sub SU$, and
\item $U$ is commutative, and $\th_{uv} = \th_u\vee\th_v$ for all $u,v\in U$.
\eit
We again fix presentations
\[
\Mon\pres{X_U}{R_U} \AND \Sgp\pres{X_S}{R_S} 
\]
for $U$ and $S$, respectively, via surmorphisms
\[
\phi_U:X_U^*\to U \AND \phi_S:X_S^+\to S.
\]
We again assume that $X_U$ and $X_S$ are disjoint, that no element of $X_U$ maps to $1$, and that $\phi_U$ and $\phi_S$ are injective on $X_U$ and $X_S$.  We also fix normal form functions 
\[
N_U:U\to X_U^* \AND N_S:S\to X_S^+,
\]
with analogous assumptions concerning $1$ and generators from $X_U$ and $X_S$.  Even though $1$ is not an element of $S$, and $\io$ is not an element of $X_S^+$, we also write $N_S(1) = \io$.  

We define the set of relations $R_1$ exactly as in \eqref{eq:R1}.  Note that it is possible to have $x^y = \io$ for some $x\in X_U$ and $y\in X_S$, in which case the corresponding relation from $R_1$ then has the form $(xy,y)$.  
We also define
\[
R_2 = R_2(X_U\phi_U),
\]
as in \eqref{eq:R2XU}.  Note that some set $\Om_{x\phi_U}$ (for $x\in X_U$) might contain a pair of the form $(s,1)$ or $(1,s)$, in which case the corresponding relation from $R_2$ then has the form $(N_S(s)x,x)$ or $(x,N_S(s)x)$.

\begin{thm}\label{thm:SU2}
With the above assumptions and notation, the semigroup $SU$ has presentation
\[
\Sgp\pres{X_U\cup X_S}{R_U\cup R_S\cup R_1\cup R_2}
\]
via
\[
\phi:(X_U\cup X_S)^+\to SU:x\mt \begin{cases}
x\phi_U &\text{if $x\in X_U$,}\\
x\phi_S &\text{if $x\in X_S$.}
\end{cases}
\]
\end{thm}

\pf
This follows from (the left-right dual of) of \cite[Theorem 6.44(ii)]{CDEGZ2023}.  Note that the set $R_3$ of relations from the stated result of \cite{CDEGZ2023} is empty because the action pair $(U,S)$ is strong.
\epf

\section{Ehresmann and restriction semigroups}\label{sect:Ehresmann}

The full-domain partition monoids $\Pnfd$ arose in \cite{EG2021} in the context of right restriction submonoids of Ehresmann monoids.  Right restriction monoids also contain natural action pairs, as discussed in \cite{CDEGZ2023}, and these will be very useful in our investigations.  Here we gather the relevant definitions and basic properties we will need.  See \cite{Lawson1991,Lawson2021} for more details on Ehresmann monoids.

\begin{defn}
An \emph{Ehresmann semigroup} is a tuple $(S,\cdot,D,R)$, where $(S,\cdot)$ is a semigroup, and where $D$ and $R$ are unary operations on $S$ satifsying the identities:
\begin{enumerate}[label=\textup{\textsf{(E\arabic*)}},leftmargin=9mm]\bmc2
\item \label{E1} \quad $D(a)a = a$,
\item \label{E2} \quad $D(a)D(b) = D(b)D(a)$,
\item \label{E3} \quad $D(ab) = D(aD(b))$,
\item \label{E4} \quad $D(ab) = D(a)D(ab)$,
\item \label{E5} \quad $R(D(a)) = D(a)$,
\item[] $aR(a) = a$,
\item[] $R(a)R(b) = R(b)R(a)$,
\item[] $R(ab) = R(R(a)b)$,
\item[] $R(ab) = R(ab)R(b)$,
\item[] $D(R(a)) = R(a)$.
\emc
\end{enumerate}
An \emph{Ehresmann monoid} is an Ehresmann semigroup that happens to be a monoid.  In this case it follows from \ref{E1} that $D(1)=R(1)=1$.
\end{defn}

It is well known, and easy to show, that \ref{E1}--\ref{E5} imply the following additional identities:
\begin{enumerate}[label=\textup{\textsf{(E\arabic*)}},leftmargin=9mm]\bmc2\addtocounter{enumi}{5}
\item \label{E6} \quad $D(D(a)) = D(a)$,
\item \label{E7} \quad $D(a)^2=D(a)$,
\item \label{E8} \quad $D(a)D(b) = D(D(a)D(b))$,
\item[] $R(R(a)) = R(a)$,
\item[] $R(a)^2 = R(a)$,
\item[] $R(a)R(b) = R(R(a)R(b))$.
\emc
\end{enumerate}
In particular, the set
\[
P = P(S) = \set{D(a)}{a\in S} = \set{R(a)}{a\in S}
\]
is a semilattice, i.e.~a semigroup of commuting idempotents.  The order on $P$ is given by
\[
p\leq q \iff p = pq \iff p = qp \iff p = pq = qp.
\]
We also have
\[
P = \set{p\in S}{p^2=p=D(p)=R(p)},
\]
as observed in \cite[Proposition 1.2]{Stokes2015}.  The elements of $P$ are called \emph{projections} in the literature, hence the notation.

In what follows, we write $a=_1b$ to indicate that elements $a$ and $b$ of an Ehresmann semigroup are equal as a consequence of \ref{E1}, with similar meanings for $=_2$, $=_3$ and so on.

\begin{lemma}\label{lem:Rab=Rb}
If $S$ is an Ehresmann semigroup, and if $a,b\in S$ are such that $R(a) \geq D(b)$, then $R(ab) = R(b)$.
\end{lemma}

\pf
The assumption gives $R(a)D(b) = D(b)$, and we then calculate
\[
R(ab) =_3 R(R(a)b) =_1 R(R(a)D(b)b) = R(D(b)b) =_1 R(b).  \qedhere
\]
\epf

\begin{cor}\label{cor:Rap}
If $S$ is an Ehresmann semigroup, and if $a\in S$ and $p\in P(S)$ are such that $R(a) \geq p$, then $R(ap) = p$.
\end{cor}

\pf
Take $b=p$ in Lemma \ref{lem:Rab=Rb}, noting that $D(p)=p=R(p)$.
\epf

The next result identifies a number of important substructures of an Ehresmann monoid.

\newpage

\begin{prop}\label{prop:MT}
If $M$ is an Ehresmann monoid, then 
\ben
\item \label{MT1} $T(M) = \set{a\in M}{R(a) = 1}$ is a submonoid of $M$,
\item \label{MT2} $I(M) = \set{a\in M}{D(a) \not= 1}$ is a right ideal of $M$,
\item \label{MT3} $\Tf(M) = \set{a\in M}{D(a) \not= 1 = R(a)}$ is a subsemigroup of $M$.
\een
\end{prop}

\pf
\firstpfitem{\ref{MT1}}  Since $R(1) = 1$, it suffices to show that $R(a)=R(b)=1$ implies $R(ab)=1$.  But here we clearly have $R(a) = 1 \geq D(b)$, and so it follows from Lemma \ref{lem:Rab=Rb} that $R(ab) = R(b) = 1$.

\pfitem{\ref{MT2}}  Let $a,b\in M$ with $D(a) \not= 1$.  We then have $D(ab) =_4 D(a)D(ab) \leq D(a) <1$.

\pfitem{\ref{MT3}}  This follows from \ref{MT1} and \ref{MT2}, as $\Tf(M) = T(M) \cap I(M)$.
\epf

The next lemma gives an important kind of cancellation property satisfied in Ehresmann monoids, which provides a link to the action pairs of Section \ref{sect:AP}.

\begin{lemma}\label{lem:SA2}
If $M$ is an Ehresmann monoid, and if $a,b\in T(M)$ and $p,q\in P(M)$, then
\[
ap = bq \implies p=q.
\]
\end{lemma}

\pf
Since $R(a) = R(b) = 1 \geq p,q$, Corollary \ref{cor:Rap} gives $R(ap) = p$ and $R(bq) = q$.  Since $ap=bq$ implies $R(ap) = R(bq)$, the result follows.
\epf

We have already mentioned that $\P_n$ is an Ehresmann monoid \cite{EG2021}.  The full-domain partition monoid $\Pnfd$ is an Ehresmann submonoid, but also satisfies the stronger property of being a right restriction monoid.

\begin{defn}
An Ehresmann semigroup is \emph{left restriction} if it additionally satisfies:
\begin{enumerate}[label=\textup{\textsf{(L)}},leftmargin=9mm]
\item \label{L} $aD(b) = D(ab)a$,
\end{enumerate}
or \emph{right restriction} if it satisfies the dual axiom:
\begin{enumerate}[label=\textup{\textsf{(R)}},leftmargin=9mm]
\item \label{R} $R(a)b = bR(ab)$.
\end{enumerate}
\end{defn}

It will be convenient to have an alternative characterisation of the right restriction property.  The third condition in the next result again links back to action pairs.

\begin{lemma}\label{lem:RR}
If $S$ is an Ehresmann semigroup, then the following are equivalent:
\ben
\item \label{RR1} $S$ is right restriction,
\item \label{RR2} $pa = aR(pa)$ for all $a\in S$ and $p\in P(S)$,
\item \label{RR3} $Pa\sub aP$ for all $a\in S$.
\een
\end{lemma}

\pf
\firstpfitem{\ref{RR1}$\iff$\ref{RR2}}  Given \ref{E3}, note that \ref{R} is equivalent to $R(a)b = bR(R(a)b)$.  Since ${P(S) = \set{R(a)}{a\in S}}$, this is equivalent to $pb = bR(pb)$, quantified over $b\in S$ and $p\in P(S)$.

\pfitem{\ref{RR2}$\iff$\ref{RR3}}  With the forward direction being trivial, suppose $Pa\sub aP$ for all $a\in S$.  Also let $a\in S$ and $p\in P(S)$.  By assumption we have $pa=aq$ for some $q\in P(S)$, so it follows that
\[
R(pa) = R(aq) =_3 R(R(a)q) =_8 R(a)q.
\]
This then gives $aR(pa) = aR(a)q =_1 aq = pa$.
\epf

\begin{cor}\label{cor:PS}
Let $M$ be a right restriction Ehresmann monoid, and let $P = P(M)$ and ${T = T(M)}$.  Then $(P,S)$ is a strong (right) action pair in $M$, for any subsemigroup $S$ of $T$, and the action of $S$ on $P$ from \eqref{eq:act} is given by
\[
p^a = R(pa) \qquad\text{for $p\in P$ and $a\in S$.}
\]
\end{cor}

\pf
Conditions \ref{A1} and \ref{A2} follow from Lemmas \ref{lem:RR} and \ref{lem:SA2}, respectively.  The formula for the action comes from the identity in Lemma \ref{lem:RR}\ref{RR2}.
\epf

We will be particularly interested in the action pairs $(P,T)$ and $(P,\Tf)$, and the corresponding products $TP$ and $\Tf P$.

\section{Partition monoids}\label{sect:Pn}

Now that we have reviewed the background we need on Ehresmann monoids in general, we look at the situation for the partition monoids $\P_n$.  After giving some further general background on~$\P_n$ in Section \ref{subsect:prelimPn}, we then recall in Section \ref{subsect:EhresmannPn} the Ehresmann structure on $\P_n$ discovered in \cite{EG2021}.  In Section~\ref{subsect:Pnfd} we focus on the full-domain partition monoid $\Pnfd$ and its singular ideal $\Sing(\Pnfd)$, identifying natural substructures that lead to action pairs and product decompositions.  Finally, in Section \ref{subsect:thu} we prove a number of results concerning the left congruences $\th_u$, which will be used to establish presentations for $\Pnfd$ and $\Sing(\Pnfd)$ in Sections \ref{sect:SingPnfd} and \ref{sect:Pnfd}.

\subsection{Further definitions and basic properties}\label{subsect:prelimPn}

A block of a partition $a\in\P_n$ is called a \emph{transversal} if it contains elements of both~$\bn$ and~$\bn'$, or a \emph{non-transversal} otherwise; these can be either \emph{upper} or \emph{lower}, with obvious meanings.  Partitions can be conveniently written in tabular notation.  Specifically, we write
\[
a = \begin{partn}{6} A_1&\cdots&A_r&C_1&\cdots&C_s\\ \hhline{~|~|~|-|-|-} B_1&\cdots&B_r&D_1&\cdots&D_t\end{partn}
\]
to indicate that $a$ has transversals $A_i\cup B_i'$ (for $1\leq i\leq r$), upper non-transversals $C_i$ (for $1\leq i\leq s$), and lower non-transversals $D_i'$ (for $1\leq i\leq t$).  It is possible for any of $r$, $s$ or $t$ to be~$0$ (but not all three, unless $n=0$).  At times we will use variations of this notation.  For example, if $a$ has no upper non-transversals then $a = \begin{partn}{6} A_1&\cdots&A_r& \multicolumn{3}{c}{} \\ \hhline{~|~|~|-|-|-} B_1&\cdots&B_r&D_1&\cdots&D_t\end{partn}$.  Or if all blocks of $a$ are transversals then $a = \begin{partn}{3} A_1&\cdots&A_r \\ B_1&\cdots&B_r \end{partn}$.  Sometimes we will omit singleton blocks from this notation; so for example $b\in\P_6$ from Figure \ref{fig:P6} can be written as $b = \begin{partn}{3} 3,4&5&1,2 \\ \hhline{~|~|-} 1&5,6&2,3\end{partn}$.

For a partition $a\in\P_n$ and an element $x\in\bn\cup\bn'$, we write $[x]_a$ for the block of $a$ containing~$x$.  The \emph{(co)domain} and \emph{(co)kernel} of $a$ are defined by
\begin{align*}
\dom(a) &= \set{x\in\bn}{[x]_a\cap\bn'\not=\es}, & \ker(a) &= \set{(x,y)\in\bn\times\bn}{[x]_a=[y]_a}, \\
\codom(a) &= \set{x\in\bn}{[x']_a\cap\bn\not=\es}, & \coker(a) &= \set{(x,y)\in\bn\times\bn}{[x']_a=[y']_a}.
\end{align*}
The \emph{rank} of $a$, denoted $\rank(a)$, is the number of transversals of $a$.
If $a = \begin{partn}{6} A_1&\cdots&A_r&C_1&\cdots&C_s\\ \hhline{~|~|~|-|-|-} B_1&\cdots&B_r&D_1&\cdots&D_t\end{partn}$, then 
\[
\rank(a) = r \COMMA \dom(a) = A_1\cup\cdots\cup A_r \AND \codom(a) = B_1\cup\cdots\cup B_r,
\]
while
\[
\bn/\ker(a) = \{A_1,\ldots,A_r,C_1,\ldots,C_s\} \AND \bn/\coker(a) = \{B_1,\ldots,B_r,D_1,\ldots,D_t\}.
\]
As we noted earlier, the group of units of $\P_n$ is (an isomorphic copy of) the \emph{symmetric group} $\S_n$, which can be characterised in terms of the above parameters as
\begin{align}
\nonumber \S_n &= \set{a\in\P_n}{\rank(a)=n} \\
\nonumber &= \set{a\in\P_n}{\dom(a)=\codom(a)=\bn,\ \ker(a)=\coker(a)=\De_\bn}\\
\label{eq:Sn} &= \set{a\in\P_n}{\dom(a)=\bn,\ \ker(a)=\De_\bn}. 
\intertext{The partition monoid also contains (isomorphic copies of) the \emph{full transformation monoid}}
\label{eq:Tn}
\T_n &= \set{a\in\P_n}{\dom(a)=\bn,\ \coker(a)=\De_\bn},
\end{align}
and the \emph{symmetric} and \emph{dual symmetric inverse monoids}, 
\[
\I_n = \set{a\in\P_n}{\ker(a) = \coker(a) = \De_\bn} \ANd \J_n = \set{a\in\P_n}{\dom(a) = \codom(a) = \bn}.
\]
A partition $a$ from $\T_n$ is identified with the mapping $\bn\to\bn$ that sends $x$ to $y$ where $y'$ is the unique element of $\bn'$ belonging to the block of $a$ containing $x$.  A typical element of $\T_n$ has the tabular form $\begin{partn}{3} A_1&\cdots&A_r \\ b_1&\cdots&b_r \end{partn}$, where we recall the above convention about omitting singleton blocks.

Elements of $\J_n$ are called \emph{block bijections}, and are the partitions from $\P_n$ whose blocks are all transversals.  Such a block bijection $a\in\J_n$ is \emph{uniform} if every transversal $A\cup B'$ of $a$ satisfies $|A|=|B|$.  The set
\begin{equation}\label{eq:Fn}
\F_n = \set{a\in\J_n}{a\text{ is uniform}}
\end{equation}
is a submonoid of $\J_n$, and hence of $\P_n$.  The monoid $\F_n$ will play an important role in Section \ref{sect:Pnfd}.

\subsection{Ehresmann structure}\label{subsect:EhresmannPn}

For an equivalence $\ve\in\Eq_n$, with $\bn/\ve = \{A_1,\ldots,A_r\}$, we define $\id_\ve = \begin{partn}{3} A_1&\cdots&A_r \\ A_1&\cdots&A_r \end{partn}$.  Note that all blocks of $\id_\ve$ are transversals.  It is easy to see that
\[
\id_\ve\id_\eta = \id_{\ve\vee\eta} \qquad\text{for $\ve,\eta\in\Eq_n$.}
\]
In particular, the set
\[
\E_n = \set{\id_\ve}{\ve\in\Eq_n}
\]
is a submonoid of $\P_n$, and is isomorphic to the join semilattice $(\Eq_n,\vee)$.  Note that $\id_{\De_{\bn}} = 1$ is the identity element of both $\E_n$ and $\P_n$.

For $a\in\P_n$, we define
\begin{equation}\label{eq:DR}
D(a) = \id_{\ker(a)} \AND R(a) = \id_{\coker(a)}.
\end{equation}
In tabular notation, if $a = \begin{partn}{6} A_1&\cdots&A_r&C_1&\cdots&C_s\\ \hhline{~|~|~|-|-|-} B_1&\cdots&B_r&D_1&\cdots&D_t\end{partn}$, then
\begin{equation}\label{eq:DaRa}
D(a) = \begin{partn}{6} A_1&\cdots&A_r&C_1&\cdots&C_s \\ A_1&\cdots&A_r&C_1&\cdots&C_s \end{partn}
\AND
R(a) = \begin{partn}{6} B_1&\cdots&B_r&D_1&\cdots&D_t \\ B_1&\cdots&B_r&D_1&\cdots&D_t \end{partn}.
\end{equation}
It follows from \cite[Theorem 4.8]{EG2021} that these operations give $\P_n$ the structure of an Ehresmann monoid, i.e.~that \ref{E1}--\ref{E5} hold.  See also \cite{MS2021}.  Note then that the set
\[
P(\P_n) = \set{D(a)}{a\in\P_n}= \set{R(a)}{a\in\P_n}
\]
is precisely the semilattice $\E_n$.

\subsection{The full-domain partition monoid}\label{subsect:Pnfd}

The Ehresmann monoid $\P_n$ is neither left nor right restriction when $n\geq2$.  For example, it is easy to check that \ref{L} and \ref{R} both fail for the partitions $a = \custpartn{1,2}{1,2}{}$ and $b = \custpartn{1,2}{1,2}{\stline11\stline22\uarc12\darc12}$ from~$\P_2$.

On the other hand, $\P_n$ contains a very natural right restriction Ehresmann submonoid, namely the \emph{full-domain partition monoid:}
\[
\Pnfd = \set{a\in\P_n}{\dom(a)=\bn},
\]
which will be the focus of our study.
(Dually, the full-codomain partition monoid is a left restriction Ehresmann monoid.)  Note that $\Pnfd$ contains the semilattice $\E_n$ and the group $\S_n$.  Several structural results concerning $\Pnfd$ were proved in \cite{EG2021}, including a description of Green's relations, and establishing its regularity and the right restriction property.

The right restriction monoid $\Pnfd$ contains the natural action pairs $(P,T)$ and $(P,\Tf)$, where
\[
P=P(\Pnfd)=\E_n \COMMA T=T(\Pnfd) \AND \Tf=\Tf(\Pnfd),
\]
are the substructures from Proposition \ref{prop:MT}.  The next result characterises these substructures, and also identifies the resulting products $TP$ and $\Tf P$.

\begin{prop}\label{prop:TPnfd}
We have
\ben
\item \label{TPnfd1} $T(\Pnfd) = \T_n$,
\item \label{TPnfd2} $I(\Pnfd) = \Sing(\Pnfd) = \Pnfd\sm\S_n$,
\item \label{TPnfd3} $\Tf(\Pnfd) = \Sing(\T_n) = \T_n\sm\S_n$.
\een
Further, we have the product decompositions
\[
\Pnfd = \T_n\E_n \AND \Sing(\Pnfd) = \Sing(\T_n)\E_n.
\]
\end{prop}

\pf
\firstpfitem{\ref{TPnfd1}}  By definition, and using \eqref{eq:Tn}, we have
\begin{align*}
T(\Pnfd) &= \set{a\in\Pnfd}{R(a)=1} \\
&= \set{a\in\Pnfd}{\coker(a)=\De_\bn} \\
&= \set{a\in\P_n}{\dom(a)=\bn,\ \coker(a)=\De_\bn} = \T_n.
\end{align*}

\pfitem{\ref{TPnfd2}}  By definition, and using \eqref{eq:Sn}, we have
\begin{align*}
I(\Pnfd) &= \set{a\in\Pnfd}{D(a)\not=1} = \set{a\in\P_n}{\dom(a)=\bn,\ \ker(a)\not=\De_\bn}= \Pnfd\sm\S_n .
\end{align*}

\pfitem{\ref{TPnfd3}}  This follows from \ref{TPnfd1} and \ref{TPnfd2}, since $\Tf(\Pnfd) = T(\Pnfd)\cap I(\Pnfd)$.

\aftercases

For the final assertion, let $a\in\Pnfd$, and write $a = \begin{partn}{6} A_1&\cdots&A_r& \multicolumn{3}{c}{} \\ \hhline{~|~|~|-|-|-} B_1&\cdots&B_r&D_1&\cdots&D_t\end{partn}$.  For each $1\leq i\leq r$, fix some $b_i\in B_i$, and let $b = \begin{partn}{3} A_1&\cdots&A_r \\ b_1&\cdots&b_r \end{partn}\in\T_n$.  Keeping \eqref{eq:DaRa} in mind, it is then easy to see that
\[
a = bR(a) \in \T_n\E_n.
\]
If $a$ belongs to $\Sing(\Pnfd)$, then $\rank(b) = \rank(a)<n$, so in fact $b\in\Sing(\T_n)$.
\epf

\begin{rem}
The monoid $\F_n$ in \eqref{eq:Fn} also has a product decomposition, as $\F_n = \E_n\S_n = \S_n\E_n$.
\end{rem}

\begin{rem}
Green's relations on $\P_n$ have simple descriptions in terms of co/domain, co/kernel and rank parameters; see for example \cite{Wilcox2007,FL2011}.  In particular, for partitions $a,b\in\P_n$ we have
\begin{equation}\label{eq:RinPn}
a\R b \IFF \dom(a)=\dom(b) \text{ and }\ker(a)=\ker(b).
\end{equation}
It follows that the elements of $\Pnfd$ can be characterised as the partitions $\R$-related to an element of $\E_n = P(\P_n)$.  Note, however, that for an arbitrary Ehresmann semigroup $S$, the set
\[
\set{a\in S}{a\R p\ (\exists p\in P(S))} = \set{a\in S}{a\R D(a)}
\]
need not be a subsemigroup of $S$ in general.  Indeed, small (counter)examples can be quickly found using Mace4 \cite{Mace4}.  On the other hand, the set
\[
\set{a\in S}{D(a)\R a\L R(a)}
\]
is an \emph{inverse} subsemigroup, as shown in \cite[Lemma 5.15]{Stein2017}.
\end{rem}

\subsection{The left congruences}\label{subsect:thu}

The next three results concern the left congruences $\th_u$ ($u\in \E_n$) from \eqref{eq:th} associated to the action pairs $(P,T)=(\E_n,\T_n)$ and $(P,\Tf)=(\E_n,\Sing(\T_n))$ in $M = \Pnfd$.  
\bit
\item Lemma \ref{lem:thu} characterises the congruences themselves.
\item Lemma \ref{lem:thuv} shows that the congruences have the required properties in order to apply Theorems~\ref{thm:SU1} and \ref{thm:SU2}.
\item Lemma \ref{lem:thuij} gives a generating set for the congruences $\th_u$ for certain simple elements $u$.
\eit

\begin{lemma}\label{lem:thu}
If $S = \T_n$ or $\Sing(\T_n)$, and if $u=\id_\ve$ for some $\ve\in\Eq_n$, then
\[
\th_u = \set{(a,b)\in S^1\times S^1}{(xa,xb)\in\ve\ \text{ for all } x\in\bn}.
\]
\end{lemma}

\pf
Let $a,b\in S^1$, and write $\bn/\ve = \{A_1,\ldots,A_k\}$.  We then have
\[
au = \set{A_ia^{-1}\cup A_i'}{i\in\bk},
\]
where as usual we write $A_ia^{-1} = \set{x\in\bn}{xa\in A_i}$.  It follows that
\[
(a,b)\in\th_u \IFF au=bu \IFF A_ia^{-1} = A_ib^{-1} \text{ for all $i\in\bk$.}
\]
It is easy to see that this latter condition is equivalent to having $(xa,xb)\in\ve$ for all $x\in\bn$.
\epf

\begin{lemma}\label{lem:thuv}
If $S = \T_n$ or $\Sing(\T_n)$, then
\[
\th_{uv} = \th_u\vee\th_v \qquad\text{for all $u,v\in\E_n$.}
\]
\end{lemma}

\pf
Let $u = \id_\ve$ and $v = \id_\eta$, where $\ve,\eta\in\Eq_n$.  Since $\E_n$ is commutative, we have $\th_u\vee\th_v\sub\th_{uv}$ (see \cite[p.~52]{CDEGZ2023}), so it remains to show that
\begin{equation}\label{eq:joinsub}
\th_{uv}\sub\th_u\vee\th_v.
\end{equation}
To do so, let $(a,b)\in\th_{uv}$.  Since $uv=\id_{\ve\vee\eta}$, it follows from Lemma \ref{lem:thu} that $(xa,xb)\in\ve\vee\eta$ for all $x\in\bn$.  We therefore have sequences
\begin{equation}\label{eq:chains}
\begin{matrix}
1a & = & h_{1,0} & \xra\ve & h_{1,1} & \xra\eta & h_{1,2} & \xra\ve & \cdots & \xra\eta & h_{1,k} & = & 1b , \\
2a & = & h_{2,0} & \xra\ve & h_{2,1} & \xra\eta & h_{2,2} & \xra\ve & \cdots & \xra\eta & h_{2,k} & = & 2b , \\
\vdots &&&&&&&&&&&& \vdots \\
xa & = & h_{x,0} & \xra\ve & h_{x,1} & \xra\eta & h_{x,2} & \xra\ve & \cdots & \xra\eta & h_{x,k} & = & xb , \\
\vdots &&&&&&&&&&&& \vdots \\
na & = & h_{n,0} & \xra\ve & h_{n,1} & \xra\eta & h_{n,2} & \xra\ve & \cdots & \xra\eta & h_{n,k} & = & nb .
\end{matrix}
\end{equation}
Since $\ve$ and $\eta$ are reflexive, we can indeed assume that the above sequences have the same length, and that they alternate between $\ve$- and $\eta$-relationships as shown.
Moreover, if we have $xa=ya$ and $xb=yb$ for some $x,y\in\bn$, then we assume that the sequences in rows $x$ and $y$ of \eqref{eq:chains} are identical.
We then use the above columns to define a sequence of transformations $c_0,c_1,\ldots,c_k\in\T_n$.  Specifically, for $0\leq j\leq k$, we define
\[
x c_j = h_{x,j} \qquad\text{for each $x\in\bn$.}
\]
If each $c_j$ belongs to $S^1$, then it again follows from Lemma \ref{lem:thu} that
\[
(c_j,c_{j+1}) \in \begin{cases}
\th_u &\text{if $j$ is even,}\\
\th_v &\text{if $j$ is odd.}
\end{cases}
\]
Thus, in this case we have
\[
\begin{matrix}
a &=& c_0 &\xra{\th_u}& c_1& \xra{\th_v} &c_2 &\xra{\th_u}& \cdots &\xra{\th_v} &c_k &=& b,
\end{matrix}
\]
meaning that $(a,b)\in\th_u\vee\th_v$, completing the proof of \eqref{eq:joinsub}.  In particular, there is nothing more to be done in the case that $S=\T_n(=S^1)$.

Thus, for the rest of the proof we assume that $S=\Sing(\T_n)$, and deal with the case that some transformation $c_j$ constructed above does not belong to $S^1 = \Sing(\T_n)\cup\{1\}$.  In this case we have $c_j\in\S_n\sm\{1\}$, but we note that $a,b\in S^1$.  Now, if $a=b=1$, then clearly $(a,b)\in\th_u\vee\th_v$.  So we assume without loss of generality that $b\not=1$.  In particular, $b$ is not injective, so without loss of generality we can assume that $1b=2b$.

\pfcase1  If we also had $1a=2a$, then by our assumption above, rows $1$ and $2$ of \eqref{eq:chains} would be identical.  This would imply that $1c_j=2c_j$ for each $0\leq j\leq k$, i.e.~that each $c_j\in\Sing(\T_n)\sub S^1$, completing the proof in this case.

\pfcase2  We are now left to deal with the case that $1a\not=2a$.  Swapping $1$ and $2$ if necessary, we can assume that either $a=1$ is the identity map, or else there exist distinct $y,z\in\{2,\ldots,n\}$ such that $ya=za$.  Since $(1a,1b)$ and $(2a,2b)$ are both in $\ve\vee\eta$, and since $1b=2b$, it follows that $(1a,2a)\in\ve\vee\eta$.  We therefore have a sequence
\[
\begin{matrix}
1a & = & h_0 & \xra\ve & h_1 & \xra\eta & h_2 & \xra\ve & \cdots & \xra\eta & h_l & = & 2a .
\end{matrix}
\]
We use this to define another sequence of transformations $d_0,d_1,\ldots,d_l\in\T_n$.  For $0\leq j\leq l$, we define
\[
xd_j = \begin{cases}
h_j &\text{if $x=1$,}\\
xa &\text{if $2\leq x\leq n$,}
\end{cases}
\]
and we note that each $d_j\in S^1 = \Sing(\T_n)\cup\{1\}$.  Indeed, if $a=1$, then $d_j$ is either $1$ (if $h_j=1$) or else singular (as then $1d_j=h_j=h_jd_j$).  On the other hand, if $a\not=1$, then with $2\leq y<z\leq n$ as above, we have $yd_j=ya=za=zd_j$, so that $d_j$ is singular.  We now define $a' = d_l$.  Another appeal to Lemma~\ref{lem:thu} gives
\[
\begin{matrix}
a &= &d_0 &\xra{\th_u} &d_1& \xra{\th_v} &d_2 &\xra{\th_u}& \cdots& \xra{\th_v} &d_l &= &a',
\end{matrix}
\]
so that $(a,a')\in\th_u\vee\th_v$.  Since $(a',a) \in \th_u\vee\th_v \sub \th_{uv}$ and $(a,b)\in\th_{uv}$, it follows by transitivity that $(a',b)\in\th_{uv}$.  Since $1a'=2a'$ and $1b=2b$, we see that $a'$ and $b$ satisfy the conditions of Case 1, so it follows that $(a',b)\in\th_u\vee\th_v$.  Combining this with $(a,a')\in\th_u\vee\th_v$ and transitivity, it follows that indeed $(a,b)\in\th_u\vee\th_v$.  This completes the proof.
\epf

The next result refers to the partitions $\eb_{ij}$ and $\tb_{ij}$ from \eqref{eq:tijeij}, and the proof uses the equivalences $\eta_{ij} \in \Eq_n$ introduced at the end of Section \ref{subsect:S}.

\begin{lemma}\label{lem:thuij}
If $S = \T_n$ or $\Sing(\T_n)$, and if $\ve\in\Eq_n$, then $\th_{\id_\ve} = (1,f)\lc$ for any idempotent $f\in E(\T_n)$ with $\ker(f)=\ve$.
In particular, for any $\oijn$ we have $\th_{\eb_{ij}} = (1,\ol t_{ij})\lc$.
\end{lemma}

\pf
The first assertion follows from Lemma \ref{lem:1s}, since $f\id_\ve = \id_\ve$ and $\id_\ve f = f$.
The second follows from the first, since $\eb_{ij} = \id_{\eta_{ij}}$ and $\ker(\ol t_{ij}) = \eta_{ij}$.
\epf

\section{\texorpdfstring{\boldmath Presentation for $\Sing(\Pnfd)$}{Presentation for $Sing(Pnfd)$}}\label{sect:SingPnfd}

To prove Theorem \ref{thm:SingPnfd}, we now work towards an application of Theorem \ref{thm:SU2}, utilising the strong action pair $(P,\Tf) = (\E_n,\Sing(\T_n))$ in the right restriction monoid $M = \Pnfd$.  First we note that:
\bit
\item $\E_n\sm\{1\}$ is a subsemigroup of $\E_n$ (as is the case for any monoid semilattice), 
\item $\E_n\sm\{1\} \sub \Sing(\Pnfd) = \Sing(\T_n)\E_n$ (cf.~Proposition \ref{prop:TPnfd}),
\item $\E_n$ is commutative, and $\th_{uv}=\th_u\vee\th_v$ for all $u,v\in\E_n$ (cf.~Lemma \ref{lem:thuv}).
\eit
Thus, Theorem \ref{thm:SU2} does indeed apply.  To use it, we need presentations $\Mon\pres{X_{\E_n}}{R_{\E_n}}$ and $\Sgp\pres{X_{\Sing(\T_n)}}{R_{\Sing(\T_n)}}$ for $\E_n$ and $\Sing(\T_n)$, via appropriate surmorphisms.  By \cite[Theorem~2]{FitzGerald2003} and \cite[Theorem 6]{East2013}, such presentations exist with:
\bit
\item $X_{\E_n} = \set{e_{ij}=e_{ji}}{\oijn}$, and $R_{\E_n}$ the set of relations \eqref{EE1}--\eqref{EE3},
\item $X_{\Sing(\T_n)} = \set{t_{ij},t_{ji}}{\oijn}$, and $R_{\Sing(\T_n)}$ the set of relations \eqref{TT1}--\eqref{TT6}.
\eit
Noting that $X_{\E_n}\cup X_{\Sing(\T_n)}$ is the alphabet $X$ from \eqref{eq:X}, Theorem \ref{thm:SU2} then tells us that $\Sing(\Pnfd) = \Sing(\T_n)\E_n$ has presentation
\begin{equation}\label{eq:XR...}
\Sgp\pres{X}{R_{\E_n} \cup R_{\Sing(\T_n)}\cup R_1\cup R_2},
\end{equation}
where $R_1$ is as in \eqref{eq:R1}, and where $R_2 = R_2(X_U\phi_U)$ is as in \eqref{eq:R2XU}.  Given Lemma \ref{lem:thuij}, we can take~$R_2$ to consist of relations \eqref{ET4}.  For $R_1$, we need words $e_{kl}^{t_{ij}} \in X_{\E_n}^*$ that map to $\ol e_{kl}^{\ol t_{ij}} = R(\eb_{kl}\tb_{ij})$.  One can easily check that
\[
 \eb_{kl}^{\tb_{ij}} = R(\eb_{kl}\tb_{ij}) = \begin{cases}
1 &\text{if $k=i$ and $l=j$,}\\
1 &\text{if $k=j$ and $l=i$,}\\
\eb_{il} &\text{if $k=j$ and $l\not=i$,}\\
\eb_{ik} &\text{if $l=j$ and $k\not=i$,}\\
\eb_{kl} &\text{if $\{i,j\}\cap\{k,l\}=\es$,}\\
\eb_{il} &\text{if $k=i$ and $l\not=j$,}\\
\eb_{ki} &\text{if $l=i$ and $k\not=j$.}
\end{cases}
\]
In these seven cases, and with the obvious choice of words $e_{kl}^{t_{ij}}$ (being empty or a single letter), the resulting relations $e_{kl} t_{ij} = t_{ij} e_{kl}^{t_{ij}}$ from $R_1$ respectively become
\begin{align}
\label{te1} e_{ij}t_{ij} &= t_{ij} &&\text{for distinct $i,j\in\bn$,}\\
\label{te2} e_{ji}t_{ij} &= t_{ij} &&\text{for distinct $i,j\in\bn$,}\\
\label{te3} e_{jl}t_{ij} &= t_{ij}e_{il} &&\text{for distinct $i,j,l\in\bn$,}\\
\label{te4} e_{kj}t_{ij} &= t_{ij}e_{ik} &&\text{for distinct $i,j,k\in\bn$,}\\
\label{te5} e_{kl}t_{ij} &= t_{ij}e_{kl} &&\text{for distinct $i,j,k,l\in\bn$,}\\
\label{te6} e_{il}t_{ij} &= t_{ij}e_{il} &&\text{for distinct $i,j,l\in\bn$,}\\
\label{te7} e_{ki}t_{ij} &= t_{ij}e_{ki} &&\text{for distinct $i,j,k\in\bn$.}
\end{align}
Keeping in mind the symmetric notation $e_{ij} = e_{ji}$, we see that:
\bit
\item \eqref{te1} and \eqref{te2} are the same sets of relations, and are exactly \eqref{ET1},
\item \eqref{te3} and \eqref{te4} are the same sets of relations, and are exactly \eqref{ET2},
\item \eqref{te5} is exactly \eqref{ET3},
\item \eqref{te6} and \eqref{te7} are the same sets of relations, but are not included in the set $R$ of relations \eqref{TT1}--\eqref{ET4}.
\eit
Thus, the presentation for $\Pnfd$ in \eqref{eq:XR...} contains the presentation $\Sgp\pres XR$ from Theorem \ref{thm:SingPnfd} plus relations \eqref{te7}.  It therefore remains to show that relations \eqref{te7} follow from $R$.  For this we see that for distinct $i,j,k\in\bn$ we can transform
\[
e_{ki}t_{ij} \to_{\eqref{ET1}} e_{ki}e_{ij}t_{ij} 
\to_{\eqref{EE3}} e_{ij}e_{jk}t_{ij} 
\to_{\eqref{ET2}} e_{ij}t_{ij}e_{ik} 
\to_{\eqref{ET1}} t_{ij}e_{ik} = t_{ij}e_{ki} .
\]
This completes the proof of Theorem \ref{thm:SingPnfd}.

\section{\texorpdfstring{\boldmath Presentation for $\mathcal P_n^\fd$}{Presentation for $Pnfd$}}\label{sect:Pnfd}

There are several ways one could prove Theorem \ref{thm:Pnfd}:
\bit
\item We could apply Theorem \ref{thm:SU1}, utilising the strong action pair $(P,T) = (\E_n,\T_n)$ in the right restriction monoid $M = \Pnfd = \T_n\E_n$, again using Lemma~\ref{lem:thuv} to verify the key assumption on joins of left congruences.  This would require a presentation $\Mon\pres{X_{\T_n}}{R_{\T_n}}$ for $\T_n$ (see for example \cite[Theorem 9]{Lavers1997}), and would result in a presentation for $\Pnfd$ over the alphabet $X_{\E_n}\cup X_{\T_n}$, which would then have to be transformed into $\Mon\pres YQ$.
\item We could apply \cite[Theorem 7.1]{JEinsn2}, which allows us to build a presentation for ${\Pnfd = \S_n \cup \Sing(\Pnfd)}$ from presentations for $\S_n$ and $\Sing(\Pnfd)$, using the fact that $\S_n$ is a sub(semi)group of $\Pnfd$, and $\Sing(\Pnfd)$ an ideal.  This would again involve a rather large generating set, and require a transformation.
\item We could use the approach of \cite{JErook,JEgrpm2}, where we begin with the desired presentation ${\Mon\pres YQ}$, and in some sense `simulate' the presentation $\Sgp\pres XR$ inside $\Mon\pres YQ$.
\eit
We take the latter approach here, as it appears to be the most direct.  

For the rest of this section, we fix
the alphabet $Y = \{s_1,\ldots,s_{n-1},e,t\}$, 
the set $Q$ of relations~\eqref{P1}--\eqref{P14}, and 
the morphism $\psi:Y^*\to\Pnfd:y\mt\ul y$, as in~\eqref{eq:siet}.
For $w\in Y^*$ we write $\ul w = w\psi \in \Pnfd$.  We also write ${\sim} = Q^\sharp$ for the congruence on $Y^*$ generated by the relations $Q$.  When demonstrating $\sim$-equivalence of words, we use subscripts to indicate which relation we are using.

We now define three sub-alphabets of $Q$:
\[
\XSn = \{s_1,\ldots,s_{n-1}\} \COMMA \XFn = \XSn\cup\{e\} \AND \XTn = \XSn\cup\{t\}.
\]

\begin{lemma}\label{lem:FT}\leavevmode
\ben
\item \label{FT1} For any $u,v\in\XFn^*$ we have $\ul u=\ul v \implies u\sim v$.
\item \label{FT2} For any $u,v\in\XTn^*$ we have $\ul u=\ul v \implies u\sim v$.  
\een
\end{lemma}

\pf
By \cite[Theorem 3]{FitzGerald2003} and \cite[Theorem 9]{Lavers1997}, the monoids $\F_n$ and $\T_n$ have presentations
\[
\Mon\pres{\XFn}{\RFn} \AND \Mon\pres{\XTn}{\RTn}
\]
via the restrictions $\psi\rest_{\XFn^*}$ and $\psi\rest_{\XTn^*}$, where:
\bit
\item $\RFn$ consists of all relations from $Q$ that involve only words over $\XFn$, and
\item $\RTn$ consists of all relations from $Q$ that involve only words over $\XTn$, together with the relation $ts_2ts_2 = ts_2t$.
\eit
It therefore suffices to show that $ts_2ts_2 \sim ts_2t$, and for this we have
\begin{align*}
ts_2ts_2 &\Psim4 ts_2s_1ts_2 \Psim1 ts_2s_1s_2s_2ts_2 \Psim3 ts_1s_2s_1s_2ts_2 \\
&\Psim8 ts_1s_2ts_2ts_2 
\Psim9 ts_1s_2s_2ts_2t \Psim1 ts_1ts_2t \Psim4 tts_2t \Psim4 ts_2t.  \qedhere
\end{align*}
\epf

The results of \cite{FitzGerald2003} and \cite{Lavers1997} mentioned in the above proof also imply that the image of~$\psi$ contains both $\F_n$ and $\T_n$.  Since $\F_n$ contains $\E_n$, and since $\T_n\E_n = \Pnfd$, it follows that $\psi$ is surjective.  To complete the proof of Theorem \ref{thm:Pnfd} it therefore remains to show that $\ker(\psi) = {\sim}$.  One inclusion is easily established:

\begin{lemma}\label{lem:sub}
We have ${\sim} \sub \ker(\psi)$.
\end{lemma}

\pf
This follows from $Q\sub\ker(\psi)$, which is established by showing diagrammatically that $\ul u=\ul v$ for every relation $(u,v)$ from $Q$.
\epf

The inclusion $\ker(\psi)\sub{\sim}$ is the subject of the rest of this section.

For a word $w = s_{i_1}\cdots s_{i_k}\in\XSn^*$, we write $w^{-1} = s_{i_k}\cdots s_{i_1}$, noting that ${ww^{-1} \sim \io \sim w^{-1}w}$ by~\eqref{P1}.  For $\oijn$ we define words
\begin{align*}
c_{ij} = (s_2\cdots s_{j-1})(s_1\cdots s_{i-1}) \COMMa
\ve_{ij} = \ve_{ji} = c_{ij}^{-1}ec_{ij} \COMMa
\tau_{ij} = c_{ij}^{-1}tc_{ij} \ANd
\tau_{ji} = c_{ij}^{-1}ts_1c_{ij}.
\end{align*}
Note that $c_{12} = \io$, $\ve_{12} = e$ and $\tau_{12} = t$.
In Figure \ref{fig:words} we illustrate the partitions $\ul\ve_{ij} = \ul c_{ij}^{-1}\ul e\ \ul c_{ij}$.  Together with similar diagrams for $\ul\tau_{ij}$ and $\ul\tau_{ji}$, it follows that
\begin{equation}\label{eq:olul}
\ul\ve_{ij} = \ol e_{ij} \AND \ul\tau_{ij} = \ol t_{ij} \qquad\text{for distinct $i,j\in\bn$.}
\end{equation}

\begin{figure}[ht]
\begin{center}
\begin{tikzpicture}[scale=.5]
\begin{scope}[shift={(0,4)}]	
\uvs{1,3,4,5,7,8,9,11}
\lvs{1,2,3,5,6,8,9,11}
\stline41
\stline82
\stline13
\stline35
\stline56
\stline78
\stline99
\stline{11}{11}
\ldotted35
\ldotted68
\ldotted9{11}
\udotted13
\udotted57
\udotted9{11}
\draw[|-] (0,2)--(0,0);
\draw(0,1)node[left]{\small $\ul c_{ij}^{-1}$};
\node () at (1,2.5) {\footnotesize $1$};
\node () at (4,2.5) {\footnotesize $i$};
\node () at (8,2.5) {\footnotesize $j$};
\node () at (11,2.5) {\footnotesize $n$};
\end{scope}
\begin{scope}[shift={(0,2)}]	
\stline11
\stline22
\stline33
\stline55
\stline66
\stline88
\stline99
\stline{11}{11}
\uuline12
\ddline12
\draw[|-] (0,2)--(0,0);
\draw(0,1)node[left]{\small $\ul e$};
\draw(14,1)node{$=$};
\end{scope}
\begin{scope}[shift={(0,0)}]	
\uvs{1,2,3,5,6,8,9,11}
\lvs{1,3,4,5,7,8,9,11}
\stline14
\stline28
\stline31
\stline53
\stline65
\stline87
\stline99
\stline{11}{11}
\udotted35
\udotted68
\udotted9{11}
\ldotted13
\ldotted57
\ldotted9{11}
\draw[|-|] (0,2)--(0,0);
\draw(0,1)node[left]{\small $\ul c_{ij}$};
\end{scope}
\begin{scope}[shift={(16,2)}]	
\uvs{1,3,4,5,7,8,9,11}
\lvs{1,3,4,5,7,8,9,11}
\stline11
\stline33
\stline44
\stline55
\stline77
\stline88
\stline99
\stline{11}{11}
\uarc48
\darc48
\ldotted13
\ldotted57
\ldotted9{11}
\udotted13
\udotted57
\udotted9{11}
\draw[|-|] (12,2)--(12,0);
\draw(12,1)node[right]{\small $\ul\ve_{ij}(=\ol e_{ij})$};
\node () at (1,2.5) {\footnotesize $1$};
\node () at (4,2.5) {\footnotesize $i$};
\node () at (8,2.5) {\footnotesize $j$};
\node () at (11,2.5) {\footnotesize $n$};
\end{scope}
\end{tikzpicture}
\caption{The partition $\ul\ve_{ij} = \ul c_{ij}^{-1}\ul e\ \ul c_{ij}$ for $\oijn$.}
\label{fig:words}
\end{center}
\end{figure}

\begin{lemma}\label{lem:ijw}
For $w\in\XSn^*$, and for distinct $i,j\in\bn$, we have
\[
w^{-1} \ve_{ij} w \sim \ve_{i\ul w,j\ul w} \AND w^{-1} \tau_{ij} w \sim \tau_{i\ul w,j\ul w}.
\]
\end{lemma}

\pf
By (both parts of) Lemma \ref{lem:FT}, it suffices to show that 
\[
\ul w^{-1} \ul\ve_{ij} \ul w = \ul\ve_{i\ul w,j\ul w} \AND \ul w^{-1} \ul\tau_{ij} \ul w = \ul\tau_{i\ul w,j\ul w}.
\]
By induction, it suffices to do this in the case where $w = s_k$ is a single letter from~$\XSn$.  For this we need to consider separate cases depending on the possible equalities between the elements of the sets $\{i,j\}$ and $\{k,k+1\}$.  Specifically, one can check diagrammatically that
\[
\ul s_k\ul \tau_{ij}\ul s_k = 
\begin{cases}
\ul \tau_{i+1,j} &\text{if $i=k$ and $j\not=k+1$,} \\
\ul \tau_{i-1,j} &\text{if $i=k+1$ and $j\not=k$,} \\
\ul \tau_{i,j+1} &\text{if $j=k$ and $i\not=k+1$,} \\
\ul \tau_{i,j-1} &\text{if $j=k+1$ and $i\not=k$,} \\
\ul \tau_{ji} &\text{if $\{i,j\}=\{k,k+1\}$,} \\
\ul \tau_{ij} &\text{if $\{i,j\}\cap\{k,k+1\}=\es$,}
\end{cases}
\]
with analogous formulae for $\ul s_k\ul \ve_{ij}\ul s_k$.  
\epf

We now also want to find a copy of the presentation $\Sgp\pres XR$ for $\Sing(\Pnfd)$ from Theorem~\ref{thm:SingPnfd} within $\Mon\pres YQ$.  For $w\in X^+$, we write $\ol w = w\phi \in\Sing(\Pnfd)$.  We also define a morphism
\[
\nu:X^+\to Y^*:w\mt\wh w \BY \wh e_{ij} = \ve_{ij} \ANd \wh t_{ij} = \tau_{ij}.
\]
It follows from \eqref{eq:olul} that $\nu\psi = \phi:X^+\to\Sing(\Pnfd)$, meaning that $\ul{\ \wh w\ } = \ol w$ for all $w\in X^+$.

\begin{lemma}\label{lem:approxsim}
If $u,v\in X^+$ are such that $\ol u = \ol v$, then $\wh u \sim \wh v$.
\end{lemma}

\pf
By Theorem \ref{thm:SingPnfd}, $\ol u=\ol v$ implies that $u$ can be transformed into $v$ by a sequence of applications of relations from $R$.  Thus, it suffices to prove the lemma in the case that $(u,v)$ is itself a relation from $R$.  If $(u,v)$ is any of \eqref{TT1}--\eqref{TT6}, then we note that $\wh u$ and $\wh v$ both belong to~$\XTn^*$, so since $\ul{\ \wh u\ } = \ol u = \ol v = \ul{\ \wh v\ }$, we obtain $\wh u \sim \wh v$ from part \ref{FT2} of Lemma \ref{lem:FT}.  Part \ref{FT1} of the same lemma deals with relations \eqref{EE1}--\eqref{EE3}.  This leaves us to deal with \eqref{ET1}--\eqref{ET4}.

\pfitem{\eqref{ET1}}  Here we have $\wh u = \ve_{ij}\tau_{ij}$ and $\wh v = \tau_{ij}$ (for distinct $i,j\in\bn$), so we need to show that $\ve_{ij}\tau_{ij} \sim \tau_{ij}$ .  To do so, let $w\in\XSn^*$ be such that $(1,2)\ul w = (i,j)$.  (Here for $s,t\in\bn$ and $\pi\in\S_n$ we write $(s,t)\pi=(s\pi,t\pi)$, and similarly allow permutations to act on other tuples.)  Then using both parts of Lemma \ref{lem:ijw} we have
\[
\ve_{ij} = \ve_{1\ul w,2\ul w} \sim w^{-1} \ve_{12} w = w^{-1}ew \ANDSIM \tau_{ij} \sim w^{-1}tw.
\]
It then follows that $\ve_{ij}\tau_{ij} \sim w^{-1}ew w^{-1}tw \Psim1 w^{-1}etw \Psim4 w^{-1}tw \sim \tau_{ij}$.

\pfitem{\eqref{ET2}}  This time we fix $w\in\XSn^*$ with $(1,2,3)\ul w=(i,j,k)$, and we have
\[
\wh u = \ve_{jk}\tau_{ij} \sim w^{-1}\ve_{23}\tau_{12}w = w^{-1} s_1s_2es_2s_1tw \ANDSIM \wh v \sim w^{-1}ts_2es_2w,
\]
so we need to show that $s_1s_2es_2s_1t  \sim ts_2es_2$.  For this we have
\[
s_1s_2es_2s_1t \Psim4 s_1s_2es_2t \Psim{11} s_1ts_2es_2 \Psim4 ts_2es_2.
\]

\pfitem{\eqref{ET3}}  Taking $w\in\XSn^*$ with $(1,2,3,4)\ul w = (i,j,k,l)$, we have $\tau_{ij} \sim w^{-1}\tau_{12}w = w^{-1}tw$ and
\[
\ve_{kl} \sim w^{-1}\ve_{34}w \sim w^{-1}(s_2s_1s_3s_2)e(s_2s_3s_1s_2)w \Psim{2} w^{-1}(s_2s_3s_1s_2)e(s_2s_3s_1s_2)w,
\]
so this time $\wh u\sim \wh v$ follows from \eqref{P14}, which says that $(s_2s_3s_1s_2)e(s_2s_3s_1s_2) t \sim t(s_2s_3s_1s_2)e(s_2s_3s_1s_2)$.

\pfitem{\eqref{ET4}}  The proof here is almost identical to that of \eqref{ET1}.
\epf

The next lemma will be used to show that any word over $Y$ with at least one occurrence of~$e$ or~$t$ is $\sim$-equivalent to an element of $\im(\nu)$.

\begin{lemma}\label{lem:ests}
If $i,j\in\bn$ are distinct, and if $1\leq k<n$, then
\ben
\item \label{ests1} $\tau_{ij}s_k \sim u$ and $s_k\tau_{ij}\sim v$ for some $u,v\in\im(\nu)$,
\item \label{ests2} $\ve_{ij}s_k \sim u$ and $s_k\ve_{ij}\sim v$ for some $u,v\in\im(\nu)$.
\een
\end{lemma}

\pf
\firstpfitem{\ref{ests1}}  Since $\ul\tau_{ij}\ul s_k \in \Sing(\T_n)$ we have $\ul\tau_{ij}\ul s_k = \ul\tau_{p_1q_1}\cdots\ul\tau_{p_lq_l}$ for suitable $p_t,q_t\in\bn$.  It then follows from Lemma \ref{lem:FT}\ref{FT2} that $\tau_{ij}s_k \sim \tau_{p_1q_1}\cdots\tau_{p_lq_l}$.
The proof for $s_k\tau_{ij}$ is essentially identical.

\pfitem{\ref{ests2}}  By part \ref{ests1} we have $s_k\tau_{ij} \sim w$ for some $w\in\im(\nu)$.  Since $\ve_{ij}\sim\tau_{ij}\ve_{ij}$, as shown in the proof of Lemma \ref{lem:approxsim}, it follows that $s_k\ve_{ij} \sim s_k\tau_{ij}\ve_{ij} \sim w\ve_{ij}$, so we can take $v = w\ve_{ij}$.

We also have $\ve_{ij}s_k \Psim1 s_ks_k\ve_{ij}s_k \sim s_k\ve_{i\ul s_k,j\ul s_k}$ (using Lemma \ref{lem:ijw} in the last step), and by the previous paragraph, this is $\sim$-equivalent to a word in $\im(\nu)$.
\epf

In the next proof, we write $\la \Si\ra$ for the subsemigroup of $Y^*$ generated by a subset $\Si\sub Y^*$.

\begin{lemma}\label{lem:imnu}
If $u\in Y^* \sm \XSn^*$, then $u\sim v$ for some $v\in\im(\nu)$.
\end{lemma}

\pf
Let $\Si = X\nu = \set{\ve_{ij} = \ve_{ji}}{\oijn} \cup \set{\tau_{ij}, \tau_{ji}}{\oijn}$, noting that $\im(\nu) = \la\Si\ra$.  Since $e=\ve_{12}$ and $t=\tau_{12}$ both belong to $\Si$, we have $Y^*\sm \XSn^* \sub \la\Si\cup\XSn\ra$.  Thus, we can prove the lemma by showing that any element~$u$ of $\la\Si\cup\XSn\ra\sm\XSn^*$ is $\sim$-equivalent to an element of $\la\Si\ra$.  For this, we write $u = z_1\cdots z_k$, where each $z_i\in\Si\cup\XSn$, and we use Lemma~\ref{lem:ests} and a simple induction on the number of terms $z_i$ that belong to $\XSn$.
\epf

We can now tie together the loose ends.

\pf[\bf Proof of Theorem \ref{thm:Pnfd}.]
It remains to show that $\ker(\psi)\sub{\sim}$, so let $(u,u')\in\ker(\psi)$, and write $a = \ul u = \ul u'$.  If $a\in\S_n$ then we must have $u,u'\in\XSn^*$ (as $\Sing(\Pnfd) = \Pnfd\sm\S_n$ is an ideal of $\Pnfd$), and then $u\sim u'$ by either part of Lemma \ref{lem:FT}.  So we now assume that $a\not\in\S_n$, in which case $u,u'\not\in\XSn^*$, so it follows from Lemma \ref{lem:imnu} that $u\sim v$ and $u'\sim v'$ for some $v,v'\in\im(\nu)$.  By definition, this means that $v = \wh{w}$ and $v' = \wh{w}'$ for some $w,w'\in X^+$.  Since ${\sim}\sub\ker(\psi)$ by Lemma \ref{lem:sub}, we also have
\[
\ol{w} = \ul{\ \wh{w}\ } = \ul v = \ul u = \ul u' = \ul v' = \ul{\ \wh{w}'\ } = \ol{w}',
\]
so it follows from Lemma \ref{lem:approxsim} that $\wh{w}\sim\wh{w}'$.  Thus, $u\sim v = \wh{w}\sim\wh{w}' = v' \sim u'$.
\epf

\section{\boldmath Presentation for $\PPnfd$}\label{sect:PPnfd}

In this section we prove Theorem \ref{thm:PPnfd}, which gives a presentation for the planar full-domain partition monoid~$\PPnfd$.  The proof is much more involved than that of Theorems \ref{thm:SingPnfd} and \ref{thm:Pnfd}.  It utilises Theorem \ref{thm:SU1}, and a certain strong (right) action pair in~$\PPnfd$.  One might wonder if this is simply the `planar analogue' of the action pair $(\E_n,\T_n)$ in the monoid $\Pnfd = \T_n\E_n$.  Specifically, one might wonder if $(\sP\E_n,\sP\T_n)$ is an action pair in~$\PPnfd$, and/or if we have the product decomposition $\PPnfd = \sP\T_n\cdot\sP\E_n$.  Both of these turn out not to be the case.  Indeed, while the pair $(\sP\E_n,\sP\T_n)$ satisfies \ref{A2} (by virtue of $\sP\E_n$ and~$\sP\T_n$ being submonoids of $\E_n$ and $\T_n$, which do form an action pair in $\Pnfd$), it does not satisfy \ref{A1} for $n\geq3$.  For example, consider $f = \custpartn{1,2,3}{1,2,3}{\stline11\stline21\stline33}\in\sP\T_3$ and $u = \custpartn{1,2,3}{1,2,3}{\stline11\stline22\stline33\uuline23\ddline23}\in\sP\E_3$.  Then $uf = \custpartn{1,2,3}{1,2,3}{\stline11\uuline13\stline33\darc13}$ is not equal to~$gv$ for any $g\in\sP\T_3$ and $v\in\sP\E_3$.  Another complication comes from the fact that $\PPnfd$ is not an \emph{Ehresmann-submonoid} of $\Pnfd$; it is closed under $D$, but not under $R$.  For example, with~$u$ and~$f$ as above, we have $uf\in\PP_3^\fd$, but $R(uf) = \custpartn{1,2,3}{1,2,3}{\stline11\stline22\stline33\uarc13\darc13} \not\in \PP_3^\fd$.

On the other hand, although the submonoids $\sP\E_n$ and $\sP\T_n$ do not lead to a product decomposition for $\PPnfd$, their union forms a generating set for $\PPnfd$.  Indeed, Theorem \ref{thm:PPnfd} claims (among other things) that $\PPnfd$ is generated by the elements
\[
\ol f_i = \custpartn{1,3,4,5,6,8}{1,3,4,5,6,8}{\stline11\stline33\stline44\stline54\stline66\stline88\udotted13\udotted68\ldotted13\ldotted68 \vertlab11\vertlab4i\vertlab8n}
\COMMa
\ol g_i = \custpartn{1,3,4,5,6,8}{1,3,4,5,6,8}{\stline11\stline33\stline45\stline55\stline66\stline88\udotted13\udotted68\ldotted13\ldotted68 \vertlab11\vertlab4i\vertlab8n}
\ANd
\ol h_i = \custpartn{1,3,4,5,6,8}{1,3,4,5,6,8}{\stline11\stline33\stline44\stline55\uuline45\ddline45\stline66\stline88\udotted13\udotted68\ldotted13\ldotted68 \vertlab11\vertlab4i\vertlab8n}
\qquad\text{for $1\leq i<n$,}
\]
and we have $\ol f_i,\ol g_i\in\sP\T_n$ and $\ol h_i\in\sP\E_n$.  One might wonder how such elements could generate partitions from $\PPnfd$ with nested lower blocks, such as $a = \custpartn{1,2,3,4,5}{1,2,3,4,5}{\draw(1,2)--(5,2);\stline11\darcx25{.6}\darcx34{.3}} \in \PP_5^\fd$.  But for this example we have
\[
\begin{tikzpicture}[scale=.35]
\begin{scope}[shift={(0,-3)}]	
\uvs{1,...,5}
\lvs{1,...,5}
\draw(1,2)--(5,2);
\stline11
\darcx25{.6}
\darcx34{.3}
\node () at (7,1) {$=$};
\node () at (-1,1) {$=$};
\node () at (-3,1) {$a$};
\end{scope}
\begin{scope}[shift={(8,0)}]	
\uvs{1,...,5}
\lvs{1,...,5}
\draw(1,2)--(5,2);
\stline11
\end{scope}
\begin{scope}[shift={(8,-2)}]	
\uvs{1,...,5}
\lvs{1,...,5}
\draw(3,2)--(4,2);
\draw(3,0)--(4,0);
\stline11
\stline22
\stline33
\stline44
\stline55
\end{scope}
\begin{scope}[shift={(8,-4)}]	
\uvs{1,...,5}
\lvs{1,...,5}
\stline11
\stline22
\stline32
\stline45
\stline55
\end{scope}
\begin{scope}[shift={(8,-6)}]	
\uvs{1,...,5}
\lvs{1,...,5}
\draw(3,2)--(4,2);
\draw(3,0)--(4,0);
\stline11
\stline22
\stline33
\stline44
\stline55
\node () at (6,3) {,};
\end{scope}
\end{tikzpicture}
\]
which leads to the factorisation $a = \ol f_4\ol f_3\ol f_2\ol f_1 \cdot \ol h_3 \cdot \ol f_2\ol g_4 \cdot \ol h_3$.

As we will see, however, the planar submonoid $\sP\T_n$ \emph{is} part of an action pair in~$\PPnfd$.  Before we discuss this, we note that in fact $\sP\T_n = \O_n$ is (an isomorphic copy of) the monoid of order-preserving transformations of $\bn$, as introduced in Section \ref{subsect:S}.  The pair $(\D_n,\O_n)$ we use here will involve a new submonoid $\D_n$.  The definition of $\D_n$ is somewhat involved, and will be given in Section \ref{subsect:Dn}, where we also prove a number of useful technical results, including a natural generating set.  Although $\D_n$ is not a semilattice, it is still a band, indeed a \emph{right regular} band, meaning that it satisfies the identity $xyx=yx$; equivalently, it is $\L$-trivial.  We then give a presentation for $\D_n$ in Section \ref{subsect:presDn}, and verify that $(\D_n,\O_n)$ is indeed an action pair with $\O_n\D_n=\PPnfd$ in Section \ref{subsect:DnOn}.  We then apply Theorem \ref{thm:SU1} to obtain a presentation for $\PPnfd$ in Section \ref{subsect:presPPnfd} in terms of a large generating set, and then use this to prove Theorem \ref{thm:PPnfd}.  Finally, we make some observations on the structure of $\PPnfd$ in Section~\ref{subsect:structurePPnfd}, showing in particular that it is an example of a \emph{grrac monoid}, in the sense of Branco, Gomes and Gould \cite{BGG2010}.

\subsection[The monoid $\D_n$]{\boldmath The monoid $\D_n$}\label{subsect:Dn}

We begin with some terminology and notation concerning planar partitions, and their (co)kernels.
Consider non-empty subsets $A=\{i_1<\cdots<i_k\}$ and $B=\{j_1<\cdots<j_l\}$ of $\bn$.  We write $A<B$ if $i_k<j_1$.  If $A<B$ or $B<A$, then we say that $A$ and $B$ are \emph{separated}.  We write $A\prec B$, and say that $A$ is \emph{nested} by $B$, if we have $j_t<i_1<\cdots<i_k<j_{t+1}$ for some $1\leq t<l$.  If $A\prec B$ or $B\prec A$, we say that $A$ and $B$ are \emph{nested}.  The immediate relevance of these concepts is contained in the following result concerning the planar partition monoid $\PP_n$, which is part of \cite[Lemma 7.1]{EMRT2018}.  In what follows, we will often use this lemma without explicit reference.

\begin{lemma}\label{lem:ns}
If $a = \begin{partn}{6} A_1&\cdots&A_r&C_1&\cdots&C_s\\ \hhline{~|~|~|-|-|-} B_1&\cdots&B_r&D_1&\cdots&D_t\end{partn}\in\PP_n$, with $\min(A_1)<\cdots<\min(A_r)$, then
\ben
\item $A_1<\cdots<A_r$ and $B_1<\cdots<B_r$,
\item for all $1\leq i<j\leq s$, $C_i$ and $C_j$ are either nested or separated,
\item for all $1\leq i<j\leq t$, $D_i$ and $D_j$ are either nested or separated,
\item for all $1\leq i\leq r$ and $1\leq j\leq s$, either $A_i$ and $C_j$ are separated or else $C_j$ is nested by $A_i$, 
\item for all $1\leq i\leq r$ and $1\leq j\leq t$, either $B_i$ and $D_j$ are separated or else $D_j$ is nested by $B_i$.  
\een
\end{lemma}

An equivalence relation $\eta\in\Eq_n$ is called \emph{planar} if every pair of $\eta$-classes is either separated or nested.  It follows quickly from Lemma \ref{lem:ns} that these are precisely the kernels and cokernels of planar partitions.  They are also the cokernels of planar \emph{full-domain} partitions, as we will see below.  It also follows from Lemma \ref{lem:ns} that an equivalence $\eta\in\Eq_n$ is the \emph{kernel} of a planar full-domain partition if and only if it is \emph{convex}, in the sense that every $\eta$-class is an interval.  This is equivalent to every pair of $\eta$-classes being separated.  

We write $\PEq_n$ for the set of all planar equivalences on $\bn$.  Note that $\PEq_n$ is not a $\vee$-submonoid of $\Eq_n$ (for $n\geq4$), as $\PEq_n$ contains the generating set $\set{\eta_{ij}}{\oijn}$ for~$\Eq_n$.
On the other hand, the set of all convex equivalences of $\bn$ is a submonoid, and is isomorphic to $\PE_n = \E_n\cap\PP_n$.  

Consider a planar equivalence $\eta\in\PEq_n$.  Let the un-nested $\eta$-classes be ${B_1<\cdots<B_r}$, and let the remaining $\eta$-classes be $C_1,\ldots,C_s$.  (So each $C_i$ is nested by some $B_j$, but there may also be nesting among the $C_i$ themselves.)  For each $1\leq i\leq r$, let $x_i = \min(B_i)$ and $y_i=\max(B_i)$, and note that
\[
x_1 = 1 \COMMA x_i = y_{i-1}+1 \text{ for $2\leq i\leq r$} \AND y_r = n.
\]
We then define the intervals $A_i = [x_i,y_i] = \{x_i,x_i+1,\ldots,y_i\}$, and the partition
\[
d_\eta = \begin{partn}{6} A_1&\cdots&A_r& \multicolumn{3}{c}{} \\ \hhline{~|~|~|-|-|-} B_1&\cdots&B_r&C_1&\cdots&C_s\end{partn} \in \PPnfd.
\]
It is easy to see that $d_\eta$ is an idempotent.  
As an example, if $\eta\in\PEq_8$ is such that
\begin{equation}\label{eq:D8}
{\bf8}/\eta = \big\{ \{1,5,6\} , \{2,3\} , \{4\} , \{7,8\} \big\}, \qquad\text{then}\qquad d_\eta = \custpartn{1,...,8}{1,...,8}{\stlines{1/1,6/6,7/7,8/8}\uuline16\uuline78\ddline23\darc15\darc56\darc78}.
\end{equation}
It is immediate from the construction that
\[
\coker(d_\eta) = \eta \qquad\text{for all $\eta\in\PEq_n$.}
\]
The kernel of $d_\eta$ is not necessarily $\eta$, though.  We will denote this kernel by $\wh\eta = \ker(d_\eta)$.  With the above notation, we have $\bn/\wh\eta = \{A_1,\ldots,A_r\} = \big\{[x_1,y_1],\ldots,[x_r,y_r]\big\}$.

We now define the set
\[
\D_n = \set{d_\eta}{\eta\in\PEq_n}.
\]
We will show that $\D_n$ is a submonoid of $\PPnfd$ by a somewhat indirect method, which will also allow us to establish a natural generating set.  This will in turn feed into a presentation in Section~\ref{subsect:presDn}.  For now it is at least clear that $\D_n$ is a subset of~$\PPnfd$, and contains the identity element $1=d_{\De_\bn}$.  The next lemma gives an alternative characterisation of the elements of~$\D_n$.

\begin{lemma}\label{lem:minmax}
Let $a\in\PPnfd$.  Then $a\in\D_n$ if and only if every transversal $A\cup B'$ of $a$ satisfies $\min(A)=\min(B)$ and $\max(A)=\max(B)$, in which case we have $a = d_\eta$, where $\eta=\coker(a)$.
\end{lemma}

\pf
The forward implication follows immediately from the definition of the elements $d_\eta$.  Conversely, suppose $a = \begin{partn}{6} A_1&\cdots&A_r& \multicolumn{3}{c}{} \\ \hhline{~|~|~|-|-|-} B_1&\cdots&B_r&C_1&\cdots&C_s\end{partn}\in\PPnfd$ satisfies the stated condition on transversals, assuming that $A_1<\cdots<A_r$, and let $\eta=\coker(a)$.  We see then that $B_1,\ldots,B_r$ are precisely the un-nested $\eta$-classes, and it then follows by definition that $a=d_\eta\in\D_n$.
\epf

For $\oijn$, we define the partition
\[
\ol h_{ij} = d_{\eta_{ij}} = \custpartn{1,3,4,5,8,9,10,12}{1,3,4,5,8,9,10,12}{\stline11\stline33\stline44\stline99\stline{10}{10}\stline{12}{12} \uuline4{5.75}\uuline{7.25}9 \udotted13\udotted{5.5}{7.5}\udotted{10}{12}\ldotted13\ldotted{5}{8}\ldotted{10}{12} \vertlab11\vertlab4i\vertlab9j\vertlab{12}n}  \in \D_n.
\]
Note that $[i,j]\cup\{i',j'\}$ is a block of $\ol h_{ij}$, and so is each of $\{i+1\}',\ldots,\{j-1\}'$.  As a concrete example, in $\D_{10}$ we have
\[
\ol h_{38} =  \custpartn{1,...,10}{1,...,10}{\stlines{1/1,2/2,3/3,8/8,9/9,10/10} \uuline38 }.
\]
Our first main aim is to show that $\D_n$ is a monoid, and that these elements $\ol h_{ij}$ generate it.  Working towards this, we start by introducing some notation that will be used frequently.

For a subset $A = \{x_1<\cdots<x_k\} \sub \bn$, and for $1\leq i\leq k$, we define the \emph{$A$-successor of $x_i$} to be
\[
\suc_A(x_i) = \begin{cases}
x_{i+1} &\text{if $i<k$,}\\
x_k &\text{if $i=k$.}
\end{cases}
\]
For an equivalence $\eta\in\Eq_n$ and for $x\in\bn$, we define the \emph{$\eta$-successor of $x$} to be
\[
\suc_\eta(x) = \suc_{[x]_\eta}(x).
\]
It is clear that for $\eta,\mu\in\Eq_n$ we have
\[
\eta=\mu \IFF \suc_\eta(x) = \suc_\mu(x) \text{ for all $x\in\bn$.}
\]
The next lemma uses successors to describe the situation when we left-multiply an element of~$\D_n$ with a partition of the form $\ol h_{ij}$.  The statement also uses the equivalences $\wh\eta = \ker(d_\eta)$.

\begin{lemma}\label{lem:suc}
Let $\eta\in\PEq_n$ and $\oijn$, put $p = \max[i]_{\wh\eta}$ and $q = \min[j]_{\wh\eta}$, and let $\mu = \eta\vee\eta_{pq}$.  Then
\ben
\item \label{suc1} $\mu\in\PEq_n$, and $\ol h_{ij}d_\eta = d_\mu$,
\item \label{suc2} if $\mu\not=\eta$, then for any $x\in\bn$ we have
\[
\suc_\mu(x) = \begin{cases}
q &\text{if $x=p$,}\\
\suc_\eta(x) &\text{otherwise.}
\end{cases}
\]
\een
\end{lemma}

\pf
\firstpfitem{\ref{suc1}}
Write $d_\eta = \begin{partn}{6} A_1&\cdots&A_r& \multicolumn{3}{c}{} \\ \hhline{~|~|~|-|-|-} B_1&\cdots&B_r&C_1&\cdots&C_s\end{partn}$, where $A_1<\cdots<A_r$, and assume that $i\in A_k$ and $j\in A_l$.  Since $i<j$, we have $k\leq l$.  Using Lemma \ref{lem:minmax}, we have
\[
p = \max[i]_{\wh\eta} = \max(A_k) = \max(B_k), \ANDSIM q = \min(B_l).
\]
Consequently, $B_k\cup B_l$ is a $\mu = \eta\vee\eta_{pq}$-class, and the remaining $\mu$-classes are all $\eta$-classes.  Since~$B_k$ and~$B_l$ are un-nested $\eta$-classes, it is easy to see that $\mu$ is planar.

Now consider the product graph $\Pi(\ol h_{ij},d_\eta)$, and note that the only non-trivial connection in the middle row coming from~$\ol h_{ij}$ is $\{i'',j''\}$.  We now consider two cases.

\begin{figure}[ht]
\begin{center}
\begin{tikzpicture}[scale=.44]
\begin{scope}[shift={(0,0)}]	
\hlines06
\darcs{0/2,2/4,4/6}
\hlines7{14}
\darcs{7/9,9/10,10/11,11/14}
\hlines{15}{18}
\darcs{15/18}
\hlines{19}{27}
\darcs{19/21,21/24,24/27}
\hlines{28}{34}
\darcs{28/30,30/32,32/33,33/34}
\draw(0,1)node[left]{\footnotesize $d_\eta$};
\draw[|-|] (7,-.5)--(14,-.5);
\node () at (10.5,-1) {\footnotesize $A_k=A_l$};
\end{scope}
\begin{scope}[shift={(0,2)}]	
\vlines06
\vlines7{9.5}
\vlines{12}{14}
\vlines{15}{18}
\vlines{19}{27}
\vlines{28}{34}
\hlinesx{9.5}{12}{very thick,red}
\uvs{7,9.5,12,14}
\vertlab7q
\vertlab{9.5}i
\vertlab{12}j
\vertlab{14}p
\draw(0,1)node[left]{\footnotesize $\ol h_{ij}$};
\end{scope}
\begin{scope}[shift={(0,-8)}]	
\hlines06
\darcs{0/2,2/4,4/6}
\hlines7{14}
\darcs{7/9,9/10,10/11,11/14}
\hlines{15}{18}
\darcs{15/18}
\hlines{19}{27}
\darcs{19/21,21/24,24/27}
\hlines{28}{34}
\darcs{28/30,30/32,32/33,33/34}
\draw(0,1)node[left]{\footnotesize $d_\eta$};
\vertlab{14}p
\vertlab{19}q
\draw[|-|] (7,-.5)--(14,-.5);
\node () at (10.5,-1) {\footnotesize $A_k$};
\draw[|-|] (19,-.5)--(27,-.5);
\node () at (23,-1) {\footnotesize $A_l$};
\end{scope}
\begin{scope}[shift={(0,-6)}]	
\vlines06
\vlines7{9.5}
\vlines{24}{27}
\vlines{28}{34}
\hlinesx{9.5}{24}{very thick,red}
\uvs{9.5,24}
\lvs{19,14}
\vertlab{9.5}i
\vertlab{24}j
\draw(0,1)node[left]{\footnotesize $\ol h_{ij}$};
\end{scope}
\end{tikzpicture}
\caption{Schematic illustrations of products $\ol h_{ij}d_\eta$ from the proof of Lemma \ref{lem:suc}, in the cases ${k=l}$~(top) and $k<l$ (bottom).  Only the transversals of $d_\eta$ are shown.  See the proof for more details.}
\label{fig:suc}
\end{center}
\end{figure}
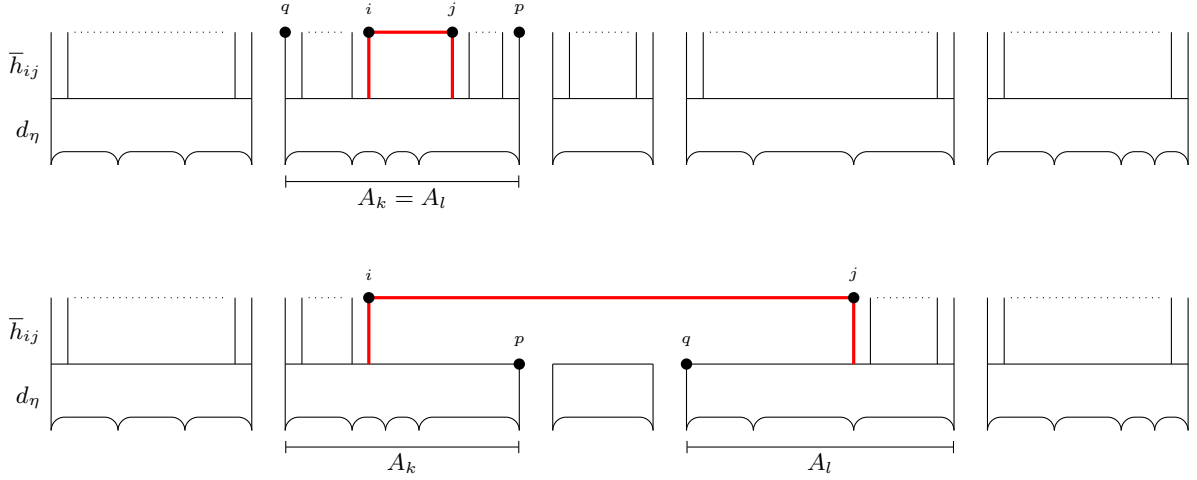

\pfcase1  Suppose first that $k=l$.  Here it is easy to see that $\ol h_{ij}d_\eta = d_\eta$.  Moreover, in this case~$p$ and~$q$ belong to the same $\eta$-class (i.e.~$B_k$), and so $\mu = \eta\vee\eta_{pq} = \eta$.  See Figure \ref{fig:suc} (top) for a schematic illustration of this case.

\pfcase2  Now we assume that $k<l$.  Again considering the product graph $\Pi(\ol h_{ij},d_\eta)$, the connection $\{i'',j''\}$ from $\ol h_{ij}$ joins the blocks $A_k''\cup B_k'$ and $A_l''\cup B_l'$, and we see that
\[
\ol h_{ij} d_\eta = \begin{partn}{13} A_1&\cdots&A_{k-1}&A_k\cup\cdots\cup A_l&A_{l+1}&\cdots&A_r& \multicolumn{6}{c}{} \\ \hhline{~|~|~|~|~|~|~|-|-|-|-|-|-} B_1&\cdots&B_{k-1}&B_k\cup B_l&B_{l+1}&\cdots&B_r&B_{k+1}&\cdots&B_{l-1}&C_1&\cdots&C_s\end{partn}.
\]
See Figure \ref{fig:suc} (bottom).  It is easy to see that this partition satisfies the conditions of Lemma~\ref{lem:minmax}, so we have $\ol h_{ij} d_\eta = d_\nu$ by that lemma, where $\nu = \coker(\ol h_{ij}d_\eta)$.  But $\nu$ is obtained from $\eta$ by joining the $\eta$-classes $B_k$ and $B_l$ into a single $\nu$-class (and leaving the other $\eta$-classes intact).  Since $p\in B_k$ and $q\in B_l$, this means that $\nu$ is exactly $\eta\vee\eta_{pq} = \mu$.

\pfitem{\ref{suc2}}  We have already noted that when $\mu$ is not equal to $\eta$, it is obtained by replacing the two $\eta$-classes $B_k$ and $B_l$ by the single $\mu$-class $B_k\cup B_l$.  Thus, the only $x\in\bn$ for which $\suc_\eta(x) \not= \suc_\mu(x)$ is $x= \max(B_k) = p$, and hence $\suc_\mu(p) = \min(B_l) = q$.
\epf

It is convenient at this point to define an abstract alphabet
\[
\XDn = \set{h_{ij}}{\oijn}.
\]
For a word $w = h_{i_1j_1}\cdots h_{i_kj_k}\in\XDn^*$ we write $\ol w = \ol h_{i_1j_1}\cdots \ol h_{i_kj_k}$.  In what follows, we will be interested in a specific family of words over $\D_n$.  In order to define these, we will first extend the~$h_{ij}$ notation, by defining
\[
h_{ii} = \io \AND \ol h_{ii} = 1 \qquad\text{for all $1\leq i\leq n$.}
\]
For a planar equivalence $\eta\in\PEq_n$ we define the word
\[
w_\eta = h_{1k_1}\cdots h_{nk_n}, \WHERE k_x = \suc_\eta(x)\text{ for all $x\in\bn$.}
\]
For example, with $\eta\in\D_8$ as in \eqref{eq:D8}, we have
\begin{equation}\label{eq:D82}
w_\eta = h_{15}h_{23}h_{33}h_{44}h_{56}h_{66}h_{78}h_{88} = h_{15}h_{23}h_{56}h_{78}.
\end{equation}

\begin{lemma}\label{lem:weta}
For any $\eta\in\PEq_n$ we have $d_\eta = \ol w_\eta$.
\end{lemma}

\pf
We use descending induction on $m = |\bn/\eta|$.  When $m=n$, we have $\eta = \De_\bn$, in which case $d_\eta = 1$ and $w_\eta = h_{11}\cdots h_{nn} = \io$, so the result is clear.  We now assume that $m<n$.  Let
\[
i = {\min}\bigset{x\in\bn}{[x]_\eta\not=\{x\}}.
\]
Since each of $1,\ldots,i-1$ belongs to a singleton $\eta$-class, we have $h_{1k_1}=\cdots=h_{i-1,k_{i-1}}=\io$, and so
\[
w_\eta = h_{ik_i}\cdots h_{nk_n}, \WHERE k_x = \suc_\eta(x) \text{ for each $i\leq x\leq n$.}
\]
Let $\mu\in\PEq_n$ be obtained from $\eta$ by replacing the class $[i]_\eta$ by $\{i\}$ and $[i]_\eta\sm\{i\}$.  We then have
\begin{align*}
d_\mu &= \ol w_\mu &&\text{by induction, as $|\bn/\mu|=m+1$}\\
&= \ol h_{i+1,l_{i+1}}\cdots \ol h_{nl_n} &&\text{where each $l_i = \suc_\mu(x)$}\\
&= \ol h_{i+1,k_{i+1}}\cdots \ol h_{nk_n} &&\text{by Lemma \ref{lem:suc}\ref{suc2}, as $\eta=\mu\vee\eta_{ik_i}$.}
\end{align*}
It follows that
\[
\ol w_\eta = \ol h_{ik_i}\cdot\ol h_{i+1,k_{i+1}}\cdots h_{nk_n} = \ol h_{ik_i}d_\mu.
\]
By Lemma \ref{lem:suc}\ref{suc1}, we have $\ol h_{ik_i}d_\mu = d_\nu$, where $\nu = \mu\vee\eta_{pq}$ for 
\[
p = \max[i]_{\wh\mu} =  i \AND q = \min[k_i]_{\wh\mu} = k_i.
\]
But this says exactly that $\nu = \mu\vee\eta_{ik_i} = \eta$, so that $d_\eta = d_\nu = \ol h_{ik_i}d_\mu = \ol w_\eta$, as required.
\epf

As an example, consider again $\eta\in\PEq_8$ from \eqref{eq:D8}.  We then have the following factorisation, verifying that $d_\eta = \ol w_\eta$; cf~\eqref{eq:D82}.
\[
\begin{tikzpicture}[scale=.35]
\begin{scope}[shift={(0,-3)}]	
\uvs{1,...,8}
\lvs{1,...,8}
\stlines{1/1,6/6,7/7,8/8}\uuline16\uuline78\ddline23\darc15\darc56\darc78
\node () at (10,1) {$=$};
\node () at (-1,1) {$=$};
\node () at (-3,1) {$d_\eta$};
\end{scope}
\begin{scope}[shift={(11,0)}]	
\uvs{1,...,8}
\lvs{1,...,8}
\stlines{1/1,5/5,6/6,7/7,8/8}
\uuline15
\node()at(9.25,1){\footnotesize$\ol h_{15}$};
\end{scope}
\begin{scope}[shift={(11,-2)}]	
\uvs{1,...,8}
\lvs{1,...,8}
\stlines{1/1,2/2,3/3,4/4,5/5,6/6,7/7,8/8}
\uuline23
\node()at(9.25,1){\footnotesize$\ol h_{23}$};
\end{scope}
\begin{scope}[shift={(11,-4)}]	
\uvs{1,...,8}
\lvs{1,...,8}
\stlines{1/1,2/2,3/3,4/4,5/5,6/6,7/7,8/8}
\uuline56
\node()at(9.25,1){\footnotesize$\ol h_{56}$};
\end{scope}
\begin{scope}[shift={(11,-6)}]	
\uvs{1,...,8}
\lvs{1,...,8}
\stlines{1/1,2/2,3/3,4/4,5/5,6/6,7/7,8/8}
\uuline78
\node () at (6,0) {,};
\node()at(9.25,1){\footnotesize$\ol h_{78}$};
\end{scope}
\end{tikzpicture}
\]
For the next statement, recall that a band $S$ is \emph{right regular} if $xyx=yx$ for all $x,y\in S$.  This is equivalent to $S$ being $\L$-trivial, i.e.~satisfying $x\L y\implies x=y$.

\begin{prop}\label{prop:Dn}
$\D_n$ is a submonoid of $\PPnfd$, and is in fact a right regular band.  Moreover, we have
\[
\D_n = \pres{\ol h_{ij}}{\oijn}.
\]
\end{prop}

\pf
The inclusion $\D_n \sub \pres{\ol h_{ij}}{\oijn}$ follows from Lemma \ref{lem:weta}, and the backward inclusion from Lemma \ref{lem:suc}\ref{suc1} and a simple induction.  Thus $\D_n = \pres{\ol h_{ij}}{\oijn}$, which also of course shows that $\D_n$ is a submonoid of $\PPnfd$.  

To demonstrate $\L$-triviality, let $\eta,\mu\in\PEq_n$.  Then using the analogue of \eqref{eq:RinPn} for the $\L$-relation, we have
\[
d_\eta \L d_\mu \Implies \coker(d_\eta) = \coker(d_\mu) \Implies \eta = \mu \Implies d_\eta = d_\mu.  \qedhere
\]
\epf

\subsection[Presentation for $\D_n$]{\boldmath Presentation for $\D_n$}\label{subsect:presDn}

We continue to fix the alphabet
\[
\XDn = \set{h_{ij}}{\oijn}.
\]
By Proposition \ref{prop:Dn}, we have a surmorphism
\[
\de:\XDn^*\to\D_n:h_{ij}\mt\ol h_{ij}.
\]
We now let $\RDn$ be the set of the following relations over $\XDn$:
\begin{align}
\label{N1} h_{ij}h_{kl} &= h_{kl} &&\text{if $k\leq i<j\leq l$,}\\
\label{N2} h_{ij}h_{kl} &= h_{kl}h_{ij} &&\text{if $i<j\leq k<l$,}\\
\label{N3} h_{ij}h_{jk} = h_{ik}h_{ij} &= h_{ik}h_{jk} &&\text{if $i<j<k$.}
\end{align}
Note that \eqref{N1} includes as a special case the relations
\[
h_{ij}^2 = h_{ij} \qquad\text{for all $\oijn$.}
\]
Our main goal of this section is to prove the following:

\begin{thm}\label{thm:Dn}
The monoid $\D_n$ has presentation $\Mon\pres \XDn \RDn$ via $\de$.
\end{thm}

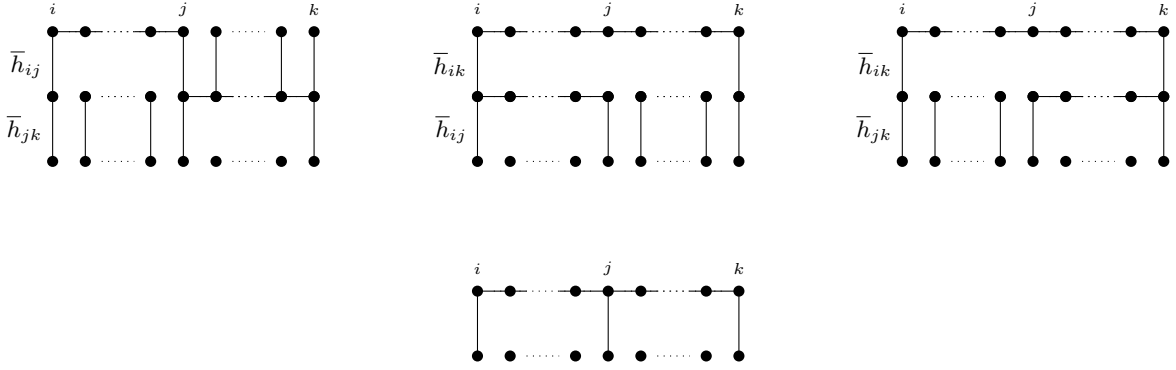
\begin{figure}[ht]
\begin{center}
\begin{tikzpicture}[scale=.43]
\begin{scope}[shift={(0,0)}]	
\uvs{1,2,4,5,6,8,9}
\lvs{1,2,4,5,6,8,9}
\stlines{1/1,5/5,6/6,8/8,9/9}
\udotted15
\uuline1{2.5}
\uuline{3.5}5
\vertlab1i
\vertlab5j
\vertlab9k
\udotted24
\ldotted24
\udotted68
\ldotted68
\draw(1,1)node[left]{\footnotesize $\ol h_{ij}$};
\end{scope}
\begin{scope}[shift={(0,-2)}]	
\uvs{1,2,4,5,6,8,9}
\lvs{1,2,4,5,6,8,9}
\stlines{1/1,2/2,4/4,5/5,9/9}
\udotted59
\uuline5{6.5}
\uuline{7.5}9
\udotted24
\ldotted24
\udotted68
\ldotted68
\draw(1,1)node[left]{\footnotesize $\ol h_{jk}$};
\end{scope}
\begin{scope}[shift={(13,0)}]	
\uvs{1,2,4,5,6,8,9}
\lvs{1,2,4,5,6,8,9}
\stlines{1/1,9/9}
\udotted15
\uuline1{2.5}
\uuline{3.5}5
\udotted59
\uuline5{6.5}
\uuline{7.5}9
\vertlab1i
\vertlab5j
\vertlab9k
\udotted24
\ldotted24
\udotted68
\ldotted68
\draw(1,1)node[left]{\footnotesize $\ol h_{ik}$};
\end{scope}
\begin{scope}[shift={(13,-2)}]	
\uvs{1,2,4,5,6,8,9}
\lvs{1,2,4,5,6,8,9}
\stlines{1/1,5/5,6/6,8/8,9/9}
\udotted15
\uuline1{2.5}
\uuline{3.5}5
\udotted24
\ldotted24
\udotted68
\ldotted68
\draw(1,1)node[left]{\footnotesize $\ol h_{ij}$};
\end{scope}
\begin{scope}[shift={(26,0)}]	
\uvs{1,2,4,5,6,8,9}
\lvs{1,2,4,5,6,8,9}
\stlines{1/1,9/9}
\udotted15
\uuline1{2.5}
\uuline{3.5}5
\udotted59
\uuline5{6.5}
\uuline{7.5}9
\vertlab1i
\vertlab5j
\vertlab9k
\udotted24
\ldotted24
\udotted68
\ldotted68
\draw(1,1)node[left]{\footnotesize $\ol h_{ik}$};
\end{scope}
\begin{scope}[shift={(26,-2)}]	
\uvs{1,2,4,5,6,8,9}
\lvs{1,2,4,5,6,8,9}
\stlines{1/1,2/2,4/4,5/5,9/9}
\udotted59
\uuline5{6.5}
\uuline{7.5}9
\udotted24
\ldotted24
\udotted68
\ldotted68
\draw(1,1)node[left]{\footnotesize $\ol h_{jk}$};
\end{scope}
\begin{scope}[shift={(13,-8)}]	
\vertlab1i
\vertlab5j
\vertlab9k
\uvs{1,2,4,5,6,8,9}
\lvs{1,2,4,5,6,8,9}
\stlines{1/1,5/5,9/9}
\udotted15
\uuline1{2.5}
\uuline{3.5}5
\udotted59
\uuline5{6.5}
\uuline{7.5}9
\ldotted24
\udotted68
\ldotted68
\end{scope}
\end{tikzpicture}
\caption{Diagrammatic verification of relation \eqref{N3}.
Top row: the product graphs $\Pi(\ol h_{ij},\ol h_{jk})$, $\Pi(\ol h_{ik},\ol h_{ij})$ and $\Pi(\ol h_{ik},\ol h_{jk})$.
Bottom row: the element $\ol h_{ij}\ol h_{jk} = \ol h_{ik}\ol h_{ij} = \ol h_{ik}\ol h_{jk}$.  In all diagrams we have pictured only vertices $i,\ldots,k$, and their (double-)dashed counterparts.}
\label{fig:N3}
\end{center}
\end{figure}

It is easy to check that $\RDn\sub\ker(\de)$.  For example, we verify relation \eqref{N3} in Figure \ref{fig:N3}.

As ever, the inclusion $\ker(\de)\sub\RDn^\sharp$ is harder to establish.
We begin by establishing some simple consequences of the relations \eqref{N1}--\eqref{N3}.  Recall the convention that $h_{ii} = \io$ for all $i$.  Throughout this section we write ${\sim}=\RDn^\sharp$.  There should be no confusion with our earlier use of this symbol.

\begin{lemma}\label{lem:N}
For $i,j,k,l\in\bn$ we have
\begin{align}
\label{N1'} h_{ij}h_{kl} &\sim h_{kl} &&\text{if $k\leq i\leq j\leq l$,}\\
\label{N2'} h_{ij}h_{kl} &\sim h_{kl}h_{ij} &&\text{if $i\leq j\leq k\leq l$,}\\
\label{N3'} h_{ij}h_{jk} &\sim h_{ik}h_{ij} \sim h_{ik}h_{jk} &&\text{if $i\leq j\leq k$,}\\
\label{N4} h_{ij}h_{kl} &\sim h_{lj}h_{kl} &&\text{if $k\leq i\leq l\leq j$,}\\
\label{N5} h_{ij}h_{kl} &\sim h_{ik}h_{kl} &&\text{if $i\leq k\leq j\leq l$.}
\end{align}
\end{lemma}

\pf
\firstpfitem{\eqref{N1'}--\eqref{N3'}}  Note that these are just \eqref{N1}--\eqref{N3} with all `$<$' signs replaced with `$\leq$'.  Thus, we just have to check that in cases with equalities of subscripts not covered in \eqref{N1}--\eqref{N3} we have $\sim$-equivalence of the stated words.  This is obvious for~\eqref{N1'} and~\eqref{N2'}.  For \eqref{N3'} it follows from the above-mentioned relations~$h_{ij}^2=h_{ij}$, which are part of \eqref{N1}.

\pfitem{\eqref{N4}}  Here we have $h_{ij}h_{kl} \Nsim1 h_{ij}h_{il}h_{kl} \Nsim3 h_{il}h_{lj}h_{kl} \Nsim2 h_{lj}h_{il}h_{kl} \Nsim1 h_{lj}h_{kl}$.

\pfitem{\eqref{N5}}  This is similar to \eqref{N4}, but a little easier.
\epf

We now prove a series of results concerning the words $w_\eta\in\XDn^*$ ($\eta\in\PEq_n$) defined above, leading up to Corollary \ref{cor:wweta}, which shows that every word over $\XDn$ is $\sim$-equivalent to some $w_\eta$.

So fix some $\eta\in\PEq_n$, and consider the word $w_\eta = h_{1k_1}\cdots h_{nk_n}$, where each $k_x=\suc_\eta(x)$, and let $p\in\bn$ be arbitrary.  Define the sequence 
\[
p_0 = p \COMMA p_1 = k_{p_0} \COMMA p_2 = k_{p_1} \COMMA \ldots,
\]
and suppose this terminates in $p_l = k_{p_l}$, where $p_0,p_1,\ldots,p_l$ are all distinct.  Also set $q = p_l$.

\begin{lemma}\label{lem:orbit}
With the above notation, and with $u = h_{pk_p}\cdots h_{qk_q}$, we have $u \sim h_{pq}u$.
\end{lemma}

\pf
We use induction on $l$.  If $l=0$ then $q=p$, so $h_{pq}=\io$, and the claim is trivial.  So now suppose $l\geq1$.  We then have
\begin{equation}\label{eq:orbit1}
u = h_{pk_p}\cdots h_{qk_q} = h_{p_0k_{p_0}} h_{p_0+1,k_{p_0+1}} \cdots h_{p_1-1,k_{p_1-1}} \cdot v,
\end{equation}
where $v = h_{p_1k_{p_1}}\cdots h_{qk_q}$.  By induction (considering that $p_1,p_2,\ldots,p_l=q$ is a shorter sequence of the relevant form), we have 
\begin{equation}\label{eq:orbit2}
v\sim h_{p_1q}v.
\end{equation}
We then calculate
\begin{align}
\nonumber u &\sim h_{p_0k_{p_0}} h_{p_0+1,k_{p_0+1}} \cdots h_{p_1-1,k_{p_1-1}} \cdot h_{p_1q}v &&\text{by \eqref{eq:orbit1} and \eqref{eq:orbit2}}\\
\label{eq:orbit3} &\sim h_{p_0k_{p_0}} h_{p_1q}\cdot h_{p_0+1,k_{p_0+1}} \cdots h_{p_1-1,k_{p_1-1}} \cdot v &&\text{by \eqref{N2}.}
\intertext{In this last step, relation \eqref{N2} applies because $k_i<p_1$ for all $p_0+1\leq i\leq p_1-1$, since the points $p_0+1,\ldots, p_1-1$ are nested by the block of $\eta$ containing $p_0$ and $p_1=k_{p_0} = \suc_\eta(p_0)$.  Keeping in mind that $k_{p_0}=p_1$, we then continue:}
\nonumber u &\sim h_{p_0p_1} h_{p_1q}\cdot h_{p_0+1,k_{p_0+1}} \cdots h_{p_1-1,k_{p_1-1}} \cdot v &&\text{by \eqref{eq:orbit3}}\\
\nonumber  &\sim h_{p_0q}h_{p_0p_1} \cdot h_{p_0+1,k_{p_0+1}} \cdots h_{p_1-1,k_{p_1-1}} \cdot v &&\text{by \eqref{N3'}}\\
\nonumber  &= h_{pq}h_{p_0k_{p_0}} \cdot h_{p_0+1,k_{p_0+1}} \cdots h_{p_1-1,k_{p_1-1}} \cdot v &&\text{as $p_0=p$ and $p_1=k_{p_0}$}\\
\nonumber &= h_{pq}u,
\end{align}
completing the proof.
\epf

Among other things, the next lemma identifies an important case of Lemma \ref{lem:orbit}.  To state the result, we first introduce some further terminology.  Consider again $\eta\in\PEq_n$, and write
\[
w_\eta = h_{1k_1}\cdots h_{nk_n}
\AND
d_\eta = \begin{partn}{6} A_1&\cdots&A_r& \multicolumn{3}{c}{} \\ \hhline{~|~|~|-|-|-} B_1&\cdots&B_r&C_1&\cdots&C_s\end{partn}.
\]
For each $1\leq i\leq r$, write $A_i = [x_i,y_i]$.  We define the \emph{bricks} of $w_\eta$ to be the words
\[
u_i = h_{x_ik_{x_i}}\cdots h_{y_ik_{y_i}} \qquad\text{for each $1\leq i\leq r$.}
\]

\begin{lemma}\label{lem:brick}
Let $\eta\in\PEq_n$, and let $u_1,\ldots,u_r$ be the bricks of $w_\eta$, as above.  Then
\ben
\item \label{brick1} $w_\eta = u_1\cdots u_r$,
\item \label{brick2} $u_i \in \set{h_{st}}{x_i\leq s\leq t\leq y_i}^*$ for all $i$,
\item \label{brick3} $u_iu_j\sim u_ju_i$ for distinct $i,j$,
\item \label{brick4} $u_i\sim h_{x_iy_i}u_i$ for all $i$.
\een
\end{lemma}

\pf
\firstpfitem{\ref{brick1}}  This follows immediately from the definition.

\pfitem{\ref{brick2}}  We must show that every subscript of every letter of the brick ${u_i = h_{x_ik_{x_i}}\cdots h_{y_ik_{y_i}}}$ belongs to $A_i=[x_i,y_i]$.  This is obviously true of the \emph{first} subscripts.  Now consider $h_{jk_j}$, where $x_i\leq j\leq y_i$.  Since $j\in[x_i,y_i]=A_i$, we see that either:
\bit
\item $j'$ belongs to the transversal $A_i\cup B_i'$ of $d_\eta$, or else
\item $j'$ belongs to a lower-nontransversal $C_t'$ of $d_\eta$ that is nested by $B_i'$.
\eit
But $k_j = \suc_\eta(j)$ also belongs to $B_i$ or $C_t$ (respectively), and since $B_i$ and $C_t$ are both contained in $A_i=[x_i,y_i]$, we are done.

\pfitem{\ref{brick3}}  This follows immediately from part \ref{brick2} and relation \eqref{N2}.

\pfitem{\ref{brick4}}  This is a special case of Lemma \ref{lem:orbit}.
\epf

We are now almost ready to show that every word over $\XDn$ is $\sim$-equivalent to one of the form $w_\eta$.  The next lemma provides the key step in the proof of this fact.

\begin{lemma}\label{lem:aww}
If $\eta\in\PEq_n$ and $\oijn$, then $h_{ij}w_\eta \sim w_\mu$, where $\mu\in\PEq_n$ is as in Lemma \ref{lem:suc}.
\end{lemma}

\pf
Write $d_\eta = \begin{partn}{6} A_1&\cdots&A_r& \multicolumn{3}{c}{} \\ \hhline{~|~|~|-|-|-} B_1&\cdots&B_r&C_1&\cdots&C_s\end{partn}$, assuming that $A_1<\cdots<A_r$, and $i\in A_k$ and $j\in A_l$.  We also write $A_k = [x,p]$ and $A_l = [q,y]$, noting that
\[
p = \max(A_k) = \max[i]_{\wh\eta} \AND q = \min(A_l) = \min[j]_{\wh\eta}
\]
have the same meaning as in Lemma~\ref{lem:suc}.  Also let $u_1,\ldots,u_r$ be the bricks of $w_\eta$.

\pfcase1  Suppose first that $k=l$, and note that $\mu=\eta$ in this case, so we must show that $h_{ij}w_\eta\sim w_\eta$.  Using \eqref{N2} and Lemma \ref{lem:brick}\ref{brick2}, we have
\[
h_{ij}w_\eta = h_{ij}u_1\cdots u_r \sim u_1\cdots u_{k-1}\cdot h_{ij}u_k \cdot u_{k+1}\cdots u_r,
\]
so we can complete the proof in this case by showing that ${h_{ij}u_k\sim u_k}$.  But this follows from Lemma \ref{lem:brick}\ref{brick4} and \eqref{N1}:
\[
h_{ij}u_k \sim h_{ij}h_{xp}u_k \sim h_{xp}u_k \sim u_k.
\]

\pfcase2  Now suppose $k\not=l$, so that $k<l$.  Using Lemma \ref{lem:brick} and \eqref{N2}, we then have
\begin{equation}\label{eq:aww1}
h_{ij}w_\eta = h_{ij} u_1\cdots u_r \sim v\cdot h_{ij}u_ku_l \cdot v',
\end{equation}
where for simplicity we write $v = u_1\cdots u_{k-1}$ and $v' = u_{k+1}\cdots u_{l-1}u_{l+1}\cdots u_r$.  Next we claim that
\begin{equation}\label{eq:aww2}
h_{ij}u_ku_l \sim u_kh_{pq}u_l.
\end{equation}
Indeed, for this we have
\begin{align*}
h_{ij}u_ku_l &\sim h_{ij}\cdot h_{xp}u_k\cdot h_{qy}u_l &&\text{by Lemma \ref{lem:brick}\ref{brick4}}\\
&\sim h_{ij} h_{xp}h_{qy} \cdot u_ku_l &&\text{by \eqref{N2} and Lemma \ref{lem:brick}\ref{brick2}}\\
&\sim h_{pj} h_{xp}h_{qy} \cdot u_ku_l &&\text{by \eqref{N4}}\\
&\sim h_{pj} h_{qy}h_{xp} \cdot u_ku_l &&\text{by \eqref{N2}}\\
&\sim h_{pq} h_{qy}h_{xp} \cdot u_ku_l &&\text{by \eqref{N5}}\\
&\sim h_{pq}  \cdot h_{xp}u_k\cdot h_{qy}u_l &&\text{by \eqref{N2} and Lemma \ref{lem:brick}\ref{brick2}}\\
&\sim h_{pq}  \cdot u_ku_l &&\text{by Lemma \ref{lem:brick}\ref{brick4}}\\
&\sim u_kh_{pq}u_l &&\text{by \eqref{N2} and Lemma \ref{lem:brick}\ref{brick2}.}
\end{align*}
We then continue:
\begin{align*}
h_{ij}w_\eta &\sim v\cdot h_{ij}u_ku_l \cdot v' &&\text{by \eqref{eq:aww1}}\\
&\sim v\cdot u_kh_{pq}u_l \cdot v' &&\text{by \eqref{eq:aww2}}\\
&= u_1\cdots u_{k-1}\cdot u_kh_{pq}u_l \cdot u_{k+1}\cdots u_{l-1}u_{l+1}\cdots u_r \\
&\sim u_1\cdots u_{k-1}\cdot u_kh_{pq} \cdot u_{k+1}\cdots u_{l-1}u_lu_{l+1}\cdots u_r &&\text{by Lemma \ref{lem:brick}\ref{brick3}}\\
&= u_1\cdots u_{k-1}\cdot u_kh_{pq} \cdot u_{k+1}\cdots u_r,
\end{align*}
and we claim that this last word this is precisely $w_\mu$.
Now, $w_\eta = u_1\cdots u_{k-1}\cdot u_k \cdot u_{k+1}\cdots u_n$ and, writing $k_i = \suc_\eta(i)$ for each $1\leq i\leq n$, we note that:
\newpage
\bit
\item $u_k = h_{xk_x}\cdots h_{p-1,k_{p-1}}h_{pk_p} = h_{xk_x}\cdots h_{p-1,k_{p-1}}$, as $k_p=p$ (since $p = \max(A_k)=\max(B_k)$), and so
\item $u_kh_{pq} = h_{xk_x}\cdots h_{p-1,k_{p-1}}h_{pq}$.
\eit
It now follows from Lemma \ref{lem:suc}\ref{suc2} that indeed $u_1\cdots u_{k-1}\cdot u_kh_{pq} \cdot u_{k+1}\cdots u_r = w_\mu$.
\epf

\begin{cor}\label{cor:wweta}
For any $w\in \XDn^*$, we have $w\sim w_\eta$, where $\eta = \coker(\ol w)$.
\end{cor}

\pf
It suffices to show that $w\sim w_\eta$ for some $\eta\in\PEq_n$, since it will then follow from Lemma~\ref{lem:weta} that
\[
\eta = \coker(d_\eta) = \coker(\ol w_\eta) = \coker(\ol w).
\]
We do this by induction on $k$, the length of the word $w$.  If $k=0$, then $w=\io = w_{\De_\bn}$, so now suppose $k\geq1$.  We then have $w = h_{ij}w'$ for some $\oijn$, and some $w'\in\XDn^*$ of length $k-1$.  By induction, we have $w'\sim w_\eta$ for some $\eta\in\PEq_n$, and it then follows from Lemma \ref{lem:aww} that 
\[
w = h_{ij}w' \sim h_{ij}w_\eta \sim w_\mu,
\]
for some $\mu\in\PEq_n$.
\epf

\pf[\bf Proof of Theorem \ref{thm:Dn}.]
We have already observed that $\de$ is surjective, and that ${\sim}\sub\ker(\de)$.  For the reverse inclusion, suppose $(u,v)\in\ker(\de)$, so that $u,v\in\XDn^*$ and $\ol u=\ol v$.  By Corollary~\ref{cor:wweta}, we have $u\sim w_\eta$, where $\eta = \coker(\ol u)$.  But since $\ol u=\ol v$, we also have $v\sim w_\eta$ by the same corollary, so $u\sim w_\eta\sim v$.
\epf

\subsection[The action pair $(\D_n,\O_n)$]{\boldmath The action pair $(\D_n,\O_n)$}\label{subsect:DnOn}

We now bring the monoid $\O_n$ back into consideration, aiming to establish that $(\D_n,\O_n)$ is a strong (right) action pair in $\PPnfd$.  We begin with a simple lemma.  The proof uses the $R$ operation from \eqref{eq:DR}, even though $\PPnfd$ is not closed under it.

\begin{lemma}\label{lem:fdeta}
For any $f\in\O_n$ and $\eta\in\PEq_n$, we have $\coker(fd_\eta) = \eta$.
\end{lemma}

\pf
Since $f\in\T_n$ we have $R(f) = 1$, so it follows from Lemma \ref{lem:Rab=Rb} that $R(fd_\eta) = R(d_\eta)$.  But then
\[
\id_{\coker(fd_\eta)} = R(fd_\eta) = R(d_\eta) = \id_{\coker(d_\eta)} = \id_\eta,
\]
which implies the claim.
\epf

\begin{lemma}\label{lem:fdmu}
If $a\in\PPnfd$, and if $\eta=\coker(a)$, then $a=fd_\eta$ for some $f\in\O_n$.
\end{lemma}

\begin{figure}[ht]
\begin{center}
\begin{tikzpicture}[scale=.5]
\begin{scope}[shift={(0,1)}]	
\uvs{1,...,20}
\lvs{1,...,20}
\uuline16
\uuline7{20}
\stline13
\stline64
\stline78
\stline{20}{20}
\ddline12
\ddline34
\ddline{10}{12}
\ddline{13}{14}
\ddline{15}{16}
\darc57
\darcx8{13}{.6}
\darcx{14}{20}{.6}
\darcx{17}{19}{.3}
\draw(0.6,1)node[left]{$a=$};
\end{scope}
\begin{scope}[shift={(0,-4)}]	
\uvs{1,...,20}
\lvs{1,...,20}
\uuline16
\uuline7{20}
\stline13
\stline63
\stline7{12}
\stline{20}{12}
\draw(0.6,0)node[left]{$=$};
\draw[|-] (21,2)--(21,0);
\draw[right](21,1)node {$f$};
\end{scope}
\begin{scope}[shift={(0,-6)}]	
\uvs{1,...,20}
\lvs{1,...,20}
\uuline12
\uuline34
\uuline57
\uuline8{20}
\ddline12
\ddline34
\stline11
\stline22
\stline33
\stline44
\stline55
\stline77
\stline88
\stline{20}{20}
\ddline{10}{12}
\ddline{13}{14}
\ddline{15}{16}
\darc57
\darcx8{13}{.6}
\darcx{14}{20}{.6}
\darcx{17}{19}{.3}
\draw[|-|] (21,2)--(21,0);
\draw[right](21,1)node {$d_\eta$};
\end{scope}
\end{tikzpicture}
\caption{A planar full-domain partition $a\in\PP_{20}^\fd$, and a factorisation $a = fd_\eta$, where $f\in\O_{20}$ and $\eta = \coker(a)$, as in the proof of Lemma \ref{lem:fdmu}.}
\label{fig:PP20fd}
\end{center}
\end{figure}
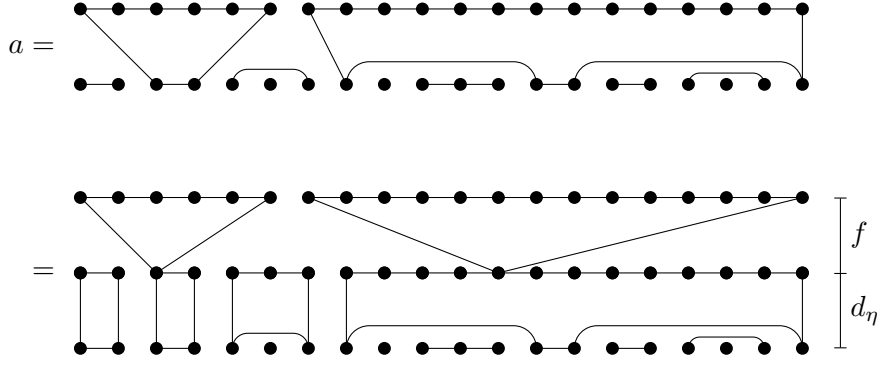

\pf
Write $a = \begin{partn}{6} A_1&\cdots&A_r& \multicolumn{3}{c}{} \\ \hhline{~|~|~|-|-|-} B_1&\cdots&B_r&C_1&\cdots&C_t\end{partn}$, where $A_1<\cdots<A_r$ and $B_1<\cdots<B_r$, and note that each $B_i$ is an un-nested $\eta$-class (but some of the $C_i$ might be un-nested as well).  In particular,~$d_\eta$ contains a transversal of the form $D_i\cup B_i'$ for each $1\leq i\leq r$, and by construction we have $D_1<\cdots<D_r$.  Thus, choosing $x_i\in D_i$ for each $i$, it follows that $f = \begin{partn}{3} A_1&\cdots&A_r \\ x_1&\cdots&x_r \end{partn}$ belongs to $\O_n$, and it is easy to see that~${a = fd_\eta\in\O_n\D_n}$.  (See Figure \ref{fig:PP20fd} for an example, with $n=20$.)
\epf

In the next statement, recall that we write $\ol h_{ii}=1$ for all $1\leq i\leq n$.

\begin{lemma}\label{lem:hfhg}
If $\oijn$ and $1\leq k\leq n-1$, then
\[
\ol h_{ij}\ol f_k = \begin{cases}
\ol f_k\ol h_{i-1,j} &\text{if $k=i-1$,}\\
\ol f_k\ol h_{i,j-1} &\text{if $k=j-1$,}\\
\ol f_k\ol h_{ij} &\text{otherwise,}
\end{cases}
\AND
\ol h_{ij}\ol g_k = \begin{cases}
\ol g_k\ol h_{i+1,j} &\text{if $k=i$,}\\
\ol g_k\ol h_{i,j+1} &\text{if $k=j$,}\\
\ol g_k\ol h_{ij} &\text{otherwise.}
\end{cases}
\]
\end{lemma}

\pf
This can be checked directly by drawing diagrams, similar to Figure \ref{fig:N3}.
\epf

The expressions in Lemma \ref{lem:hfhg} are not the simplest we could give, but will be useful for proving the next result.  For example, we actually have $\ol h_{ij}\ol f_k = \ol h_{i-1,j}$ when $k=i-1$.

\begin{prop}\label{prop:DnOn}\leavevmode
\ben
\item \label{DnOn1} $(\D_n,\O_n)$ is a strong (right) action pair in $\PPnfd$.
\item \label{DnOn2} We have the product decomposition $\PPnfd = \O_n\D_n$.
\een
\end{prop}

\pf
\firstpfitem{\ref{DnOn1}}  For \ref{A1} we need to show that for all $f\in\O_n$ and $u\in\D_n$, we have $uf=fv$ for some $v\in\D_n$.  This can be shown directly, though the argument is a little involved.  But it also follows from Lemma \ref{lem:hfhg}, and the fact that:
\bit
\item $\D_n=\pres{\ol h_{ij}}{\oijn}$, as shown in Proposition \ref{prop:Dn}, and
\item $\O_n=\pres{\ol f_k,\ol g_k}{1\leq k<n}$, as is well known; see for example \cite{Aizenstat1962,Solomon1996}.
\eit
We write $u$ as a product of $\ol h_{ij}$-s and $f$ as a product of $\ol f_k$-s and $\ol g_k$-s, and then use Lemma \ref{lem:hfhg} to move all $\ol f_k$-s and $\ol g_k$-s to the left.

For \ref{A2}, let $f,g\in\O_n$ and $\eta,\mu\in\PEq_n$.  Then using Lemma \ref{lem:fdeta}, we see that
\[
fd_\eta = gd_\mu \Implies \eta = \coker(fd_\eta) = \coker(gd_\mu) = \mu \Implies d_\eta = d_\mu.
\]

\pfitem{\ref{DnOn2}}  The forward inclusion follows from Lemma \ref{lem:fdmu}.  The reverse is clear.
\epf

\subsection[Presentation for $\PPnfd$]{\boldmath Presentation for $\PPnfd$}\label{subsect:presPPnfd}

We now come to the proof of Theorem \ref{thm:PPnfd}, which will be in two stages.  First we apply Theorem~\ref{thm:SU1}, using the action pair $(\D_n,\O_n)$, to obtain an initial presentation for $\PPnfd$.  This will involve the alphabet
\[
\set{f_i,g_i}{1\leq i<n} \cup \set{h_{ij}}{\oijn}.
\]
Then, using an argument similar to (but slightly easier than) that of Section \ref{sect:Pnfd}, we will deduce from this the presentation $\Mon\pres ZO$ stated in Theorem \ref{thm:PPnfd}.

We begin with some lemmas that will enable the application of Theorem \ref{thm:SU1}.  The first involves the left congruences $\th_u$ ($u\in\D_n$).  For a convex partition $\eta$ of $\bn$, with classes $A_1<\cdots<A_r$, let
\[
f_\eta = \begin{partn}{3} A_1&\cdots&A_r \\ x_1&\cdots&x_r \end{partn} \in \O_n, \WHERE x_i = \min(A_i) \text{ for each $1\leq i\leq k$.}
\]

\begin{lemma}\label{lem:thuDn}
For any $u\in\D_n$ we have $\th_u = (1,f_\eta)\lc$, where $\eta = \ker(u)$.
\end{lemma}

\pf
We must have $u=d_\mu$ for some $\mu\in\PEq_n$ with $\eta = \wh\mu$, and it is easy to check that $d_\mu f_\eta = f_\eta$ and $f_\eta d_\mu = d_\mu$.  The result now follows from Lemma~\ref{lem:1s}.
\epf

The next lemma will be used to verify assumption \ref{SU11} in Theorem \ref{thm:SU1}.

\begin{lemma}\label{lem:VDn}
Let $V = \set{\ol h_{i,i+1}}{1\leq i<n}$.  Then for any $u\in\D_n$ we have
\[
\th_u = \bigvee_{v\in V,\atop u\leqR \hspace{0.5mm} v} \th_v.
\]
\end{lemma}

\pf
Fix some $u\in\D_n$, write $\eta = \ker(u)$, and let the $\eta$-classes be $A_1<\cdots<A_r$.  Also let
\[
I = \bn\sm\{\max(A_1),\ldots,\max(A_r)\}
\AND
W = \set{\ol h_{i,i+1}}{i\in I} .
\]
We claim that $W$ consists of precisely the elements of $V$ that are $\geqR$-above $u$.  Indeed, since each $\ol h_{i,i+1}$ ($1\leq i<n$) is an idempotent, we have $u\leqR\ol h_{i,i+1} \iff u = \ol h_{i,i+1}u$, and it is easy to see that this is equivalent to having $i\in I$.  By Lemma \ref{lem:thuDn} we have
\[
\th_u = (1,f_\eta)\lc \AND \th_{\ol h_{i,i+1}} = (1,\ol f_i)\lc.
\]
If $I = \{i_1<\cdots<i_{n-r}\}$, then $f_\eta = \ol f_{i_{n-r}}\cdots\ol f_{i_1} \in \pres{\ol f_i}{i\in I}$, and it follows that assumptions \ref{bigvee1}--\ref{bigvee4} of Lemma \ref{lem:bigvee} are satisfied, and hence $\th_u = \bigvee_{v\in W}\th_v$, as required.
\epf

We are now ready to apply Theorem \ref{thm:SU1}.  We begin with:
\bit
\item the presentation $\Mon\pres\XDn\RDn$ for $\D_n$ from Theorem \ref{thm:Dn},
\item the presentation $\Mon\pres\XOn\ROn$ for $\O_n$, where $\XOn = \set{f_i,g_i}{1\leq i<n}$, and where~$\ROn$ is the set of all relations from $O$ involving only words over $\XOn$, namely \eqref{O2}, \eqref{O3}, \eqref{O5}--\eqref{O7}, \eqref{O10}, \eqref{O11}, and the relevant parts of \eqref{O1}.  This presentation for $\O_n$ is stated in \cite[Section~3]{Lavers1997}, where it was attributed to A\u{\i}zen\v{s}tat \cite{Aizenstat1962}; a more recent proof (in a more general context) can be found in \cite{Solomon1996}.
\eit

\begin{prop}\label{prop:PPnfd}
The monoid $\PPnfd$ has presentation
\[
\Mon\pres{\XOn\cup\XDn}{\ROn\cup\RDn\cup R_1\cup R_2(V)}
\]
via
\[
\Xi:(\XOn\cup\XDn)^*\to\PPnfd:f_i\mt\ol f_i,\ g_i\mt\ol g_i,\ h_{ij}\mt\ol h_{ij},
\]
where $R_1$ consists of relations
\begin{equation}\label{eq:R1PP}
h_{ij}f_k = \begin{cases}
f_kh_{i-1,j} &\text{if $k=i-1$,}\\
f_kh_{i,j-1} &\text{if $k=j-1$,}\\
f_kh_{ij} &\text{otherwise,}
\end{cases}
\AND
h_{ij}g_k = \begin{cases}
g_kh_{i+1,j} &\text{if $k=i$,}\\
g_kh_{i,j+1} &\text{if $k=j$,}\\
g_kh_{ij} &\text{otherwise,}
\end{cases}
\end{equation}
and $R_2(V)$ consists of relations
\begin{equation}\label{eq:R2VPP}
f_ih_{i,i+1} = h_{i,i+1} \qquad\text{for all $1\leq i<n$.}
\end{equation}
\end{prop}

\pf
This follows from Theorem \ref{thm:SU1}, using Lemma \ref{lem:VDn} to verify assumption \ref{SU11} of the theorem, and Lemmas \ref{lem:hfhg}, \ref{lem:thuDn} and \ref{lem:VDn} to verify the forms of the relations $R_1\cup R_2(V)$.
\epf

We now embark on the second stage of our proof of Theorem \ref{thm:PPnfd}, in which we will be aided by Proposition \ref{prop:PPnfd}.  A key role in all that follows will be played by the words
\[
\al_{ij} = h_ig_{i+1}\cdots g_{j-1} \in Z^* \qquad\text{for $\oijn$}.
\]

\begin{lemma}\label{lem:alh}
For any $\oijn$ we have $\al_{ij}\xi = \ol h_{ij} (= h_{ij}\Xi)$.  Consequently, $\xi$ is surjective.
\end{lemma}

\pf
The first assertion is easily checked diagrammatically.  The second follows because $\PPnfd$ is generated by $\set{\ol f_i,\ol g_i}{1\leq i<n}\cup\set{\ol h_{ij}}{\oijn}$, as follows from Proposition \ref{prop:PPnfd}.
\epf

So we now need to show that $\ker(\xi) = O^\sharp$.  As ever, the inclusion $O^\sharp\sub\ker(\xi)$ follows by simply checking diagrammatically that $\xi$ preserves each relation \eqref{O1}--\eqref{O12} from $O$.  The rest of this section is devoted to establishing the inclusion $\ker(\xi)\sub O^\sharp$.  To this end, we write ${\approx} = O^\sharp$.  First we make the following observation, which follows from the fact that $O$ contains the defining relations from the presentation $\Mon\pres{\XOn}{\ROn}$ for $\O_n$.

\begin{lemma}\label{lem:On}
If $u,v\in\XOn^*$ and $u\xi=v\xi$, then $u\approx v$.  
\end{lemma}

We now prove two lemmas involving the words $\al_{ij}$.  It will also be convenient to have an alternative expression for these words.  Specifically, we define
\[
\be_{ij} = h_{j-1}f_{j-2}\cdots f_i \qquad\text{for $\oijn$.}
\]

\begin{lemma}\label{lem:albe}
For any $\oijn$ we have $\al_{ij} \approx \be_{ij}$.
\end{lemma}

\pf
We fix $1\leq i<n$, and use induction on $j\in\{i+1,\ldots,n\}$.  When $j=i+1$, we have $\al_{i,i+1} = h_i = \be_{i,i+1}$.  For $i+2\leq j\leq n$, we use the inductive hypothesis, the definition of the $\al$ and $\be$ words, and the stated relations from $O$, to calculate:
\begin{align*}
\al_{ij} = \al_{i,j-1}\cdot g_{j-1} 
\approx \be_{i,j-1}\cdot g_{j-1} 
=h_{j-2}f_{j-3}\cdots f_i\cdot g_{j-1} 
&\approx_{\eqref{O7}} h_{j-2}g_{j-1}\cdot f_{j-3}\cdots f_i \\
&\approx_{\eqref{O12}} h_{j-1}f_{j-2}\cdot f_{j-3}\cdots f_i 
= \be_{ij}. \ \ \qedhere
\end{align*}
\epf

\begin{lemma}\label{lem:RDn}
For $i,j,k,l\in\bn$ we have
\begin{align}
\label{N1''} \al_{ij}\al_{kl} &\approx \al_{kl} &&\text{if $k\leq i<j\leq l$,}\\
\label{N2''} \al_{ij}\al_{kl} &\approx \al_{kl}\al_{ij} &&\text{if $i<j\leq k<l$,}\\
\label{N3''} \al_{ij}\al_{jk} \approx \al_{ik}\al_{ij} &\approx \al_{ik}\al_{jk} &&\text{if $i<j<k$.}
\end{align}
\end{lemma}

\pf
\firstpfitem{\eqref{N1''}}  First note that by the definition of the $\al$ and $\be$ words, and Lemma \ref{lem:albe}, we have
\[
\al_{kl} = \al_{kj}\cdot g_j\cdots g_{l-1} 
\approx \be_{kj}\cdot g_j\cdots g_{l-1} 
= \be_{ij}\cdot f_{i-1}\cdots f_k\cdot g_j\cdots g_{l-1} .
\]
Thus, it suffices to show that $\al_{ij}\be_{ij} \approx \be_{ij}$.  For $j=i+1$ we have $\al_{i,i+1}\be_{i,i+1} = h_ih_i \approx h_i$ by~\eqref{O1}.  For $j\geq i+2$ we have
\begin{align*}
\al_{ij}\be_{ij} &= \al_{i,j-1}g_{j-1}\cdot h_{j-1}f_{j-2}\cdots f_i \\
&\approx \be_{i,j-1}\cdot h_{j-1}f_{j-2}\cdots f_i &&\text{by Lemma \ref{lem:albe} and \eqref{O1}}\\
&= h_{j-2}\cdot f_{j-3}\cdots f_i\cdot h_{j-1}f_{j-2}\cdots f_i \\
&\approx h_{j-1}h_{j-2}\cdot f_{j-3}\cdots f_i\cdot f_{j-2}\cdots f_i &&\text{by \eqref{O4} and \eqref{O8}}\\
&\approx h_{j-1}h_{j-2}\cdot f_{j-2}\cdots f_i &&\text{by Lemma \ref{lem:On}}\\
&\approx h_{j-1}\cdot f_{j-2}\cdots f_i &&\text{by \eqref{O1}}\\
&= \be_{ij}.
\end{align*}

\pfitem{\eqref{N2''}}  By Lemma \ref{lem:albe}, it suffices to show that $\al_{ij}\be_{kl} \approx \be_{kl}\al_{ij}$ (for $i<j\leq k<l$).  Now,
\[
\al_{ij} = h_ig_{i+1}\cdots g_{j-1} \AND \be_{kl} = h_{l-1}f_{l-2}\cdots f_k,
\]
and relations \eqref{O4} and \eqref{O7}--\eqref{O9} say that each letter appearing in $\al_{ij}$ commutes with each letter appearing in $\be_{kl}$.  To verify this, it is best to consider separate cases in which $l\geq k+2$ and $l=k+1$.  In the latter we have $\be_{kl} = h_{l-1} = h_k$.

\pfitem{\eqref{N3''}}  First note that
\[
\al_{ik}\al_{jk} = \al_{ij}g_j\cdots g_{k-1}\al_{jk},
\]
so we can show that $\al_{ik}\al_{jk} \approx \al_{ij}\al_{jk}$ by showing that $g_j\cdots g_{k-1}\al_{jk}\approx\al_{jk}$, and we do this by induction on $k$.  When $k=j+1$ we have 
\[
g_j\cdots g_{k-1}\al_{jk} = g_j\cdot h_j \approx_{\eqref{O1}} h_j = \al_{jk}.
\]
For $k\geq j+2$ we have
\begin{align*}
g_j\cdots g_{k-2}g_{k-1}\al_{jk} &\approx g_j\cdots g_{k-2}g_{k-1}\be_{jk} &&\text{by Lemma \ref{lem:albe}}\\
&= g_j\cdots g_{k-2}g_{k-1}\cdot h_{k-1}f_{k-2}\cdots f_j \\
&\approx g_j\cdots g_{k-2}\cdot h_{k-1}f_{k-2}\cdots f_j &&\text{by \eqref{O1}}\\
&= g_j\cdots g_{k-2}\cdot \be_{jk} \\
&\approx g_j\cdots g_{k-2}\cdot \al_{jk} &&\text{by Lemma \ref{lem:albe}}\\
&= g_j\cdots g_{k-2}\cdot \al_{j,k-1}g_{k-1} \\
&\approx \al_{j,k-1}g_{k-1} &&\text{by induction}\\
&= \al_{jk}.
\end{align*}
It remains to show that $\al_{ik}\al_{ij} \approx \al_{ij}\al_{jk}$, and since $\al_{ij}\al_{jk}\approx\al_{jk}\al_{ij}$ by \eqref{N2''}, it suffices to show that $\al_{ik}\al_{ij} \approx \al_{jk}\al_{ij}$, which by Lemma \ref{lem:albe} is equivalent to $\be_{ik}\be_{ij} \approx \be_{jk}\be_{ij}$.  But this is done in similar fashion to the proof of $\al_{ik}\al_{jk} \approx \al_{ij}\al_{jk}$ above.  We first note that
\[
\be_{ik}\be_{ij} = \be_{jk}f_{j-1}\cdots f_i\be_{ij},
\]
and then show by induction that $f_{j-1}\cdots f_i\be_{ij} \approx \be_{ij}$.
\epf

\begin{lemma}\label{lem:R1}
For $\oijn$ and $1\leq k<n$ we have
\[
\al_{ij}f_k \approx \begin{cases}
f_k\al_{i-1,j} &\text{if $k=i-1$,}\\
f_k\al_{i,j-1} &\text{if $k=j-1$,}\\
f_k\al_{ij} &\text{otherwise,}
\end{cases}
\AND
\al_{ij}g_k \approx \begin{cases}
g_k\al_{i+1,j} &\text{if $k=i$,}\\
g_k\al_{i,j+1} &\text{if $k=j$,}\\
g_k\al_{ij} &\text{otherwise.}
\end{cases}
\]
\end{lemma}

\pf
We just prove the relations involving $f_k$, as those for $g_k$ are analogous.

For $k=i-1$ we first use Lemma \ref{lem:albe} to calculate 
\[
\al_{ij}f_{i-1} \approx \be_{ij}f_{i-1} = \be_{i-1,j} \approx \al_{i-1,j}.
\]
Since $\al_{i-1,j} = h_{i-1}g_i\cdots g_{j-1}$, it follows from \eqref{O1} that $\al_{i-1,j}\approx f_{i-1}\al_{i-1,j}$.

For $k=j-1$ we first deal with the $j=i+1$ case.  Here by \eqref{O1} we have
\[
\al_{ij}f_{j-1} = h_if_i \approx f_i = f_ih_{ii} = f_{j-1}h_{i,j-1}.
\]
On the other hand, if $j\geq i+2$ then
\begin{align*}
\al_{ij}f_{j-1} &= h_ig_{i+1}\cdots g_{j-2}g_{j-1}\cdot f_{j-1} \\
&\approx h_ig_{i+1}\cdots g_{j-2}\cdot f_{j-1} &&\text{by \eqref{O1}}\\
&\approx f_{j-1}\cdot h_ig_{i+1}\cdots g_{j-2} &&\text{by \eqref{O7} and \eqref{O8}}\\
&= f_{j-1}\al_{i,j-1}.
\end{align*}

When $k\not=i-1$ and $k\not=j-1$, one of the following holds:
\[
k\leq i-2 \COMMA i\leq k\leq j-2 \qquad\text{or}\qquad k\geq j.
\]
In the first and third cases, \eqref{O7} and \eqref{O8} say that $f_k$ commutes with each letter appearing in $\al_{ij} = h_ig_{i+1}\cdots g_{j-1}$.  This leaves us to deal with the case that $i\leq k\leq j-2$.  Here we first use Lemma~\ref{lem:On} to calculate
\[
\al_{ij}f_k = h_i\cdot g_{i+1}\cdots g_{j-1}f_k \approx h_i\cdot g_{i+1}\cdots g_{j-1} = \al_{ij}.
\]
Thus, it remains to show that $\al_{ij}\approx f_k\al_{ij}$.  This follows from \eqref{O1} when $k=i$, as the first letter of $\al_{ij}$ is $h_i$.  When $i+1\leq k\leq j-2$ we use \eqref{O8} and then Lemma \ref{lem:On} to calculate
\[
f_k\al_{ij} = f_kh_ig_{i+1}\cdots g_{j-1} \approx h_if_kg_{i+1}\cdots g_{j-1} \approx h_ig_{i+1}\cdots g_{j-1} = \al_{ij}.  \qedhere
\]
\epf

We now define a morphism
\[
(\XOn\cup\XDn)^*\to Z^*:w\mt\wh w \BY \wh f_i = f_i \COMMa \wh g_i = g_i \ANd \wh h_{ij} = \al_{ij}.
\]

\begin{lemma}\label{lem:kerXi}
For any $(u,v)\in\ker(\Xi)$, we have $\wh u\approx\wh v$.
\end{lemma}

\pf
Since $\ker(\Xi) = (\ROn\cup\RDn\cup R_1\cup R_2(V))^\sharp$ by Proposition \ref{prop:PPnfd}, we can assume that $(u,v)$ is one of the relations from this set.  We consider the relations of the various forms separately.
\bit
\item  If $(u,v)\in\ROn$, then $\wh u= u \approx v = \wh v$ since $O$ contains~$\ROn$.
\item  If $(u,v)\in\RDn$, then Lemma \ref{lem:RDn} gives $\wh u\approx\wh v$.
\item  Similarly, Lemma \ref{lem:R1} deals with relations from $R_1$.
\item  Finally, if $(u,v)\in R_2(V)$, then $\wh u = f_ih_i \approx h_i = \wh v$ by \eqref{O1}.  \qedhere
\eit
\epf

We are now ready to complete our objective.

\pf[\bf Proof of Theorem \ref{thm:PPnfd}.]
As already observed, it remains to show that $\ker(\xi)\sub{\approx}$, so let $(u,v)\in\ker(\xi)$.  Let $u'$ and $v'$ be the words over $\XOn\cup\XDn$ obtained by replacing any letter $h_i\in Z$ by $h_{i,i+1}\in\XDn$.  By Lemma \ref{lem:alh} we have $u'\Xi = u\xi$ and $v'\Xi = v\xi$.  Since $(u,v)\in\ker(\xi)$, it follows that $(u',v')\in\ker(\Xi)$.  Lemma \ref{lem:kerXi} then gives $\wh u\sgap'\approx\wh v\sgap'$.  But $\wh u\sgap' = u$ and $\wh v\sgap' = v$, and this completes the proof that $u\approx v$.
\epf

\subsection[The structure of $\PPnfd$]{\boldmath The structure of $\PPnfd$}\label{subsect:structurePPnfd}

Our study of the full-domain partition monoid $\Pnfd$ made heavy use of the $D$ and $R$ operations given in \eqref{eq:DR}, which give $\Pnfd$ the structure of a right restriction Ehresmann monoid.  We have already observed that the planar submonoid $\PPnfd$ is closed under $D$ but not $R$, so it does not inherit the property of being a right restriction monoid (or even an Ehresmann monoid).  It does however satisfy some more general conditions with respect to a different unary range operation, which we will denote by $\rho$ to distinguish it from $R$.  Specifically, we define
\[
\rho(a) = d_{\coker(a)} \qquad\text{for $a\in\PPnfd$.}
\]
So $\rho$ is a map $\PPnfd\to\D_n$.  The next result shows that $\PPnfd$ is a so-called \emph{grrac monoid}.  This is the left-right dual of the \emph{glrac monoids} introduced in \cite{BGG2010} and studied further in \cite{Jones2023}.  The acronyms stand for `\emph{g}eneralised \emph{r}ight/\emph{l}eft \emph{r}estriction monoids possessing the \emph{a}mple and \emph{c}ongruence conditions', which are \ref{G8} and \ref{G7}, respectively, in the following list.  In the proof we again make use of the $R$ operation, even though $\PPnfd$ is not closed under it.

\newpage

\begin{prop}\label{prop:grrac}
For all $a,b\in\PPnfd$ we have:
\begin{enumerate}[label=\textup{\textsf{(G\arabic*)}},leftmargin=9mm]\bmc2
\item \label{G1} $a\rho(a) = a$,
\item \label{G2} $\rho(\rho(a)) = \rho(a)$,
\item \label{G5} $\rho(a)\rho(a) = \rho(a)$,
\item \label{G4} $\rho(a)\rho(b)\rho(a) = \rho(b)\rho(a)$,
\item \label{G3} $\rho(\rho(a)\rho(b)) = \rho(a)\rho(b)$,
\item \label{G6} $\rho(ab)\rho(b) = \rho(ab)$,
\item \label{G7} $\rho(ab) = \rho(\rho(a)b)$,
\item \label{G8} $\rho(a)b = b\rho(ab)$.
\emc
\end{enumerate}
\end{prop}

\pf
We begin with some simple observations that be used throughout the proof; these hold for any $a,b\in\PPnfd$.  First, by construction, we have
\begin{equation}\label{eq:cokerrho}
\coker(\rho(a)) = \coker(d_{\coker(a)}) = \coker(a).
\end{equation}
It follows from this that
\begin{equation}\label{eq:Rrho}
R(\rho(a)) = \id_{\coker(\rho(a))} = \id_{\coker(a)} = R(a).
\end{equation}
Combining this with \ref{E3}, it follows that 
\begin{equation}\label{eq:Rrhoab}
R(ab) = R(R(a)b) = R(R(\rho(a))b) = R(\rho(a)b).
\end{equation}

\pfitem{\ref{G1}}  By Lemma \ref{lem:fdmu}, we have $a=fd_\eta$ for some $f\in\O_n $, and where $\eta=\coker(a)$.  Since $d_\eta$ is an idempotent, it follows that $a=ad_\eta = a\rho(a)$.  

\pfitem{\ref{G2}}  Using \eqref{eq:cokerrho}, we have $\rho(\rho(a)) = d_{\coker(\rho(a))} = d_{\coker(a)} = \rho(a)$.

\pfitem{\ref{G5} and \ref{G4}}  These follow from the fact that $\D_n$ is a right regular band (Proposition \ref{prop:Dn}).

\pfitem{\ref{G3}}  Since $\D_n$ is a monoid, we have $\rho(a)\rho(b) = d_\eta$ for some $\eta\in\PEq_n$, so this identity follows from $\rho(d_\eta) = d_{\coker(d_\eta)} = d_\eta$.

\pfitem{\ref{G6}}  Since $\rho(ab)\rho(b)$ and $\rho(ab)$ both belong to $\D_n$ we can show they are equal by showing that they have the same cokernel, and by definition this is equivalent to showing that
\[
R(\rho(ab)\rho(b)) = R(\rho(ab)).
\]
For this we have
\[
R(\rho(ab)\rho(b)) =_{\eqref{eq:Rrhoab}} R(ab\rho(b))
=_{\ref{G1}} R(ab)
=_{\eqref{eq:Rrho}} R(\rho(ab)).
\]

\pfitem{\ref{G7}}  Again, we can establish this by showing that $R(\rho(ab)) = R(\rho(\rho(a)b))$.  By \eqref{eq:Rrho}, this is equivalent to $R(ab) = R(\rho(a)b)$, and this is precisely \eqref{eq:Rrhoab}.

\pfitem{\ref{G8}}  Let $u = \rho(a)$, and note that by \ref{G7}, the identity \ref{G8} is equivalent to $ub = b\rho(ub)$, i.e.~to
\[
ub = bd_\eta, \WHERE \eta = \coker(ub).
\]
Now, by Lemma \ref{lem:fdmu}, we have $b = fd_\mu$ for some $f\in\O_n$, where $\mu = \coker(b)$.  Since $(\D_n,\O_n)$ is a strong action pair (Proposition \ref{prop:DnOn}), it follows from property \ref{A1} that $uf=fv$ for some $v\in\D_n$.  Since $\D_n$ is a right regular band (Proposition \ref{prop:Dn}), we have $vd_\mu = d_\mu vd_\mu$.  Also, since~$\D_n$ is a monoid, we have $vd_\mu = d_\nu$ for some $\nu\in\PEq_n$.  Putting all of this together, it follows that
\begin{equation}\label{eq:ub}
ub = ufd_\mu = fvd_\mu = fd_\mu vd_\mu = bd_\nu.
\end{equation}
So it remains to show that $\nu = \eta$.  But for this we have
\begin{align*}
\nu &= \coker(fd_\nu) &&\text{by Lemma \ref{lem:fdeta}}\\
&= \coker(fvd_\mu) &&\text{by definition of $\nu$}\\
&= \coker(ub) &&\text{by \eqref{eq:ub}}\\
&= \eta.  && \qedhere
\end{align*}
\epf

\footnotesize
\def\bibspacing{-1.1pt}
\bibliography{biblio}

\begin{thebibliography}{10}

\bibitem{Auinger2012}
K.~Auinger.
\newblock Krohn-{R}hodes complexity of {B}rauer type semigroups.
\newblock {\em Port. Math.}, 69(4):341--360, 2012.

\bibitem{Auinger2014}
K.~Auinger.
\newblock Pseudovarieties generated by {B}rauer type monoids.
\newblock {\em Forum Math.}, 26(1):1--24, 2014.

\bibitem{ACHLV2015}
K.~Auinger, Y.~Chen, X.~Hu, Y.~Luo, and M.~V. Volkov.
\newblock The finite basis problem for {K}auffman monoids.
\newblock {\em Algebra Universalis}, 74(3-4):333--350, 2015.

\bibitem{ADV2012_2}
K.~Auinger, I.~Dolinka, and M.~V. Volkov.
\newblock Equational theories of semigroups with involution.
\newblock {\em J. Algebra}, 369:203--225, 2012.

\bibitem{Aizenstat1962}
A.~J. A\u{\i}zen\v{s}tat.
\newblock The defining relations of the endomorphism semigroup of a finite
  linearly ordered set.
\newblock {\em Sibirsk. Mat. \v Z.}, 3:161--169, 1962.

\bibitem{BH2014}
G.~Benkart and T.~Halverson.
\newblock Motzkin algebras.
\newblock {\em European J. Combin.}, 36:473--502, 2014.

\bibitem{BDP2002}
M.~Borisavljevi{\'c}, K.~Do{\v{s}}en, and Z.~Petri{\'c}.
\newblock Kauffman monoids.
\newblock {\em J. Knot Theory Ramifications}, 11(2):127--143, 2002.

\bibitem{BGG2010}
M.~J.~J. Branco, G.~M.~S. Gomes, and V.~Gould.
\newblock Extensions and covers for semigroups whose idempotents form a left
  regular band.
\newblock {\em Semigroup Forum}, 81(1):51--70, 2010.

\bibitem{Brauer1937}
R.~Brauer.
\newblock On algebras which are connected with the semisimple continuous
  groups.
\newblock {\em Ann. of Math. (2)}, 38(4):857--872, 1937.

\bibitem{CDEGZ2023}
S.~Carson, I.~Dolinka, J.~East, V.~Gould, and R.-e. Zenab.
\newblock Product decompositions of semigroups induced by action pairs.
\newblock {\em Dissertationes Math.}, 587:180, 2023.

\bibitem{CHKLV2019}
Y.~Chen, X.~Hu, N.~V. Kitov, Y.~Luo, and M.~V. Volkov.
\newblock Identities of the {K}auffman monoid {$\mathcal K_3$}.
\newblock {\em Comm. Algebra}, 48(5):1956--1968, 2020.

\bibitem{CPbook}
A.~H. Clifford and G.~B. Preston.
\newblock {\em The algebraic theory of semigroups. {V}ol. {I}}.
\newblock Mathematical Surveys, No. 7. American Mathematical Society,
  Providence, R.I., 1961.

\bibitem{JEgrpm}
J.~East.
\newblock Generators and relations for partition monoids and algebras.
\newblock {\em J. Algebra}, 339:1--26, 2011.

\bibitem{East2013}
J.~East.
\newblock Defining relations for idempotent generators in finite full
  transformation semigroups.
\newblock {\em Semigroup Forum}, 86(3):451--485, 2013.

\bibitem{JEinsn2}
J.~East.
\newblock A symmetrical presentation for the singular part of the symmetric
  inverse monoid.
\newblock {\em Algebra Universalis}, 74(3-4):207--228, 2015.

\bibitem{JEgrpm2}
J.~East.
\newblock Presentations for (singular) partition monoids: a new approach.
\newblock {\em Math. Proc. Cambridge Philos. Soc.}, 165(3):549--562, 2018.

\bibitem{JErook}
J.~East.
\newblock Presentations for rook partition monoids and algebras and their
  singular ideals.
\newblock {\em J. Pure Appl. Algebra}, 223(3):1097--1122, 2019.

\bibitem{East2021}
J.~East.
\newblock Presentations for {T}emperley-{L}ieb algebras.
\newblock {\em Q. J. Math.}, 72(4):1253--1269, 2021.

\bibitem{EGM2025}
J.~East, R.~Gray, and P.~A. Azeef~Muhammed.
\newblock Topology of {G}raham--{H}oughton complexes.
\newblock {\em In preparation}.

\bibitem{EG2017}
J.~East and R.~D. Gray.
\newblock Diagram monoids and {G}raham--{H}oughton graphs: {I}dempotents and
  generating sets of ideals.
\newblock {\em J. Combin. Theory Ser. A}, 146:63--128, 2017.

\bibitem{EG2021}
J.~East and R.~D. Gray.
\newblock Ehresmann theory and partition monoids.
\newblock {\em J. Algebra}, 579:318--352, 2021.

\bibitem{EMRT2018}
J.~East, J.~D. Mitchell, N.~Ru\v{s}kuc, and M.~Torpey.
\newblock Congruence lattices of finite diagram monoids.
\newblock {\em Adv. Math.}, 333:931--1003, 2018.

\bibitem{FitzGerald2003}
D.~G. FitzGerald.
\newblock A presentation for the monoid of uniform block permutations.
\newblock {\em Bull. Austral. Math. Soc.}, 68(2):317--324, 2003.

\bibitem{FL2011}
D.~G. FitzGerald and K.~W. Lau.
\newblock On the partition monoid and some related semigroups.
\newblock {\em Bull. Aust. Math. Soc.}, 83(2):273--288, 2011.

\bibitem{Gould_notes}
V.~Gould.
\newblock Notes on restriction semigroups and related structures.
\newblock \\ \url{https://www-users.york.ac.uk/~varg1/restriction.pdf}, 2010.

\bibitem{Howie1995}
J.~M. Howie.
\newblock {\em Fundamentals of semigroup theory}, volume~12 of {\em London
  Mathematical Society Monographs. New Series}.
\newblock The Clarendon Press, Oxford University Press, New York, 1995.
\newblock Oxford Science Publications.

\bibitem{Jones2023}
P.~R. Jones.
\newblock On the structure of glrac semigroups.
\newblock {\em Semigroup Forum}, 106(1):169--183, 2023.

\bibitem{Jones1994_2}
V.~F.~R. Jones.
\newblock The {P}otts model and the symmetric group.
\newblock In {\em Subfactors ({K}yuzeso, 1993)}, pages 259--267. World Sci.
  Publ., River Edge, NJ, 1994.

\bibitem{KV2019}
N.~V. Kitov and M.~V. Volkov.
\newblock Identities of the {K}auffman monoid {$\mathcal K_4$} and of the
  {J}ones monoid {$\mathcal J_4$}.
\newblock In {\em Fields of logic and computation. {III}}, volume 12180 of {\em
  Lecture Notes in Comput. Sci.}, pages 156--178. Springer, Cham, [2020.

\bibitem{KV2023}
N.~V. Kitov and M.~V. Volkov.
\newblock Identities in twisted {B}rauer monoids.
\newblock In {\em Semigroups, algebras and operator theory}, volume 436 of {\em
  Springer Proc. Math. Stat.}, pages 79--103. Springer, Singapore, 2023.

\bibitem{Koenig2008}
S.~Koenig.
\newblock A panorama of diagram algebras.
\newblock In {\em Trends in representation theory of algebras and related
  topics}, EMS Ser. Congr. Rep., pages 491--540. Eur. Math. Soc., Z\"urich,
  2008.

\bibitem{KM2006}
G.~Kudryavtseva and V.~Mazorchuk.
\newblock On presentations of {B}rauer-type monoids.
\newblock {\em Cent. Eur. J. Math.}, 4(3):413--434 (electronic), 2006.

\bibitem{LF2006}
K.~W. Lau and D.~G. FitzGerald.
\newblock Ideal structure of the {K}auffman and related monoids.
\newblock {\em Comm. Algebra}, 34(7):2617--2629, 2006.

\bibitem{Lavers1997}
T.~G. Lavers.
\newblock The theory of vines.
\newblock {\em Comm. Algebra}, 25(4):1257--1284, 1997.

\bibitem{Lawson1991}
M.~V. Lawson.
\newblock Semigroups and ordered categories. {I}. {T}he reduced case.
\newblock {\em J. Algebra}, 141(2):422--462, 1991.

\bibitem{Lawson2021}
M.~V. Lawson.
\newblock On {E}hresmann semigroups.
\newblock {\em Semigroup Forum}, 103(3):953--965, 2021.

\bibitem{MM2007}
V.~Maltcev and V.~Mazorchuk.
\newblock Presentation of the singular part of the {B}rauer monoid.
\newblock {\em Math. Bohem.}, 132(3):297--323, 2007.

\bibitem{MS2021}
S.~Margolis and I.~Stein.
\newblock Ehresmann semigroups whose categories are {EI} and their
  representation theory.
\newblock {\em J. Algebra}, 585:176--206, 2021.

\bibitem{Martin1994}
P.~Martin.
\newblock Temperley-{L}ieb algebras for nonplanar statistical mechanics---the
  partition algebra construction.
\newblock {\em J. Knot Theory Ramifications}, 3(1):51--82, 1994.

\bibitem{Martin2008}
P.~Martin.
\newblock On diagram categories, representation theory and statistical
  mechanics.
\newblock In {\em Noncommutative rings, group rings, diagram algebras and their
  applications}, volume 456 of {\em Contemp. Math.}, pages 99--136. Amer. Math.
  Soc., Providence, RI, 2008.

\bibitem{MM2014}
P.~Martin and V.~Mazorchuk.
\newblock On the representation theory of partial {B}rauer algebras.
\newblock {\em Q. J. Math.}, 65(1):225--247, 2014.

\bibitem{Maz2002}
V.~Mazorchuk.
\newblock Endomorphisms of {$\mathfrak B_n$, $\mathcal P \mathfrak B_n$}, and
  {$\mathfrak C_n$}.
\newblock {\em Comm. Algebra}, 30(7):3489--3513, 2002.

\bibitem{Mace4}
W.~McCune.
\newblock Prover9 and {M}ace4.
\newblock \url{http://www.cs.unm.edu/~mccune/prover9/}.

\bibitem{Solomon1996}
A.~Solomon.
\newblock Catalan monoids, monoids of local endomorphisms, and their
  presentations.
\newblock {\em Semigroup Forum}, 53(3):351--368, 1996.

\bibitem{Stein2017}
I.~Stein.
\newblock Algebras of {E}hresmann semigroups and categories.
\newblock {\em Semigroup Forum}, 95(3):509--526, 2017.

\bibitem{Stokes2015}
T.~Stokes.
\newblock Domain and range operations in semigroups and rings.
\newblock {\em Comm. Algebra}, 43(9):3979--4007, 2015.

\bibitem{Stokes2022b}
T.~Stokes.
\newblock Left restriction monoids from left {$E$}-completions.
\newblock {\em J. Algebra}, 608:143--185, 2022.

\bibitem{TL1971}
H.~N.~V. Temperley and E.~H. Lieb.
\newblock Relations between the ``percolation'' and ``colouring'' problem and
  other graph-theoretical problems associated with regular planar lattices:
  some exact results for the ``percolation'' problem.
\newblock {\em Proc. Roy. Soc. London Ser. A}, 322(1549):251--280, 1971.

\bibitem{Wilcox2007}
S.~Wilcox.
\newblock Cellularity of diagram algebras as twisted semigroup algebras.
\newblock {\em J. Algebra}, 309(1):10--31, 2007.

\end{thebibliography}
\bibliographystyle{abbrv}

\end{document}